\documentclass[11pt,a4paper]{article}

\usepackage[margin=1in]{geometry}
\usepackage{amsmath,amsthm,amssymb,mathtools}
\usepackage{mathrsfs}
\usepackage{enumitem}
\usepackage{array,tabularx}
\usepackage{microtype}
\usepackage[hidelinks]{hyperref}
\usepackage{xcolor}
\usepackage{tikz-cd}
\usepackage{graphicx}
\usepackage{needspace}

\newtheorem{thm}{Theorem}[section]
\newtheorem{pro}[thm]{Proposition}
\newtheorem{lem}[thm]{Lemma}
\newtheorem{Def}[thm]{Definition}
\newtheorem{cor}[thm]{Corollary}

\newtheorem{mainthmx}{Theorem}

\newtheorem{maintheoremi}{Theorem}

\newtheorem{mainproa}{Theorem}

\newtheorem{mainprob}{Theorem}

\newtheorem{mainproii}{Theorem}

\newtheorem{defpro}[thm]{Proposition--Definition}

\theoremstyle{definition}
\newtheorem{convention}[thm]{Convention}
\theoremstyle{remark}
\newtheorem{rem}[thm]{Remark}

\newcommand{\Z}{\mathbb Z}
\newcommand{\N}{\mathbb N}
\newcommand{\Ztwo}{\mathbb Z_2}
\newcommand{\Dinf}{D_\infty}
\newcommand{\Aalt}{\mathbf A}
\newcommand{\NekExpHyp}{\ensuremath{(\mathrm{N}_{\mathrm{exp}})}}

\newcommand{\Fix}{\operatorname{Fix}}
\newcommand{\Int}{\operatorname{Int}}
\newcommand{\Clop}{\operatorname{Clop}}

\newcommand{\MLD}{\mathrm{MLD}}

\newcommand{\diam}{\operatorname{diam}}

\newcommand{\Mov}{\operatorname{Mov}}

\newcommand{\Xc}{\mathcal X}

\title{Orbital graph structure of alternating full groups\\
and simple torsion amenable groups of exponential growth}
\author{Yitong Zhang\thanks{Laboratoire de Probabilit\'es, Statistique et
Mod\'elisation (LPSM), Sorbonne Universit\'e, 4 place Jussieu,
75005 Paris, France. Email: \href{mailto:zhang.yitong.math@gmail.com}{\texttt{zhang.yitong.math@gmail.com}}.}}
\date{September 2026}
\hypersetup{pdftitle={Orbital graph structure of alternating full groups and simple torsion amenable groups of exponential growth},pdfauthor={Yitong Zhang}}

\hypersetup{bookmarksdepth=2}
\begin{document}

\maketitle

\begin{abstract}
We develop a framework for constructing simple alternating full groups
$\Aalt(G,X)$, as defined in \cite{NekDyn}, associated with finite fragmentations of minimal expansive
dihedral Cantor actions, with controlled orbital-graph codings and purely
non-Hausdorff singularities.

We also prove that positive entropy of a naturally labeled admissible
orbital line forces exponential word growth of the acting group.

As an application, we construct continuum many pairwise abstractly
nonisomorphic, finitely generated, infinite, simple, torsion, amenable groups
of exponential word growth. They arise as alternating full groups of
fragmented minimal expansive dihedral Cantor actions with linearly repetitive
generic regular orbital graphs and purely non-Hausdorff singularities.

Under the same minimality, expansivity, linear repetitivity and
singularity assumptions, if every finite clopen coding is determined
by a bounded fragmentation window, the alternating full group has
intermediate word growth.
\end{abstract}

\tableofcontents

\clearpage
\section{Introduction}\label{sec:introduction}

Let $\Dinf=\langle a,b\mid a^2=b^2=1\rangle$ act on an infinite Cantor
set $X$, and let $G=\langle A\cup B\rangle\curvearrowright X$ be a finite
fragmentation of this action. We write \NekExpHyp\ when the dihedral action
is minimal and expansive, the fragmentation has a purely non-Hausdorff
(PNH) singularity, and the $(A\cup B)$-labeled orbital graphs at generic
regular points are linearly repetitive (LR). These are the assumptions
in the linearly repetitive, expansive case of
\cite[Theorem~1.2]{Nek2018}.

Let $T_2$ be the rooted binary tree.  For
$\omega\in\{0,1,2\}^{\N}$, write $G_\omega$ for the generalized
Grigorchuk group in its standard boundary action on $\partial T_2$;
$G_{(012)^\infty}$ is the classical first Grigorchuk group
\cite{Gri85}.  For every parameter $\omega$ used here, each of the
symbols $0,1,2$ occurs infinitely often.  In the dyadic Gray-code model of
Section~\ref{sec:preliminaries}, this action is generated by one
unfragmented reflection $a$ and a Klein four fragmentation of the other
reflection $b$, with the level-$k$ choice prescribed by $\omega_k$;
see also \cite[Example~3.3]{Nek2018}.  When $G_\omega$ acts on an
extension below, we use the lifts of these same generators.

Call $\omega\in\{0,1,2\}^{\N}$ an \emph{LR parameter} if all three symbols
occur with bounded gaps.  Among parameters in which all three symbols
occur infinitely often, this condition is equivalent to LR of the orbital
graphs of $G_\omega\curvearrowright\partial T_2$ at regular points
$z\in\Ztwo\setminus\Z$; see Proposition~\ref{pro:valuationLR}.

\subsection*{Main results}
\Needspace{8\baselineskip}
\paragraph{Orbital graph of the alternating full group of a fragmentation.}
Recall that a Cantor action of a finitely generated group $\Lambda$ is expansive
if and only if it admits a finite clopen generating partition.
For a minimal expansive dihedral action $\Dinf=\langle a,b\rangle\curvearrowright X$,
put $T=ab$.  By Lemma~\ref{lem:reflection-invariant-generator}, we may
choose a finite clopen generating partition $\mathcal Q_X$ whose atoms
are individually $b$-invariant.  Write $f_X:X\to\mathcal B_X$ for its symbol map and set
\[
 w_X(x)(n)=f_X(T^n x),\qquad n\in\Z.
\]
The map $w_X$ conjugates $T\curvearrowright X$ to a subshift in
$\mathcal B_X^\Z$, and $w_X(bx)(n)=w_X(x)(-n)$, so $b$ acts by ordinary
reversal.  The full dihedral coding $\chi_X=c_{f_X}$ is determined by
$\chi_X(x)(T^n)=\chi_X(x)(bT^n)=w_X(x)(n)$.
We give the coding conventions in Section~\ref{sec:preliminaries}.
For a bi-infinite sequence $u$, write
$\Omega(u)=\overline{\{\sigma^n u:n\in\Z\}}$ for its shift orbit closure.

\Needspace{10\baselineskip}
\begin{mainthmx}[Containing fragmentations and joint orbital codings]
\label{thm:main-orbital-coding}
Let $G=\langle A\cup B\rangle\curvearrowright X$ be a finite fragmentation
of a minimal expansive Cantor action
$\Dinf=\langle a,b\rangle\curvearrowright X$, and put $L=\Aalt(G,X)$.
Fix a $b$-invariant generating coding $w_X$ as above.  There is a further finite
fragmentation $H$ of the same dihedral action, with the same germ groupoid
as $G$, such that
\[
 G\leq H\leq\mathbf F(G\curvearrowright X),
 \qquad L\leq H,
 \qquad [H:L]<\infty,
\]
where $\mathbf F(G\curvearrowright X)$ denotes the topological full group.

For every Nekrashevych-regular point $x\in X$, the natural orbital graphs
of $G$ and $H$ are decorated two-sided lines.  Let $\gamma_G(x)$ and
$\gamma_H(x)$ be the two-edge cell presentations of their edge-label
sequences, using the common dihedral indexing.  Then there is a
shift-equivariant homeomorphism
\[
 \Omega\bigl(\gamma_H(x)\bigr)
 \cong
 \Omega\bigl(\gamma_G(x)\vee w_X(x)\bigr)
\]
that sends $\gamma_H(x)$ to
$\gamma_G(x)\vee w_X(x)$, where $\vee$ denotes the
coordinatewise joint of sequences.
\end{mainthmx}

We call this relation between sequences \emph{mutual local derivability}
(MLD); see Section~\ref{sec:preliminaries}.  The proof of
Theorem~\ref{thm:main-orbital-coding} combines the containing fragmentation
of Theorem~\ref{thm:germ-preserving-containment}, the local comparison of
Proposition~\ref{pro:expansive-full-group-MLD}, and the finite-index
assertion of Proposition~\ref{pro:virtuality}.  We also use
Lemma~\ref{lem:free-universally-conull}, which identifies the regular
locus with the free dihedral locus.

Theorem~\ref{thm:main-orbital-coding} is the cornerstone of the constructions
below.  To describe its coding comparisons precisely and use them in the
constructions, Section~\ref{sec:preliminaries} introduces a framework of
codings that generalizes the clopen codings defined above.  This framework
is developed further in the subsequent sections and serves two purposes:
\begin{enumerate}[label=(\arabic*),leftmargin=2em,itemsep=0.25em,topsep=0.4em]
\item Using Definition~\ref{def:coding-singular} and
Convention~\ref{conv:maximal-continuous-extension},
Proposition~\ref{pro:fragmentation-coding-general} describes the orbital
graphs of fragmentations in terms of coding names.  In the setting of
Theorem~\ref{thm:main-orbital-coding}, this provides a common symbolic
framework for $\gamma_G$ and $\gamma_H$.
\item The graph semicocycle construction introduced in
Proposition--Definition~\ref{pro:U} provides tools for constructing almost
one-to-one extensions of the fragmented action that preserve the regular
orbital graphs of the original fragmentation
(Lemma~\ref{lem:regular-graph-extension}), while allowing us to prescribe
expansivity and repetitivity properties of the extended ambient action.
In particular, suitable non-clopen codings can be resolved by passing to a
corresponding extension in which the lifted closures of their coding pieces
are clopen, so that the natural lifted coding becomes a clopen coding on the
whole extension.  This mechanism is used in
Propositions~\ref{pro:resolved-entropy} and~\ref{pro:W-entropy} to provide a
common ambient space for the general entropy comparison, in
Proposition~\ref{pro:general-expansive-marker} and
Theorem~\ref{thm:continuum-simple-torsion} for the sparse-marker constructions,
and in Proposition~\ref{pro:Toeplitz} for the two-piece Toeplitz construction.
In the setting of Theorem~\ref{thm:main-orbital-coding}, these constructions allow us to control
the one-track joint $\gamma_G\vee w_X$, and hence $\gamma_H$, up to MLD.
\end{enumerate}

\paragraph{Positive entropy and exponential growth.}
Positive entropy of the ambient action gives exponential word growth
of the alternating full group.

\Needspace{11\baselineskip}
\begin{maintheoremi}[Positive entropy and exponential growth]
\label{thm:positive-entropy-criterion}
Let $\Dinf=\langle a,b\rangle\curvearrowright X$ be a minimal expansive
Cantor action, and let $G=\langle A\cup B\rangle\curvearrowright X$ be a
finite fragmentation.  If
\[
 h_{\rm top}(\Dinf\curvearrowright X)
 =\tfrac12 h_{\rm top}(ab\curvearrowright X)>0,
\]
then $\Aalt(G,X)$ is finitely generated, simple, amenable, and has exponential
word growth.

Moreover, if $G\curvearrowright X$ has a purely non-Hausdorff singularity,
then $\Aalt(G,X)$ is also torsion.
\end{maintheoremi}
The proof is given in Subsection~\ref{subsec:entropy-growth}.
Theorem~\ref{thm:positive-entropy} gives examples satisfying these hypotheses.

\Needspace{15\baselineskip}
\begin{mainproa}[Continuum many simple amenable torsion groups of exponential growth]
\label{thm:positive-entropy}
For every LR parameter $\omega$ and every $F:\N\to\N$, there are $2^{\aleph_0}$ minimal expansive Toeplitz
$\Dinf$-systems $X_\lambda$ with fragmentations
$G=G_\omega\curvearrowright X_\lambda$ that are almost
one-to-one extensions of the generalized Grigorchuk boundary action
$G\curvearrowright\partial T_2$.  Their labeled orbital
graphs at regular points are the same as those of
$G\curvearrowright\partial T_2$, and hence are LR; each action also has a PNH
singularity.  Moreover, the groups
\[
 L_\lambda=\Aalt(G_\omega,X_\lambda)
\]
are pairwise abstractly nonisomorphic, finitely generated,
simple, torsion, and amenable, with exponential word growth.
For every finite generating set of $L_\lambda$, its generic orbital
graph has repetitivity uniformly larger than $F$ up to linear rescaling.
\end{mainproa}

Fix an LR parameter $\omega$.  We start with $G=G_\omega$ acting on
$\partial T_2$, viewed as a fragmentation of the $\Dinf$-odometer, and
lift this action to a positive-entropy Toeplitz system $X$.
Lemmas~\ref{lem:regular-graph-extension} and~\ref{lem:canonical-fragmentation-lift}
give the conditions for lifting a fragmentation through an almost
one-to-one extension.  We construct the positive-entropy extension in
Lemma~\ref{lem:self-contained-positive-Toeplitz} and
Proposition~\ref{pro:positive-entropy-action} in
Subsection~\ref{sec:positive-growth}.

By Lemma~\ref{lem:regular-graph-extension}, the lift leaves the regular
orbital graphs of the fragmentation unchanged, so they are still LR.
On the other hand, let $H\supseteq\Aalt(G,X)$ be the containing
fragmentation of Theorem~\ref{thm:main-orbital-coding}.  Its orbital graph
is MLD with the joint of the orbital graph of $G\curvearrowright X$
and the Toeplitz sequence defining $X$.  The Toeplitz coordinate has
positive entropy, so this joint has positive entropy and exponential
pattern complexity; see Proposition~\ref{pro:joint-coding-entropy}.

Proposition~\ref{pro:entropy-growth} bounds this pattern complexity
in terms of word growth, so $H$ has exponential growth.  Since
$\Aalt(G,X)$ has finite index in $H$ by
Theorem~\ref{thm:main-orbital-coding}, it also has exponential growth.
This is the argument in Theorem~\ref{thm:positive-entropy-criterion};
Proposition~\ref{pro:positive-entropy-action} applies it to actions
satisfying \NekExpHyp.

To obtain continuum many pairwise abstractly nonisomorphic groups, we
use the sparse-marker coding rules of
Proposition~\ref{pro:general-expansive-marker} in
Subsection~\ref{sec:positive-growth}.
Theorem~\ref{thm:continuum-simple-torsion} constructs the further almost
one-to-one extensions required for Theorem~\ref{thm:positive-entropy}.
For every finite generating set, their regular orbital graphs have
repetitivity uniformly larger than the prescribed function $F$ up to
linear rescaling.  The same theorem describes the fibers, invariant
measures and entropy of the extensions.  To the best of our knowledge,
these are the first examples of finitely generated, infinite, simple,
amenable torsion groups of exponential word growth.

\paragraph{What the fragmentation coding determines.}
A partial coding may determine a point without a finite window
determining its clopen neighborhood.  In
Subsection~\ref{subsec:fragmentation-recovery}, we distinguish these two
notions by conditions on the fragmentation.  They describe what the
orbital graph of $G\curvearrowright X$ determines about the orbital
graph of $\Aalt(G,X)\curvearrowright X$ or of a finite containing
fragmentation $H\curvearrowright X$.  Condition $(W)$ requires the map taking a point
to its (germ) orbital graph to be injective on the set of regular and
Hausdorff singular points; in particular, it distinguishes regular points
(in the sense of Nekrashevych) by their orbital graphs.  Condition $(U)$
requires every finite clopen coding to be recoverable from a bounded
fragmentation window in the orbital graph.

Proposition~\ref{pro:W-entropy} shows that $(W)$ identifies the
$\Z$-system $ab\curvearrowright X$ with the $\Z$-subshift associated with
the orbital graph of $G\curvearrowright X$ (viewed as a sequence) at the
level of measured dynamics.  When $\Gamma(G,X)$ is LR, condition $(W)$
therefore forces $h_{\rm top}(\Dinf\curvearrowright X)=0$.
It does not, however, imply LR of either
$\Gamma(\Aalt(G,X)\curvearrowright X)$ or $\Gamma(H\curvearrowright X)$.

\paragraph{Zero entropy and non-transfer of linear repetitivity.}
Both the fragmentation's regular orbital graphs and the underlying
subshift can be LR without the orbital graphs of a containing
fragmentation being LR.  In the following construction, the latter have
arbitrarily large repetitivity.  The joint coding has zero entropy;
we do not determine the growth of the alternating full group.

\Needspace{9\baselineskip}
\begin{mainprob}[Failure of LR transfer despite LR of both the orbital graph and the subshift]
\label{thm:main-transfer-strong}
For every LR parameter $\omega$ and every prescribed function $F$, there is
an almost one-to-one binary Toeplitz extension
$\Dinf=\langle a,b\rangle\curvearrowright X$ of the dyadic odometer
$\Dinf\curvearrowright\partial T_2$, with a faithful lifted fragmentation
$G_\omega\curvearrowright X$, such that the underlying subshift
$ab\curvearrowright X$ and the orbital graphs at regular points are both LR,
and the fragmentation has a PNH singularity and satisfies $(W)$.
Nevertheless, the generic orbital graph associated with every finite
generating set of $\Aalt(G_\omega,X)$, as well as that of every finite
containing fragmentation, has repetitivity uniformly larger than $F$
up to linear rescaling.
\end{mainprob}
Theorem~\ref{thm:refined-full} in Subsection~\ref{sec:construction} gives
the full statement: part~(ii) concerns LR of the underlying subshift,
part~(iii) concerns LR of the fragmentation's orbital graphs and
condition $(W)$, and part~(v) gives the uniform repetitivity bounds.

Starting with an LR parameter $\omega$ and a prescribed function $F$,
we choose, in Lemma~\ref{lem:q-construction}, a binary level sequence
$q=(q_k)$ in which both symbols occur with bounded gaps, while
every pair $(i,\varepsilon)\in\{0,1,2\}\times\{0,1\}$ occurs infinitely
often as $(\omega_k,q_k)$.  For each regular phase
$z\in\Ztwo\setminus\Z$, the valuation words
\[
 u_\omega^z(n)=\omega_{\nu_2(z+n)},
 \qquad
 u_q^z(n)=q_{\nu_2(z+n)}
\]
are MLD with the orbital graphs
$\Gamma(G_\omega\curvearrowright\partial T_2,z)$ and
$\Gamma(G_q\curvearrowright\partial T_2,z)$,
respectively.\footnote{Here $G_q$ denotes the auxiliary fragmentation
obtained from the same rule with
$q\in\{0,1\}^{\N}\subset\{0,1,2\}^{\N}$; the third fragmentation
piece is empty.}
Both words are LR, but their joint has repetitivity uniformly larger
than $F$ up to linear rescaling; see
Proposition~\ref{pro:joint-nonLR}.

The binary sequence defines the two semicocycle pieces
\[
 Q_\varepsilon
 =\bigcup_{\{k:q_k=\varepsilon\}}
   \bigl(2^k\Ztwo\setminus2^{k+1}\Ztwo\bigr),
 \qquad \varepsilon\in\{0,1\}.
\]

The two pieces define the almost one-to-one Toeplitz extension
$\Dinf\curvearrowright X$ of Proposition~\ref{pro:Toeplitz}, analogous
to the point-explosion construction in
\cite[Introduction and Section~3]{NekSub25}.  It is the smallest minimal
extension in which the closures of the pieces lifted from
$\Ztwo\setminus\Z$ form a clopen partition; this partition generates
on $X$.  We lift the fragmentation to $G_\omega\curvearrowright X$
in Proposition~\ref{pro:lifted-all}.  For every regular
base point $z\in\Ztwo\setminus\Z$ with lift $\bar z$, the symbolic
coding of $X$ is MLD with $u_q^z$, whereas the fragmentation orbital
graph is MLD with $u_\omega^z$.  Theorem~\ref{thm:main-orbital-coding}
therefore identifies the containing fragmentation's orbital graph, up to
MLD, with $u_\omega^z\vee u_q^z$.
Proposition~\ref{pro:expansive-full-group-MLD} gives the same comparison
for the alternating full group's orbital graph.
The recurrence of every pair $(i,\varepsilon)$ ensures that both points
above $0^\infty$ are purely non-Hausdorff.  Finally,
Proposition~\ref{pro:joint-nonLR} bounds the joint's number of length-$n$
patterns by $12n$.  Thus the joint has zero entropy, although its repetitivity is
arbitrarily large.

\paragraph{Linear repetitivity and subexponential growth under $(U)$.}
Under \NekExpHyp, condition $(U)$ makes linear repetitivity pass to a
germ-preserving containing fragmentation.  The subexponential-growth result of \cite{Nek2018} then applies to this
fragmentation, and hence also to its alternating full subgroup.

\Needspace{7\baselineskip}
\begin{mainproii}[Subexponential growth from LR and condition $(U)$]\label{thm:compatibility}
Let $G=\langle A\cup B\rangle\curvearrowright X$ satisfy \NekExpHyp. If the
fragmentation also satisfies $(U)$, there is a germ-preserving
finite fragmentation $H\supseteq\Aalt(G,X)$ with LR generic orbital graphs.
Consequently, $\Aalt(G,X)$ is a finitely generated, infinite, simple,
torsion, amenable group of subexponential word growth; in particular,
it has intermediate growth.
\end{mainproii}
The proof is given at the end of Section~\ref{sec:machinery}.

\subsection*{Comparison of the hypotheses}
Under \NekExpHyp\ alone, the ambient action can have positive entropy and
its alternating full group can have exponential word growth, as in
Theorem~\ref{thm:positive-entropy}.  Adding $(W)$ forces both the ambient
entropy and the joint coding entropy to vanish by
Proposition~\ref{pro:W-entropy}, but the joint repetitivity can still be
uniformly larger than any prescribed function up to linear rescaling,
as shown by Theorem~\ref{thm:main-transfer-strong}.  With the stronger
condition $(U)$, LR passes to a containing fragmentation, and the
alternating full group has intermediate word growth by
Theorem~\ref{thm:compatibility}.

This distinction also concerns the expansive conclusion of
\cite[Theorem~1.2]{Nek2018}: although its hypotheses require LR
of the generic orbital graphs of the fragmentation $G$, this does not
by itself imply LR for a germ-preserving fragmentation containing
$\Aalt(G,X)$.  By Theorem~\ref{thm:main-orbital-coding}, the latter
orbital line is MLD with the joint coding $\gamma_G\vee w_X$, rather
than with $\gamma_G$ alone.  Theorem~\ref{thm:positive-entropy}
shows that LR of the $G$-orbital graphs alone does not force
intermediate growth of $\Aalt(G,X)$.  Under $(U)$, the additional
coding $w_X$ is recoverable from a bounded fragmentation window,
giving the intermediate-growth conclusion of
Theorem~\ref{thm:compatibility}.

\begin{rem}[The middle case]\label{rem:middle-case}
We have not determined the growth of $\Aalt(G,X)$ in
Theorem~\ref{thm:main-transfer-strong}.  We expect intermediate growth.
One possible approach is a joint traverse/portrait argument using the
two LR coordinates without requiring their joint to be LR.  Zero
orbital entropy alone does not imply subexponential word growth.
\end{rem}

\subsection*{Structure of the paper}
In Section~\ref{sec:preliminaries}, we recall codings, local derivations,
repetitivity and dyadic Toeplitz systems.  In Section~\ref{sec:machinery},
we prove the local comparisons and construct a containing fragmentation
$H\supseteq\Aalt(G,X)$.  Section~\ref{sec:orbital-entropy} compares
entropy and invariant measures, and proves that positive orbital
entropy implies exponential word growth.
Section~\ref{sec:extensions} begins with the lifting lemmas for almost
one-to-one extensions.  We then construct the positive-entropy examples
in Subsection~\ref{sec:positive-growth} and the zero-entropy examples
with arbitrarily large joint repetitivity in
Subsection~\ref{sec:construction}.

\section{Preliminaries}\label{sec:preliminaries}
\paragraph{Group notation.}
We use $\Lambda$ and $K$ for arbitrary acting groups; our applications
of statements about $K$ concern $\Z$ or $\Dinf$.
We reserve $G$ for a fragmentation of a dihedral action, $H$ for a
containing fragmentation of the alternating full group, and $L$ for
an alternating full group.

\subsection*{Coding conventions and orbital graphs at regular points}
\begin{Def}[Coding rules, codings, and coding names]
\label{def:orbit-coding}\label{def:coding}
For an action $K\curvearrowright X$ and a finite alphabet $\mathcal A$,
a \emph{coding rule}, or \emph{central-symbol map}, is a continuous map
$p:D\to\mathcal A$ such that
\[
 D_p(K)=\{x:K x\subseteq D\}=\bigcap_{g\in K}g^{-1}D
\]
is nonempty.  Its associated \emph{coding} is the map
\[
 c_p:D_p(K)\longrightarrow\mathcal A^K,
 \qquad c_p(x)(g)=p(gx).
\]
For $x\in D_p(K)$, the value $c_p(x)\in\mathcal A^K$ is the
\emph{coding name at $x$}; its central symbol is $c_p(x)(e)=p(x)$,
where $e$ denotes the identity of $K$.

We call $p$ \emph{total} if $D_p(K)=D=X$, and \emph{partial} otherwise.
The coding name $c_p(x)$ is constant on the cosets $gK_x$.
Changing the root from $x$ to $hx$ gives
\[
 c_p(hx)=h\cdot c_p(x),\qquad
 (h\cdot u)(g)=u(gh).
\]
Thus $c_p$ is equivariant, and it is continuous on its relative domain.

A finite clopen partition $\mathcal Q=\{Q_\alpha\}$ defines the total
coding rule $p_{\mathcal Q}(x)=\alpha$ for $x\in Q_\alpha$.
We call $c_{p_{\mathcal Q}}$ the associated \emph{$\mathcal Q$-clopen coding}
and $c_{p_{\mathcal Q}}(x)$ the \emph{$\mathcal Q$-coding name at $x$}.
The partition is \emph{generating} when $c_{p_{\mathcal Q}}$ is injective.
\end{Def}

A coding name $c_p(x)$ can be displayed by decorating each vertex
$y\in K x$ of a labeled orbital graph with the symbol $p(y)$;
we call the resulting decorated graph the \emph{coding graph at $x$}.
For two codings $c_p$ and $c_{p'}$, their \emph{joint coding} is
\[
 (c_p\vee c_{p'})(x)=\bigl(c_p(x),c_{p'}(x)\bigr)
 \qquad\bigl(x\in D_p(K)\cap D_{p'}(K)\bigr),
\]
where the product of the two name spaces is identified coordinatewise
with the name space over the product alphabet.  Its value at $x$ is
the \emph{joint coding name at $x$}.
The \emph{subshift associated with $p$} is
\[
 Z_p=\overline{c_p(D_p(K))}\subseteq\mathcal A^K,
\]
the closure of all its coding names.  It is minimal if
$K\curvearrowright X$ is minimal.

For general (possibly partial) dihedral coding rules, the following
reindexing is useful.

\begin{convention}[Two-track reindexing for general codings]
\label{conv:cell-subshift}
For $\Dinf\curvearrowright X$, put $T=ab$.  The decomposition
$\Dinf=\{T^n:n\in\Z\}\sqcup\{bT^n:n\in\Z\}$ defines a reindexing homeomorphism
\[
 J:\mathcal A^{\Dinf}\longrightarrow(\mathcal A^2)^\Z,
 \qquad J(c)(n)=\bigl(c(T^n),c(bT^n)\bigr).
\]
For $x\in D_p(\Dinf)$, set $\kappa_p(x)=J(c_p(x))$, so that
\[
 \kappa_p(x)(n)=\bigl(p(T^nx),p(bT^nx)\bigr),
 \qquad
 Z^{\mathrm{cell}}_p:=\overline{\kappa_p(D_p(\Dinf))}=J(Z_p).
\]
We call $Z^{\mathrm{cell}}_p\subseteq(\mathcal A^2)^\Z$ the \emph{cell subshift}.
The restriction of $J$ conjugates $T\curvearrowright Z_p$ with
$\sigma\curvearrowright Z_p^{\mathrm{cell}}$; this is a conjugacy between the
two symbolic presentations.  If $p$ is total and generating, $\kappa_p$
also conjugates $T\curvearrowright X$ with this cell subshift.
The reflection $b$ becomes reversal followed by exchanging the two entries:
\[
 J(b\cdot c)(n)=\operatorname{swap}\bigl(J(c)(-n)\bigr).
\]
If the coding rule $p$ is $b$-invariant, the two entries agree and the
single word $n\mapsto p(T^nx)$ contains the same information.  We retain $J$ for partial and non-$b$-invariant coding rules.
\end{convention}
\begin{Def}[Regularity and linear repetitivity]\label{def:action-LR}
For a finitely generated action $\Lambda\curvearrowright X$, a point is
\emph{regular (in the sense of Nekrashevych)} if $gx=x$ implies
$x\in\Int\Fix(g)$ for every $g\in\Lambda$.  The regular locus is an invariant
dense $G_\delta$ for a countable Cantor action: its complement is the union
of the nowhere dense sets $\Fix(g)\setminus\Int\Fix(g)$.

For a repetitive word $u$, let its \emph{repetitivity} $R_u(n)$ be the
least $R$ such that every length-$R$ subword contains every
length-$n$ word. We also call the finite subwords \emph{patterns}.  The word is
LR when $R_u(n)\leq Cn$.  For a repetitive bounded-degree labeled graph
$\Gamma$, let its \emph{repetitivity} $R_\Gamma(r)$ be the least $R$
such that every radius-$R$ ball contains
an occurrence of every rooted radius-$r$ ball type appearing in $\Gamma$,
with each occurrence understood as the full ambient rooted ball.  Edge
labels, directions, and loops are retained; distance is taken in the
underlying undirected graph.

A minimal action is LR when its orbital graph
at one (equivalently, every) regular point is LR for a fixed finite
generating set.
\end{Def}

\begin{lem}[Mutual local containment of regular orbital graphs]
\label{lem:LR-one-all}
Fix a finite generating set for a minimal Cantor action.  The orbital graphs
of any two regular points are \emph{mutually locally contained}: every
finite rooted labeled ball of either graph occurs in the other.  Thus they
have the same finite-pattern language, pattern-complexity function,
and repetitivity function.
\end{lem}
\begin{proof}
To determine a finite rooted orbital ball, it is enough to check finitely
many equalities and inequalities $ux=vx$ among bounded-length words.
We enlarge the radius by one if necessary to record boundary edges.
At a regular point, the equalities hold throughout a neighborhood;
the inequalities also hold nearby, by continuity.  Thus the rooted-graph map is continuous at every regular point.
By minimality, translates of either point approach the other.  Continuity
then gives local containment in both directions.  Every finite labeled pattern is contained in a finite ball,
so the pattern languages and their complexity functions agree as well.
Whether every radius-$R$ ball contains every radius-$r$ pattern depends only
on that language.
\end{proof}

\Needspace{7\baselineskip}
\subsection*{Local derivations and quantitative repetitivity}
\label{subsec:local-derivations}
\begin{Def}[Local derivation and MLD]\label{def:local-derivation}
For words $u,v:K\to\mathcal A$, $v$ is \emph{locally derivable}
from $u$ if there is a fixed finite set $F\subseteq K$ such that
$v(g)$ is determined by $u|_{\{fg:f\in F\}}$ for every
$g\in K$.  Equivalently, there is a map
$\phi:\mathcal A^F\to\mathcal A$ such that
$v(g)=\phi(f\mapsto u(fg))$ for every $g\in K$.
The words are \emph{mutually
locally derivable}, abbreviated MLD, if there are such rules in both
directions on this same indexing.

For labeled bounded-degree graphs $\Gamma,\Gamma'$, a local derivation
of $\Gamma$ from $\Gamma'$ \emph{via a bi-Lipschitz bijection}
$f:V(\Gamma)\to V(\Gamma')$ means that there are constants $C\geq1$
and $r\geq0$ such that, for every $R$ and $v$, the rooted ball
$B_\Gamma(v,R)$ is determined uniformly by
$B_{\Gamma'}(f(v),CR+r)$, with the vertices identified through $f$.
The rule respects rooted labeled isomorphisms and locates every recovered
edge endpoint as the endpoint of a bounded-length path in the source
graph.  Equivalently, one reads the target decorations and neighboring
vertices by one bounded local rule and then iterates it.  After a fixed
rescaling of the source metric the bound has the form $R+r$; with the
intrinsic graph metrics the factor $C$ allows for different edge lengths.
The graphs are MLD if one bi-Lipschitz bijection $h$ gives both local
derivations, from $\Gamma'$ to $\Gamma$ via $h$ and in the reverse
direction via $h^{-1}$.\footnote{The common identification is part of
MLD.  Local derivations via two independently chosen bijections $f$ and
$g$ do not provide this datum unless the identifications are compatible;
one must use a map and its inverse, not unrelated reindexings.}

Our graph-MLD comparisons concern orbital graphs of actions with
the same point-orbits in one space.  The vertex bijection is always
the identity on the common orbit.
The notation $\Gamma\overset{\MLD}{\cong}u$ regards $u$ as a decorated
nearest-neighbor line.  A two-edge cell presentation retains the two
vertices of each cell.
\end{Def}

For bi-infinite words the ambient indexing is $K=\Z$.

\begin{lem}[One-point criterion for MLD of orbit codings]
\label{lem:one-point-MLD}
Let $K\curvearrowright X$ be minimal, and let
$p:D_p\to\mathcal A$ and $q:D_q\to\mathcal B$ be continuous coding
rules.  Put $D=D_p(K)\cap D_q(K)$ and assume $D\ne\varnothing$.
If $c_p(x)$ and $c_q(x)$ are MLD at one $x\in D$ by finite-window
rules on the common $K$-indexing, then the same
two rules make $c_p(y)$ and $c_q(y)$ MLD for every $y\in D$.
\end{lem}
\begin{proof}
For the right-shift convention in Definition~\ref{def:orbit-coding},
the two rules give finite sets $F,E\subseteq K$ and maps
$\phi:\mathcal A^F\to\mathcal B$ and
$\psi:\mathcal B^E\to\mathcal A$ such that, for every $g\in K$,
\[
 q(gx)=\phi\bigl((p(fgx))_{f\in F}\bigr),\qquad
 p(gx)=\psi\bigl((q(egx))_{e\in E}\bigr).
\]
The orbit $K x$ is dense in $X$ and lies in the invariant set $D$, so it
is dense in $D$ in the relative topology.  Both sides are continuous on $D$, so the identities extend from $Kx$
to all of $D$.  Applying them at $gy$ gives the required rules for
$c_p(y)$ and $c_q(y)$.
\end{proof}

Under the hypotheses of Lemma~\ref{lem:one-point-MLD}, for codings
$c,c'$ with a nonempty intersection of domains, we write
$c\overset{\MLD}{\cong}c'$ if
$c(x)\overset{\MLD}{\cong}c'(x)$ at one, equivalently every, point
of that intersection.

Thus a local derivation between codings uses the same finite-window rule
at every point of their common domain.  When the root $x$ in a graph
comparison is understood, we write the coding symbol for its name at $x$.

\begin{lem}[Local derivations and joints]\label{lem:LD}\label{lem:join-MLD}
A local derivation of words or graphs preserves LR.  Moreover,
$a\overset{\MLD}{\cong}a'$
implies $a\vee b\overset{\MLD}{\cong}a'\vee b$ for a fixed word $b$.
\end{lem}
\begin{proof}
A radius-$r$ rule gives
$R_v(n)\leq R_u(n+2r)+2r$, hence preserves a linear bound.  For the assertion about joints, apply the local rule to the first
coordinate and leave the second unchanged.  The inverse rule gives the
other direction.
\end{proof}

\begin{pro}[Clopen codings and finite windows]\label{lem:CHL}
\label{pro:continuous}\label{rem:subshift-coding}
An expansive Cantor action of a finitely generated group admits a finite
clopen generating partition. Its name map conjugates the action to a
subshift.  Every finite clopen coding is locally derivable from a generating
clopen coding, and any two such generating codings are MLD.  Equivalently,
a continuous equivariant map between finite-alphabet subshifts is a
finite-window map; a conjugacy has finite-window rules in both directions.
\end{pro}
\begin{proof}
Choose a finite clopen partition with atom diameters smaller than an
expansivity constant.  Two points with the same coding name cannot be separated by the
action, so expansivity implies that they are equal.  The name map is
therefore injective; by compactness, it is a conjugacy onto its image.  For a continuous equivariant map between subshifts, the inverse image
of each central target symbol is clopen, so it is a finite union of
cylinders.  There are only finitely many target
symbols; taking a window large enough for all of them gives the local
rule at the identity.  Equivariance gives the rule at every coordinate.
For a conjugacy, we apply the same argument to its inverse.  We can also
apply this argument to the symbol map of any finite clopen partition
to recover its coding from a generating coding.
\end{proof}

For two bi-infinite sequences $u$ and $v$ with the same $\Z$-indexing,
MLD is equivalent to the existence of a shift-equivariant homeomorphism
$\Phi:\Omega(u)\to\Omega(v)$ satisfying $\Phi(u)=v$.
Indeed, the finite-window rules defining MLD extend continuously to the
orbit closures, and their compositions are the identity on the dense
shift orbits, hence everywhere.  Conversely, a conjugacy and its inverse
are finite-window maps by Proposition~\ref{pro:continuous}.
The condition $\Phi(u)=v$ specifies the correspondence between the two
chosen sequences and their common indexing.

\begin{lem}[A reflection-invariant generating partition]
\label{lem:reflection-invariant-generator}
Let $\Dinf=\langle a,b\mid a^2=b^2=1\rangle$ act expansively on a
Cantor set $X$, and put $T=ab$.  If $T$ has no points of period at
most two, then there exists a finite clopen partition $\mathcal Q_X$
generating for $T$ whose every atom is $b$-invariant.  Its symbol map
$f_X:X\to\mathcal B_X$ therefore gives a conjugacy
\[
 w_X:X\longrightarrow w_X(X)\subseteq\mathcal B_X^{\Z},\qquad
 w_X(x)(n)=f_X(T^nx),
\]
for which $w_X(Tx)=\sigma w_X(x)$ and
$w_X(bx)(n)=w_X(x)(-n)$.
In particular the conclusion holds for every minimal expansive
$\Dinf$-action on an infinite Cantor set.
\end{lem}
\begin{proof}
First, expansivity passes from the dihedral action to $T$.  Indeed, if
all corresponding $T$-iterates of two points are sufficiently close,
uniform continuity of $b$ makes their images under all $bT^n$ close
as well.  Expansivity of the dihedral action then forces the two points
to coincide.  Fix a finite clopen $T$-generating partition $\mathcal P$.

Let $q:X\to Q=X/\langle b\rangle$ be the quotient map.  As $b$ is an
involution of a compact zero-dimensional space, $Q$ is compact and
zero-dimensional.  Consider
\[
 F:X\longrightarrow Q\times Q,\qquad F(x)=(q(x),q(Tx)).
\]
We claim that $F$ is injective.  If $F(x)=F(y)$, the first coordinate
gives $y=x$ or $y=bx$.  Suppose that $y=bx\ne x$.  Equality of the
second coordinates gives $q(Tx)=q(Tbx)=q(T^{-1}x)$.  Thus either
$Tx=T^{-1}x$, which gives $T^2x=x$, or $Tx=bT^{-1}x=Tbx$, which gives
$bx=x$.  The first contradicts the hypothesis on $T$, and the second
contradicts $y\ne x$.  Hence $F$ is injective and, by compactness, is a
topological embedding.

Cover $F(X)$ by clopen rectangles on whose intersections with $F(X)$
the $\mathcal P$-symbol is constant.  Such rectangles exist because $F$ is
an embedding and $Q$ is zero-dimensional.  Take a finite subcover,
and let $\mathcal C$ be a common clopen refinement of all its sides.
Then the $\mathcal C$-symbols of $q(x)$ and $q(Tx)$ determine the
$\mathcal P$-symbol of $x$.  Set $\mathcal Q_X=q^{-1}\mathcal C$ and let $f_X$ be its
symbol map.  Since $q\circ b=q$, every atom of $\mathcal Q_X$ is
$b$-invariant.  Moreover, each $\mathcal P$-symbol at $T^nx$ is
read from $(f_X(T^nx),f_X(T^{n+1}x))$, so the $f_X$-itinerary is
$T$-generating.  Finally $T^nb=bT^{-n}$ gives ordinary reversal.
For a minimal dihedral action on an infinite Cantor set, a $T$-periodic
point would have a finite dihedral orbit, contradicting minimality.
\end{proof}

\begin{convention}[Expansive presentations and invariant one-track coding]
\label{conv:expansive-presentations}
Whenever an action $K\curvearrowright X$ admits a finite clopen generating
partition, fix one such partition $\mathcal Q_X$.  Let $p_X$ be its
central-symbol map and let
\[
 \chi_X=c_{p_X}:X\longrightarrow{\mathcal Q_X}^K
\]
be its coding, using the atoms as the alphabet.  Thus $\chi_X(x)$ is the
coding name at $x$ and $p_X(x)=\chi_X(x)(e)$ is its central symbol.
By Proposition~\ref{pro:continuous}, $\chi_X$ conjugates the action
with its image subshift.

When $K=\Dinf=\langle a,b\rangle$ acts minimally and expansively on an
infinite Cantor set, choose $\mathcal Q_X$ as in
Lemma~\ref{lem:reflection-invariant-generator}, so that
$p_X(bx)=p_X(x)$.  In this case we use the one-track name
\[
 w_X(x)(n):=p_X(T^nx),\qquad T=ab,
\]
which is generating for $T$ and sends $b$ to ordinary reversal.  The
full dihedral name $\chi_X$ is recovered from $w_X$ by
$\chi_X(x)(T^n)=\chi_X(x)(bT^n)=w_X(x)(n)$.

For a given subshift the zero-coordinate partition is also a generating
choice; when it is not $b$-invariant, we use an equivalent invariant
recoding for $w_X$.  The subscript records the acting space: for
$\widehat X$ and $Y$, write $p_{\widehat X},\chi_{\widehat X}$ and
$p_Y,\chi_Y$, respectively, and use $w_{\widehat X},w_Y$ in the
minimal expansive dihedral setting.  An unsubscripted $\chi$ denotes
an arbitrary finite clopen coding, not necessarily generating; its
central-symbol map is $p_\chi(x)=\chi(x)(e)$.
\end{convention}

We will use the following comparison with a prescribed repetitivity function.

\begin{Def}[Repetitivity uniformly larger than a prescribed function]
\label{def:worse-F}\label{def:uniformly-larger-F}
Let $F:\N_{\geq1}\to\N$ be nondecreasing.  We say that $R_\Gamma$ is
\emph{uniformly larger than $F$ up to linear rescaling} if
\[
 \forall C\in\N_{\geq1}\ \exists n_0(C)\quad
 \forall n\geq n_0(C),\qquad
 R_\Gamma(n)>C F(Cn+C)+C.
\]
The finitely many exceptional scales may depend on $C$ and on the
graph presentation,
in particular on its generating set.  For an arbitrary prescribed
function $F$, we use its nondecreasing envelope
\[
 F^\uparrow(n)=\max_{1\leq m\leq n}F(m).
\]

\end{Def}

\begin{lem}[Quantitative comparison and all-scale lower bounds]
\label{lem:MLD-repetitivity}
Let $\Gamma_1,\Gamma_2$ be repetitive bounded-jump graphs on a common
line-like orbit, with uniformly comparable metrics.  If $\Gamma_1$ is
locally derivable from $\Gamma_2$, there is an integer $a\geq1$ such that
\[
 R_{\Gamma_2}(an+a)\geq a^{-1}R_{\Gamma_1}(n)-a
 \qquad(n\geq1).
\]
Consequently the property in Definition~\ref{def:worse-F} passes from
$\Gamma_1$ to $\Gamma_2$; in particular it is invariant under graph-MLD.

For a prescribed $F$, put
\[
 \widehat F(n)=n^2\bigl(F^\uparrow(n^2)+1\bigr).
\]
If $R_{\Gamma_1}(n)\geq\widehat F(n)$ for every sufficiently large $n$,
then $R_{\Gamma_2}$ is uniformly larger than $F$ up to linear rescaling.
It suffices to have an increasing unbounded sequence $(n_j)$ with
\[
 R_{\Gamma_1}(n_j)>\widehat F(n_{j+1})
 \qquad\text{for every sufficiently large }j.
\]

\end{lem}
\begin{proof}
A radius-$n$ pattern in $\Gamma_1$ is recovered from a radius-$an+a$
pattern in $\Gamma_2$, for a uniform $a$.  An occurrence of the latter
produces an occurrence of the former.  Comparing the two orbital metrics
therefore gives the displayed repetitivity inequality after increasing
$a$.  Suppose that $R_{\Gamma_1}$ satisfies Definition~\ref{def:worse-F}.
Apply its lower bound with a rescaling constant large enough to account
for $a$ and the prescribed rescaling of $R_{\Gamma_2}$.  The displayed
inequality then gives the same property for $R_{\Gamma_2}$.  If the
graphs are MLD, we can interchange them to obtain the converse.

For the bound involving the majorant, fix $C\geq1$ and, for large $m$, set
$n=\lfloor(m-a)/a\rfloor$.  Then $an+a\leq m<a(n+2)$, so
\[
 R_{\Gamma_2}(m)\geq a^{-1}\widehat F(n)-a.
\]
Eventually $n^2\geq Cm+C$ and $n^2/a-a>C$.  Since $F^\uparrow$ is
nonnegative and nondecreasing, the right-hand side is then strictly
larger than $C F^\uparrow(Cm+C)+C$.  This holds for each fixed $a,C$;
the majorant itself is independent of the target graph and generating set.

Finally, if $n_j\leq n<n_{j+1}$, monotonicity gives
\[
 R_{\Gamma_1}(n)\geq R_{\Gamma_1}(n_j)
 >\widehat F(n_{j+1})\geq\widehat F(n).
\]
Thus the bound holds at every sufficiently large scale.  The factor
$n^2$ makes it superlinear even when $F$ is bounded or linear.
\end{proof}

\subsection*{Dyadic codings and Toeplitz phases}
\label{subsec:valuation}\label{subsec:toeplitz-framework}

Identify $\partial T_2$ with $\Ztwo$ by the Gray code so that
\[
 b(z)=-z,\qquad a(z)=1-z,\qquad ab(z)=z+1,
\]
and put
\[
 W_k=2^k\Ztwo\setminus2^{k+1}\Ztwo.
\]
For a finite-alphabet sequence $h=(h_k)$ and
$z\in\Ztwo\setminus\Z$, define the $\mathbb Z$-word $u_h^z$ by
\[
     u_h^z(n)=h_{\nu_2(z+n)}.
\]
Given $\omega\in\{0,1,2\}^{\N}$, define the pieces
$P_i^\omega=\bigcup_{\omega_k=i}W_k$ and fragment the reflection $b$ by
the Klein four group whose nontrivial elements act as $b$ on two of the
three pieces.  Together with the unfragmented reflection $a$ this is the
standard boundary action of $G_\omega$; its orbital graph at the
regular point $z$ is MLD with, and naturally identified with,
$u_\omega^z$.  If every symbol occurs infinitely often, all
three pieces accumulate at $0$.  Each nontrivial element of the Klein
four group is the identity on one accumulating piece and agrees with $b$
on the other two, so $0$ is PNH.  Its germ group is the Klein four group,
by \cite[Lemma~3.5]{Nek2018}; the same holds along its orbit.

\begin{pro}[Exact LR criterion for valuation words]
\label{pro:valuationLR}
The words $u_h^z$, uniformly in $z\in\Ztwo\setminus\Z$, are LR if and only
if there is $M$ such that every letter occurring in the tail
$h_rh_{r+1}\ldots$ already occurs among $h_r,\ldots,h_{r+M}$.
\end{pro}

\begin{proof}
Suppose the tail bound is $M$, and choose $r$ so that $2^{r-1}\leq L<2^r$.
An interval of length $L$ contains at most one position of valuation
at least $r$.
All other symbols are determined by the residue modulo $2^r$ and the list
$h_0,\ldots,h_{r-1}$.  If the block contains a position of valuation at least $r$, its symbol
occurs at some level between $r$ and $r+M$.  The same block therefore
reappears at centers separated by at most $2^{r+M+1}=O(L)$.
This proves LR.

Conversely, if the tail bound fails, choose a letter $\alpha$ and levels
$r_j<k_j$ such that $h_{k_j}=\alpha$, $h_s\neq\alpha$ for
$r_j\leq s<k_j$, and $k_j-r_j\to\infty$.  Center a block of radius
$2^{r_j}$ at valuation $k_j$.  A return center cannot have valuation in
$[r_j,k_j)$, because its central symbol would be wrong.  If its valuation
is $s<r_j$, equality of the central symbols forces $h_s=\alpha$.
Choose an offset $t$, $|t|\leq2^{r_j}$, at which that low center has
valuation $r_j$.  Then $\nu_2(t)=s$.  At this offset from the original center, the label
is $h_s=\alpha$; at the proposed return, it is $h_{r_j}\neq\alpha$.
The two blocks therefore differ.  Every return center must have
valuation at least $k_j$, so successive returns are at least $2^{k_j}$
apart.  Since the block has diameter $O(2^{r_j})$ and
$k_j-r_j\to\infty$, the ratio of the return gap to the block diameter
is unbounded.
\end{proof}

For a three-letter sequence in which all letters recur infinitely, this is
exactly the bounded-gap condition $(\mathrm{LR}_\omega)$.

\paragraph{Toeplitz points, period structures and phases.}
For the combinatorics and repetitivity of one-dimensional simple Toeplitz
subshifts, see also \cite{Sell20}.
For a word $v\in\mathcal C^\Z$, put
\[
 \operatorname{Per}_m(v)=\{n:v(n+\ell m)=v(n)\text{ for all }\ell\in\Z\}.
\]
The word is Toeplitz if these periodic parts cover $\Z$.  The $m$-skeleton
retains $\operatorname{Per}_m(v)$ and puts holes elsewhere; an essential
period is the least period of a nonempty skeleton.  A period structure is
a sequence of essential periods $s_j\mid s_{j+1}$ whose periodic parts
cover $\Z$.  A nonperiodic Toeplitz word has a minimal orbit closure and
an almost one-to-one factor onto its maximal odometer
\[
 O=\varprojlim_j\Z/s_j\Z,\qquad \rho(\sigma y)=\rho(y)+1.
\]
The value $\rho(y)$ is the phase; Toeplitz points are exactly those in
singleton fibers.  If all coordinate periods can be chosen as powers of
two, one may choose $s_j=2^{r_j}$.  Nonperiodicity forces $r_j\to\infty$,
so $O\cong\Ztwo$; not every power of two needs to be essential.  These are
the Toeplitz factor and period-structure theorems of
\cite[Sections~5 and~7, especially Theorems~5.1 and~7.4]{Dow05}.
Translating the odometer changes the origin of the phase map without
changing its equivariance.

\section{Fragmentation codings and orbital graphs of alternating full groups}\label{sec:machinery}
\subsection{Fragmentation codings, conditions \texorpdfstring{$(W)$ and $(U)$}{(W) and (U)}}
\label{subsec:fragmentation-recovery}

For the total generating codings of
Convention~\ref{conv:expansive-presentations}, every clopen coding is
locally derivable from the generating one.  This need not hold for
partial coding rules.  Using the local derivations of
Section~\ref{sec:preliminaries}, we distinguish injectivity of a coding
from recovery by a finite window.  We first define the singular set of
a partial coding rule.

\begin{Def}[Maximal extension and coding singularities]
\label{def:coding-singular}
Let a countable group $\Lambda$ act on a Cantor set $X$, let $\mathcal A$ be
finite, let $D^\circ\subset X$ be dense and open, and let
\[
    p^\circ:D^\circ\longrightarrow\mathcal A
\]
be continuous.  A point $x\in X$ is called \emph{locally extendable} for
$p^\circ$ if there are a neighborhood $U$ of $x$ and a continuous map
$p_U:U\to\mathcal A$ such that
\[
    p_U|_{U\cap D^\circ}=p^\circ|_{U\cap D^\circ}.
\]
Let $E_0(p^\circ)$ be the set of locally extendable points and define
\[
    \Sigma_0(p^\circ):=X\setminus E_0(p^\circ).
\]
We call $\Sigma_0(p^\circ)$ the \emph{primary coding-singular set}.

Because $D^\circ$ is dense and the alphabet is discrete, every local
extension is unique, and the local extensions glue to a unique continuous
map
\[
    p:X\setminus\Sigma_0(p^\circ)\longrightarrow\mathcal A
\]
extending $p^\circ$.  We call $p$ the \emph{maximal continuous extension}
of the raw coding rule.  Thus $\Sigma_0(p^\circ)$ is exactly the set of points to
which the coding rule cannot be continuously extended; equivalently, it is the
smallest closed subset $\Sigma\subset X$ for which $p^\circ$ extends
continuously to $X\setminus\Sigma$.  We abbreviate it to $\Sigma_0(p)$ once
the maximal extension has been fixed.

The singular set of the associated coding is the orbit saturation
\[
    \Sigma_\Lambda(p)=\bigcup_{g\in\Lambda}g\Sigma_0(p),
\]
and the corresponding coding-regular domain is
\[
    D_p=X\setminus\Sigma_\Lambda(p)
        =\{x\in X:\Lambda x\subseteq X\setminus\Sigma_0(p)\}
        =D_p(\Lambda)
\]
in the notation of Definition~\ref{def:orbit-coding}.  Since $\Lambda$ is
countable and $\Sigma_0(p)$ is closed and nowhere dense, $D_p$ is an
invariant dense $G_\delta$ set.  For $x\in D_p$, the coding name is given by
\[
    c_p(x)(g)=p(gx),
    \qquad g\in\Lambda.
\]
We call points of $\Sigma_\Lambda(p)$ \emph{$p$-singular} and points of
$D_p$ \emph{$p$-regular}.
\end{Def}

\begin{proof}[Existence and maximality of the extension]
By definition, $E_0(p^\circ)$ is open and contains $D^\circ$.
Two local extensions agree on $D^\circ$ wherever both are defined.
Since $D^\circ$ is dense, continuity makes them agree on their whole
overlap.  We can therefore glue them to a unique continuous extension
on $E_0(p^\circ)$.  If $p^\circ$ extends continuously to $X\setminus\Sigma$ for
some closed set $\Sigma$, then every point of $X\setminus\Sigma$ is locally
extendable, so
\[
    X\setminus\Sigma\subseteq E_0(p^\circ)
    \quad\text{and hence}\quad
    \Sigma_0(p^\circ)\subseteq\Sigma.
\]
This proves the minimality statement.
\end{proof}

\begin{convention}[Maximal continuous extension]
\label{conv:maximal-continuous-extension}
Henceforth, we replace every coding rule $p$ by its maximal continuous extension.

In this setting, $p$ is partial if $\Sigma_0(p)\neq \emptyset$.

\end{convention}

\Needspace{8\baselineskip}
\begin{pro}[The natural fragmentation coding and the non-Hausdorff locus]
\label{pro:fragmentation-coding-general}
Let $\Dinf=\langle a,b\rangle$ act minimally on an infinite Cantor set $X$,
and let
\[
    G=\langle A\cup B\rangle
\]
be a finite fragmentation, where $A$ fragments $a$ and $B$ fragments $b$.
On the dense open set
\[
    D^\circ=X\setminus\bigl(\Fix(a)\cup\Fix(b)\bigr),
\]
define, for $h\in A$ and $k\in B$,
\[
    hx=a^{\varepsilon_h(x)}x,
    \qquad
    kx=b^{\delta_k(x)}x,
\]
and put
\[
    p^\circ(x)
    =\bigl((\varepsilon_h(x))_{h\in A},
            (\delta_k(x))_{k\in B}\bigr).
\]
Let $p$ be its maximal continuous extension from
Definition~\ref{def:coding-singular}.  We call $p$ the \emph{natural
fragmentation coding rule}, and in this case we write its coding $c_p$ as $c_G$.  Then:
\begin{enumerate}[label=(\roman*)]
    \item At every Nekrashevych-regular point $x$, the coding graph defined
    by $p$ is canonically the labeled orbital graph of the fragmentation action
    \[
        \Gamma(G\curvearrowright X,x;A\cup B),
    \]
    This means that the underlying graphs are identified, and the vertex labels
    of the former determine the labels of the adjacent edges and loops of the
    latter, and vice versa.

    At a Hausdorff singularity, $p$ is defined at every orbit vertex and
    records, at a vertex fixed by $a$ or $b$, whether a fragment loop has
    the trivial germ or the germ of that reflection.  The ordinary orbital
    graph records only the loop.

    \item Let $\operatorname{NH}(G\curvearrowright X)$ denote the set of
    non-Hausdorff singularities in the sense of Definition~2.2 of
    \cite{Nek2018}.  Then
    \[
        \boxed{
        \Sigma_{\Dinf}(p)
        =\operatorname{NH}(G\curvearrowright X).
        }
    \]
    More precisely, if $\xi\in\Sigma_0(p)$, then exactly one of $a$ and $b$
    fixes $\xi$, and some fragment represents a non-trivial germ at $\xi$
    whose fixed-point interior accumulates on $\xi$.
\end{enumerate}
\end{pro}

\begin{proof}
First note that every non-trivial element of $\Dinf$ has fixed-point set with
empty interior.  Put $t=ba$.  If $t^n x=x$ for some $n\neq0$, then the
$\Dinf$-orbit of $x$ is finite, contradicting minimality on an infinite
space.  Thus non-trivial translations have no fixed points.  Every reflection
is conjugate to $a$ or $b$.  Suppose, for example, that a non-empty open set
$U$ were fixed pointwise by $a$, and choose $x\in U$.  Since $ax=x$, the
$\Dinf$-orbit of $x$ is its $\langle t\rangle$-orbit and is dense.  Choose
$n\neq0$ with $t^n x\in U$.  Using $at^n=t^{-n}a$, we obtain
\[
    t^n x=at^n x=t^{-n}ax=t^{-n}x,
\]
so $t^{2n}x=x$, a contradiction.  The argument for $b$ is identical.  Hence
all non-trivial reflections have fixed-point sets with empty interior.  It
follows in particular that
$D^\circ=X\setminus(\Fix(a)\cup\Fix(b))$ is dense and open.

For $x\notin\Fix(a)$, the two possible values $x$ and $ax$ are distinct.
By continuity, a fragment $h\in A$ must make the same choice, $hx=x$
or $hx=ax$, throughout a sufficiently small neighborhood of $x$.
This gives a locally constant, $a$-invariant profile.  We apply the
same argument to $B$ and $b$.  At a Nekrashevych-regular vertex, the $A$-profile records
exactly which fragment generators give loops and which cross the $a$-edge;
    the $B$-profile records the same data for the $b$-edge.  This proves the
    regular-point assertion in (i).

Let $\xi\in\Sigma_0(p)$.  If neither $a$ nor $b$ fixed $\xi$, then $\xi$
would lie in $D^\circ$, where the raw profile is already continuous.  The two
involutions cannot both fix $\xi$, since then the entire $\Dinf$-orbit of
$\xi$ would be the singleton $\{\xi\}$.  Hence exactly one fixes $\xi$;
suppose it is $a$.

The $B$-profile extends near $\xi$, whereas the $A$-profile does not.
Since the latter has finitely many coordinates, there is some $h\in A$
whose coordinate takes both values $0$ and $1$ arbitrarily close to
$\xi$ in $D^\circ$.  The value-$0$ set is open in $D^\circ$, and every one of its
points has a neighborhood on which $h$ is the identity.  Hence
$\Int\Fix(h)$ accumulates on $\xi$.  On the value-$1$ set we have
$hx=ax\neq x$, so $h$ is non-trivial on every neighborhood of $\xi$.
Finally, $h\xi=\xi$ because $a\xi=\xi$.  Thus $h$ represents a non-trivial
germ at $\xi$ whose fixed-point interior accumulates on $\xi$, and $\xi$ is
a non-Hausdorff singularity.

The $G$- and $\Dinf$-orbits coincide.  Since the non-Hausdorff locus is
$G$-invariant, the preceding paragraph gives
\[
    \Sigma_{\Dinf}(p)
    \subseteq
    \operatorname{NH}(G\curvearrowright X).
\]

For the reverse inclusion, let $x\notin\Sigma_{\Dinf}(p)$ and take
$g\in G_x$.  Write $g=s_m\cdots s_1$ with $s_i\in A\cup B$, and put
$x_i=s_i\cdots s_1x$.  Every $x_i$ is $p$-regular.  At such a point, each
fragment $s_i$ agrees on a neighborhood with one fixed element of
$\{1,a\}$ or $\{1,b\}$, according to whether $s_i\in A$ or $s_i\in B$.
Indeed, the relevant profile coordinate is locally constant there, and the
corresponding equality, first valid on the dense set $D^\circ$, extends by
continuity.  There are only finitely many points along this trajectory.  We can
therefore choose one neighborhood $U$ of $x$ on which all these local
choices remain valid.  Their product gives an element $d\in\Dinf$
such that
\[
    g|_U=d|_U.
\]
If $g$ has a non-trivial germ at $x$, then $d\neq1$.  Were
$\Int\Fix(g)$ to accumulate at $x$, it would meet $U$ in a non-empty
open set.  On this set $g=d$, so $d$ would fix an open set, contrary
to the first paragraph.  Thus $x$ is not a non-Hausdorff singularity,
which proves the reverse inclusion in (ii).

Finally, let $x$ be a Hausdorff singularity.  By (ii), $x\in D_p$, so $p$
is defined at every vertex of its $\Dinf$-orbit.  Consider a vertex
$y\in\Fix(a)$ in that orbit and $h\in A$.  The extended coordinate
$\varepsilon_h$ is locally constant at $y$, so there are a neighborhood
$U$ of $y$ and $e\in\{0,1\}$ such that $\varepsilon_h=e$ on
$U\cap D^\circ$.  Hence $h=a^e$ on $U\cap D^\circ$, and by density and
continuity the equality holds throughout $U$.  Since $\Fix(a)$ has
empty interior, the germs of $1$ and $a$ at $y$ are distinct.  Thus
$\varepsilon_h(y)=0$ precisely when $h$ has the identity germ at $y$,
whereas $\varepsilon_h(y)=1$ precisely when it has the germ of $a$.
The same argument applies to $k\in B$ at a vertex fixed by $b$.
At vertices fixed by neither reflection, the profile is already read
from the labeled orbital graph.  At an $a$-fixed vertex, however, every
$h\in A$ gives a loop regardless of its germ, and similarly for $B$ at
a $b$-fixed vertex.  This proves the Hausdorff-singular assertion in (i).
\end{proof}

\begin{rem}[Non-Hausdorff versus purely non-Hausdorff]
\label{rem:coding-vs-PNH}
The equality in Proposition~\ref{pro:fragmentation-coding-general} concerns
non-Hausdorff singularities.  Pure non-Hausdorffness is stronger: it requires
the fixed-point interior to accumulate for every non-trivial germ in the
stabilizer, not merely for one.  If $a\xi=\xi$, Lemma~3.5 of
\cite{Nek2018} identifies the germ group as the quotient
\[
    G_\xi/G_{(\xi)}\cong A/A_{(\xi)},
    \qquad A_{(\xi)}=\{h\in A:\xi\in\Int\Fix(h)\}.
\]
Here $G_{(\xi)}$ consists of elements acting identically on a neighborhood
of $\xi$.  Consequently, $\xi$ is purely non-Hausdorff precisely when
\[
    \xi\in\overline{\Int\Fix(h)}
    \qquad(h\in A\setminus A_{(\xi)}).
\]
The condition depends only on the germ represented by $h$.
A natural coding singularity gives this condition for at least one
such germ.  It therefore characterizes the non-Hausdorff locus, not
necessarily the purely non-Hausdorff locus.
\end{rem}

\begin{Def}[Condition $(W)$]\label{def:W}
A partial coding rule $p$ satisfies condition $(W)$ if its coding map
\[
    c_p:D_p\longrightarrow\mathcal A^\Lambda
\]
is injective.  Equivalently, any two distinct $p$-regular points are
distinguished by some coordinate of their coding names:
\[
    x\neq y
    \quad\Longrightarrow\quad
    \exists g\in\Lambda\text{ such that }p(gx)\neq p(gy).
\]
\end{Def}

\begin{pro}[What condition $(W)$ recovers]\label{pro:W-recovers}
Assume $(W)$ for the coding rule $p: D_p\to \mathcal A$, and let $\mathcal Q=\{Q_\beta\}_{\beta\in\mathcal B}$ be any
finite clopen partition of $X$.  Then the $\mathcal Q$-coding name of a
$p$-regular point is determined set-theoretically by its $p$-coding name:
there is a unique equivariant map
\[
    \Theta_{\mathcal Q}:c_p(D_p)\longrightarrow\mathcal B^\Lambda
\]
such that
\[
    \Theta_{\mathcal Q}(c_p(x))(g)=\beta
    \quad\Longleftrightarrow\quad
    gx\in Q_\beta.
\]
For the natural fragmentation coding rule,
Proposition~\ref{pro:fragmentation-coding-general} gives
\[
    D_p=X\setminus\operatorname{NH}(G\curvearrowright X).
\]
Thus the $p$-regular points are the Nekrashevych-regular points together with
the Hausdorff singularities.  At a Nekrashevych-regular point, the
$p$-coding graph is exactly the rooted fragmented orbital graph, so that
orbital graph determines the $\mathcal Q$-coding name at that point.  At a Hausdorff
singularity, the same conclusion holds for the germ orbital
graph; the ordinary orbital graph alone may forget the additional germ
label recorded by the maximally extended coding rule.
\end{pro}

\begin{proof}
Since $c_p$ is injective, its name determines the point $x$ and hence
its $\mathcal Q$-name.  This defines $\Theta_{\mathcal Q}$, and the
definition shows that it is equivariant.  The assertions about the
fragmentation coding follow from
Proposition~\ref{pro:fragmentation-coding-general}.
\end{proof}

\begin{rem}[Global determination versus a local rule]\label{rem:W-not-local}
Proposition~\ref{pro:W-recovers} uses the entire coding name.
Condition $(W)$ need not give a finite $p$-window which determines
the central $\mathcal Q$-symbol, let alone a uniform radius.
In other words, the orbital graph can distinguish all regular points
without determining their clopen neighborhoods from finite rooted balls.
\end{rem}

\begin{defpro}[Condition $(U)$ and its equivalent formulations]
\label{pro:U}
For an action $\Lambda \curvearrowright X$ and a coding rule $p: D_p\to \mathcal A$, choose finite subsets $F_m\nearrow\Lambda$.  For a pattern
$w\in\mathcal A^{F_m}$, put
\[
    C_m(w)=\{x\in D_p:c_p(x)|_{F_m}=w\}.
\]
Following the usual semicocycle terminology, form the \emph{graph
semicocycle extension}
\[
    \widehat X_p
    =\overline{\{(x,c_p(x)):x\in D_p\}}
    \subset X\times\mathcal A^\Lambda, \qquad g\cdot (x,c)=\bigl(g\cdot x, h\mapsto c(hg)\bigr).
\]
and the associated \emph{symbolic semicocycle extension}, or
\emph{semicocycle hull},
\[
    Z_p=\overline{c_p(D_p)}\subset\mathcal A^\Lambda.
\]
Let $\rho_p:\widehat X_p\to Z_p$ be the symbolic projection, and for
$\alpha\in\mathcal A$ put
\[
    \widehat P_\alpha
    =\{(x,c)\in\widehat X_p:c(e)=\alpha\}
    =\overline{\{(x,c_p(x)):p(x)=\alpha\}}.
\]
These are the natural resolved closures of the semicocycle pieces.

We say that the coding rule $p$ satisfies condition $(U)$ if any, and hence all, of the
following equivalent conditions hold:
\begin{enumerate}[label=(\roman*)]
    \item the finite-name cells shrink uniformly,
    \[
        \lim_{m\to\infty}
        \sup_{C_m(w)\neq\varnothing}\diam C_m(w)=0;
    \]
    \item $\rho_p:\widehat X_p\to Z_p$ is injective;
    \item the clopen partition
    $\widehat{\mathcal P}=\{\widehat P_\alpha\}_{\alpha\in\mathcal A}$
    is generating for the action on $\widehat X_p$;
    \item for every finite clopen partition $\mathcal Q$ of $X$, the
    $\mathcal Q$-coding on $D_p$ is locally derivable from the coding $c_p$.
\end{enumerate}
When a generating coding $\chi_X$ exists for $\Lambda \curvearrowright X$, these conditions are also
equivalent to:
\begin{enumerate}[label=(\roman*),resume]
    \item $\chi_X$ is locally derivable from the coding $c_p$ on $D_p$.
\end{enumerate}
Moreover, $(U)$ implies $(W)$.
\end{defpro}

\begin{proof}
Assume (i) and suppose $(x,c),(y,c)\in\widehat X_p$ have the same symbolic
coordinate.  Choose regular graph points
$(x_n,c_p(x_n))\to(x,c)$ and
$(y_n,c_p(y_n))\to(y,c)$.  For every fixed $m$, the two coding names agree on coordinates
$F_m$ for all sufficiently large $n$.  Hence $x_n,y_n$ lie in one common
cell $C_m(w)$, and (i) gives $x=y$.  Thus (ii) holds.

Conversely, if (i) fails, there are $\epsilon>0$, integers $m_j\to\infty$,
and points $x_j,y_j\in D_p$ with
\[
    c_p(x_j)|_{F_{m_j}}=c_p(y_j)|_{F_{m_j}},
    \qquad d(x_j,y_j)\geq\epsilon.
\]
After passing to subsequences, the corresponding marked graph points converge
to $(x,c)$ and $(y,c)$ with $x\neq y$ and the same symbolic coordinate.
Hence (ii) fails.  This proves (i)$\Leftrightarrow$(ii).

The complete $\widehat{\mathcal P}$-name of $(x,c)$ is exactly the symbolic
coordinate $c$, so (ii) and (iii) are equivalent.

Assume (i), and let $\mathcal Q$ be finite and clopen.  Distinct atoms have
a positive minimum distance.  For sufficiently large $m$, every non-empty
cell $C_m(w)$ has smaller diameter and is therefore contained in one
$\mathcal Q$-atom.  The central $\mathcal Q$-symbol is then determined by one
finite $p$-window, and equivariance gives the complete local derivation.
Thus (i) implies (iv).

Conversely, choose for every $\epsilon>0$ a finite clopen partition of $X$
whose atoms have diameter less than $\epsilon$.  By (iv), its central symbol
is determined by a finite $p$-window.  Every sufficiently fine cell $C_m(w)$
is therefore contained in one atom and has diameter less than $\epsilon$.
Thus (iv) implies (i).

To prove (v)$\Rightarrow$(i), fix $\epsilon>0$.  Since $\chi_X$ is a
homeomorphism onto its image, agreement on a sufficiently large finite
$\chi_X$-window forces two points of $X$ to be $\epsilon$-close.
Local derivability from $c_p$ determines this window from a finite
$p$-window.  Thus, for all sufficiently large $m$, every $C_m(w)$ has
small diameter, as required in (i).  The converse follows from (iv),
which applies to every finite clopen generating partition.  This also
shows that the choice of generating coding in (v) does not matter.

Finally, if two points have the same complete $p$-name, then they belong to
the same cell $C_m(w)$ for every $m$.  Condition (i) forces them to coincide,
so $(U)$ implies $(W)$.
\end{proof}

\begin{rem}[Schreier-graph hull and amenability]\label{rem:matte-bon-hull}
Let $p$ be the natural coding rule associated with a finite fragmentation
$G\curvearrowright X$ of a minimal action $\Dinf\curvearrowright X$ on an
infinite Cantor set.  If $Z_p$ is finite, its cell sequences are periodic.
Let $d\geq1$ be a period of the associated one-edge labeled regular orbital
line.  Every fragment generator commutes with the period translation
$\tau_d(n)=n+d$ on that line.  Since the action on a dense regular orbit is
faithful, we obtain an embedding
\[
 G\hookrightarrow C_{\operatorname{Sym}(\Z)}(\tau_d)
   \cong \Z^d\rtimes\operatorname{Sym}(d),
\]
which is virtually abelian and hence amenable.  Here the symmetric group
permutes the $d$ residue classes, and each integer coordinate records a
translation within one residue class.

Now assume that $Z_p$ is infinite.  Then $Z_p$ may be viewed as the hull of
the labeled fragmented regular orbital line.  After passing to
the equivalent one-edge presentation, the fragment labels define clopen
cylinder sets, and every fragment acts by moving the origin one step along
the line, or fixing it, according to the label at the origin.  Hence the
fragmentation embeds in the topological full group of the one-step shift on
this line hull.  This is the Schreier-graph subshift construction used by
Matte Bon in \cite[Section~3, especially Proposition~3.7]{MatteBon15}; for fragmentations, the same argument is recalled in
\cite[Proposition~5.1]{Nek2018}.  By the theorem of Juschenko--Monod that the
topological full group of a minimal Cantor $\Z$-system is amenable
\cite[Theorem~A]{JM13}, this also gives a proof of amenability of the
fragmentation group.
\end{rem}

Suppose that $p$ satisfies $(U)$, and let $\chi$ be a finite clopen
coding.  We recover $c_p\vee\chi$ locally from $c_p$ by keeping its
central symbol and reading the central $\chi$-symbol from a bounded
$c_p$-window; see Lemma~\ref{lem:join-MLD}.  Now restrict the original
fragment generators to atoms determined by finitely many
$\chi$-coordinates.  The same window information determines these
restrictions, so the new orbital graph is locally derivable from the
original fragmentation coding.  Hence LR passes to this further
fragmentation.

With only $(W)$, the coding distinguishes points on the dense regular
graph, but distinct points added in its closure may have the same name.
Thus injectivity on the regular domain need not give a generating
partition on the whole extension.

\begin{rem}\label{rem:germs-vs-MLD}
Groups of germs identify maps which agree on a neighborhood.  They do
not describe how the fragment labels are distributed along an orbit,
and so do not by themselves control linear repetitivity.
Condition $(U)$ allows us to determine clopen membership from finite
fragmentation windows; condition $(W)$ only separates points by their
complete names.  In Theorem~\ref{thm:main-transfer-strong}, $(W)$ holds
but $(U)$ fails, and the joint $c_p\vee \chi_X$ is not LR.
\end{rem}

\subsection{Containing fragmentations of alternating full groups and their orbital graphs}\label{subsec:full-groups}

Let $\Lambda\curvearrowright X$ be a minimal action on an infinite Cantor set, let
$\mathcal H$ be its groupoid of germs, and write
\[
    \Aalt_{\mathcal H}=\Aalt(\mathcal H).
\]
For $\alpha\in\operatorname{Homeo}(X)$, define its \emph{movement set} by
\[
    \Mov(\alpha)
    :=\{x\in X:\alpha(x)\neq x\}
    =X\setminus\Fix(\alpha).
\]
This set is open, and its closure is the usual support of $\alpha$.
We use movement sets because the orbital graph determines whether
$\alpha$ moves a given point of $X$.

\begin{pro}[Local three-cycles for a minimal action]
\label{pro:local-three-cycles}
The following hold.
\begin{enumerate}[label=(\roman*)]
    \item If $\gamma\in\mathcal H$ has source $x$, range $y$, and
    $x\neq y$, then there exists $\alpha\in\Aalt_{\mathcal H}$ whose germ
    at $x$ is $\gamma$.
    \item For every non-empty clopen set $U\subset X$ and every $x\in U$,
    there exists $\alpha\in\Aalt_{\mathcal H}$ such that
    \[
        x\in\Mov(\alpha)\subset U.
    \]
\end{enumerate}
\end{pro}

\begin{proof}
For (i), choose a third point $z$ in the $\Lambda$-orbit of $x$, distinct from
$x$ and $y$.  Take pairwise disjoint clopen neighborhoods of $x,y,z$.
After shrinking the bisections representing the relevant orbit arrows, they
form a degree-three multisection for which one matrix entry contains
$\gamma$.  One of the two corresponding $3$-cycles belongs to
$\Aalt_{\mathcal H}$ and has germ $\gamma$ at $x$.

For (ii), the set $\Lambda x\cap U$ is infinite by minimality.  We can
therefore choose two distinct further points $y,z$ in it and repeat the construction
in (i), this time taking all three diagonal components inside $U$.
The resulting $3$-cycle moves $x$ and is the identity outside $U$.
\end{proof}

\begin{cor}[The alternating full group has the same orbits and germ groupoid]
\label{pro:same-germ-groupoid}
For every minimal action $\Lambda \curvearrowright X$ as above,
\[
    \Aalt_{\mathcal H}x=\Lambda x
    \qquad\text{for every }x\in X,
\]
and the groupoid of germs of the action
$\Aalt_{\mathcal H}\curvearrowright X$ is exactly $\mathcal H$.
\end{cor}

\begin{proof}
Every germ of an element of $\Aalt_{\mathcal H}$ already belongs to
$\mathcal H$.  Proposition~\ref{pro:local-three-cycles}(i) realizes every
arrow with different source and range by a germ of an element of
$\Aalt_{\mathcal H}$, and therefore already gives equality of the orbits.

It remains to consider an isotropy germ $\gamma\in\mathcal H$ at $x$.  Choose an arrow
$\delta:x\to y$ with $y\neq x$.  By
Proposition~\ref{pro:local-three-cycles}(i), both $\delta$ and
$\delta\gamma$ are germs of elements $\alpha,\beta\in\Aalt_{\mathcal H}$.
Then the germ of $\alpha^{-1}\beta$ at $x$ is
\[
    \delta^{-1}(\delta\gamma)=\gamma.
\]
Thus every arrow of $\mathcal H$ occurs as a germ of the alternating-full-
group action.
\end{proof}

\begin{pro}[Clopen sets from movement sets and orbital graphs]
\label{pro:clopen-readable}
For every clopen set $U\subset X$, there exist
$\alpha_1,\ldots,\alpha_r\in\Aalt_{\mathcal H}$ such that
\[
    \Mov(\alpha_i)\subset U
    \qquad\text{and}\qquad
    U=\bigcup_{i=1}^{r}\Mov(\alpha_i).
\]
If $\Aalt_{\mathcal H}$ is finitely generated and $S$ is any finite
generating set, then membership in $U$ is determined by a uniformly bounded
rooted ball in the $S$-labeled orbital graph.
\end{pro}

\begin{proof}
For every $x\in U$, Proposition~\ref{pro:local-three-cycles}(ii) gives an
$\alpha_x\in\Aalt_{\mathcal H}$ with
\[
    x\in\Mov(\alpha_x)\subset U.
\]
These movement sets cover the compact set $U$, so a finite subcover gives
the first assertion.

Suppose now that $S$ is a finite generating set, and write each
$\alpha_i$ as a word $w_i$ in $S$.  To decide whether
$x\in\Mov(\alpha_i)$, follow the path labeled by $w_i$ from $x$
and check whether it returns to $x$.  A ball of radius $|w_i|$ contains
this path.  Taking the maximum of these finitely many lengths lets us
perform every test and hence decide whether $x\in U$.
\end{proof}

\begin{rem}
Proposition~\ref{pro:clopen-readable} explains in what sense a finite
generating set of an alternating full group determines all clopen sets.
It does not generate the Boolean algebra $\Clop(X)$ directly.  Instead,
minimality gives local three-cycles whose movement sets cover a given
clopen set.  Compactness gives a finite subcover, and each three-cycle
is a word in the chosen generators.  This applies to every minimal
Cantor action whose alternating full group is finitely generated.
\end{rem}

\begin{lem}[Elementary restricted fragments]
\label{lem:elementary-restricted-fragment}
Let $h\in A$ and let $E\subset X$ be clopen and $h$-invariant.  Define
\[
    h_E(x)=
    \begin{cases}
       h(x),&x\in E,\\
       x,&x\notin E.
    \end{cases}
\]
Then $h_E$ is an involutive homeomorphism, it is a fragment of $a$, and all
of its germs belong to the germ groupoid of $G$.  The analogous statement
holds for $h\in B$ and the involution $b$.
\end{lem}

\begin{proof}
On $E$ the map $h_E$ agrees with $h$, and on the clopen complement it
is the identity.  Since $E$ is $h$-invariant, $h_E$ is a homeomorphism
of order at most two.  It is a fragment of $a$ because
$h(x)\in\{x,ax\}$.  Its germ is that of $h$ on $E$ and is trivial
outside $E$, so all its germs already belong to the germ groupoid of $G$.
\end{proof}

\begin{lem}[Fragmented routing of compact-open transpositions]
\label{lem:fragmented-routing}
Let \(\mathcal G\) be the groupoid of germs of a finite fragmentation

\[
G=\langle A\cup B\rangle
\]

of a minimal dihedral action. Let \(F\subset\mathcal G\) be a compact-open bisection with

\[
s(F)\cap r(F)=\varnothing,
\]

and let \(\tau_F\) be the full-group involution which exchanges \(s(F)\) and \(r(F)\) through \(F\) and is the identity elsewhere. Then \(\tau_F\) is a finite product of elementary restricted fragments from Lemma~\ref{lem:elementary-restricted-fragment}.

For a finite family of such bisections, the required restricted fragments can be chosen from one finite family. Moreover, after fixing finite clopen local data recording the cocycle choices for the bisections and the restriction domains used in the construction, the labels of these newly added restricted fragments are determined by a uniformly bounded window of the natural fragmentation coding name at the current vertex for the original fragmentation \(G\), together with those clopen data.

\end{lem}

\begin{proof}
Start with a simple path whose successive labels are
$t_1,\ldots,t_m$.  Choose a clopen neighborhood $E$ of its initial
point small enough that the sets $E_i=t_i\cdots t_1E$, with $E_0=E$,
are pairwise disjoint.  Put $\tau_i=(t_i)_{E_{i-1}\sqcup E_i}$.
The transposition of the endpoints is then
\[
 \tau_{p|_E}=\tau_1\cdots\tau_{m-1}\tau_m
                 \tau_{m-1}\cdots\tau_1,
 \qquad p=t_m\cdots t_1.
\]
It acts as $p$ on $E$, as $p^{-1}$ on $p(E)$, and fixes the
complement.  We will use this factorization for the path joining the
two endpoints, and then separately for the isotropy that may remain.

Fix $x\in s(F)$ and a word $w$ representing $F$ near $x$.  Choose a
simple point path from $x$ to $w(x)$, represented by $p\in G$.
The element $k=p^{-1}w$ fixes $x$.  If its germ is trivial, then
$F=p$ near $x$, and the preceding palindrome suffices.

Suppose now that the germ of $k$ is nontrivial.  Choose a simple point
path $v$ from $x$ to the reflection-fixed endpoint $\xi$ of its orbital
ray.  We transport the germ along $v$.  By \cite[Lemma~3.5]{Nek2018},
there is an element
$h\in A$, respectively $h\in B$, according as $a\xi=\xi$ or
$b\xi=\xi$, such that
\[
 [vkv^{-1}]_\xi=[h]_\xi.
\]
Thus the involution $k_0=v^{-1}hv$ fixes $x$ and has the same germ as
$k$, so $F=pk_0$ on a neighborhood of $x$.  Choose a sufficiently
small clopen $k_0$-invariant neighborhood $U\subset s(F)$ on which this
equality holds and whose images along each of the two simple paths
$p$ and $v$ are pairwise disjoint.  Such a neighborhood is obtained
by intersecting a small clopen neighborhood with its $k_0$-image.
In particular $U\cap p(U)=\varnothing$, and $v(U)$ is $h$-invariant.

If $v$ has positive length, let $\tau_v$ be its endpoint
transposition, constructed by the palindrome above, and put
\[
 \kappa=\tau_v\,h_{v(U)}\,\tau_v^{-1}.
\]
Then $\kappa$ is a product of elementary restricted fragments, agrees with
$k_0$ on $U$, and is the identity outside $U$.  If $v$ has length zero,
use $\kappa=h_U$ instead.  This localizes the isotropy at the source.
Let $\tau_p$ be the endpoint transposition for $p|_U$.  Since $\kappa(U)=U$ and $\kappa$ fixes $p(U)$,
\[
             \tau_{F|_U}=\kappa^{-1}\tau_p \kappa.
\]
On $U$ this map is $p\kappa=F$, on $p(U)=F(U)$ it is
$\kappa^{-1}p^{-1}=F^{-1}$, and it fixes the complement.

The same factorization is available over any clopen $V\subset U$:
keep $\kappa$ and use the simple-path transposition for $p|_{\kappa(V)}$ in
place of $\tau_p$.  Its conjugate by $\kappa^{-1}$ is exactly
$\tau_{F|_V}$.  Consequently, compactness of $s(F)$ and a finite
clopen refinement of these neighborhoods give disjoint source pieces
whose associated transpositions multiply to $\tau_F$.  This also
handles a finite family of bisections using only finitely many
restricted fragments.

All restrictions used above have clopen invariant domains.  Choose a
finite clopen coding which determines membership in these domains.
On a regular orbit, its coding name and the original fragment edge
determine each restricted fragment.  Their products are then determined
along uniformly bounded paths, proving locality.
\end{proof}

For the remainder of this section, write $\mathbf F(G)$ for the topological
full group of the action $G\curvearrowright X$.

\begin{thm}[Germ-preserving containment in a finite fragmentation]
\label{thm:germ-preserving-containment}
Let $\Dinf=\langle a,b\rangle$ act minimally on a Cantor set $X$, let
\[
    G=\langle A\cup B\rangle
\]
be a finite fragmentation, and let
\[
    \mathcal G=\operatorname{Germ}(G\curvearrowright X).
\]
Let $L\leq\Aalt(\mathcal G)$ be finitely generated.  Then there exist finite
fragment groups
\[
    \widehat A\supseteq A,
    \qquad
    \widehat B\supseteq B,
\]
fragmenting $a$ and $b$, respectively, such that, for
\[
    H_L=\langle\widehat A\cup\widehat B\rangle,
\]
one has
\[
    G\leq H_L,
    \qquad
    L\leq H_L,
\]
and
\[
    \operatorname{Germ}(H_L\curvearrowright X)=\mathcal G.
\]
In particular,
\[
    H_L\leq\mathbf F(G).
\]

Fix a finite generating set $S_L$ of $L$.  Let $\chi$ be a finite clopen
coding which records the finitely many compact-open cocycle decisions used
by the elements of $S_L$.\footnote{For each $s\in S_L$, fix a finite clopen
partition $X=\bigsqcup_i U_{s,i}$ and words $w_{s,i}$ in $A\cup B$ with
$s|_{U_{s,i}}=w_{s,i}|_{U_{s,i}}$.  Also fix the multisections and the
clopen domains used in their routing in the proof of Lemma~\ref{lem:fragmented-routing}.  Take a common finite
clopen refinement $\mathcal Q$ of all these domains, including the
invariant restriction sets in Lemma~\ref{lem:elementary-restricted-fragment},
and let $\chi=c_{p_{\mathcal Q}}$ be its coding.  The central symbol $\chi(x)(e)=p_{\mathcal Q}(x)$
then specifies the chosen local word and every restriction-domain test;
the fragmentation coding name specifies how that word acts. This auxiliary $\chi$ is not
assumed generating.  Any further finite clopen refinement is allowed.}  On Nekrashevych-regular points,
\[
    c_G\vee \chi
    \longrightarrow
    \Gamma(H_L\curvearrowright X,\cdot)
\]
is a graph-local derivation.  If $L=\Aalt(\mathcal G)$, then the reverse
local derivation also holds, and therefore
\[
    \Gamma(H_L\curvearrowright X,\cdot)
    \overset{\MLD}{\cong}
    c_G\vee \chi.
\]
\end{thm}

\begin{proof}
Choose a finite generating set $S_L$ of $L$.  By the definition of the
alternating full group, every $s\in S_L$ is a finite product of $3$-cycles
associated with compact-open degree-three multisections of $\mathcal G$.
A $3$-cycle is a product of two transpositions between two diagonal
components of the multisection.  Each of these transpositions is of the form
$\tau_F$ for a compact-open bisection with disjoint source and range.
Lemma~\ref{lem:fragmented-routing} therefore writes every element of $S_L$
as a product of finitely many elementary restricted fragments of elements of
$A\cup B$.

Let $\mathcal F_a$ and $\mathcal F_b$ be the finite collections of all
$a$- and $b$-fragments which occur in these factorizations, and put
\[
    \widehat A=\langle A\cup\mathcal F_a\rangle,
    \qquad
    \widehat B=\langle B\cup\mathcal F_b\rangle.
\]
On each two-point orbit of $a$, a fragment acts either as the identity
or as the transposition; it fixes every point fixed by $a$.  Thus any
two such fragments commute and have order at most two.  The same holds
for $b$.  Since we have added only finitely many fragments,
$\widehat A$ and $\widehat B$ are finite elementary abelian $2$-groups.  Since they contain
$A$ and $B$, they still satisfy
\[
    \widehat A x=\{x,ax\},
    \qquad
    \widehat B x=\{x,bx\}.
\]
Thus $H_L$ is a finite fragmentation of the same dihedral action.  The
chosen factorizations give $L\leq H_L$, while $A\subseteq\widehat A$ and
$B\subseteq\widehat B$ give $G\leq H_L$.

Every newly adjoined fragment is an elementary restricted fragment of an
original fragment, so all of its germs lie in $\mathcal G$.  Therefore
\[
    \operatorname{Germ}(H_L\curvearrowright X)\subseteq\mathcal G.
\]
The reverse inclusion follows from $G\leq H_L$, proving equality of the germ
groupoids.  Hence every element of $H_L$ is locally equal to elements of
$G$, and therefore $H_L\leq\mathbf F(G)$.

For locality, apply the last assertion of
Lemma~\ref{lem:fragmented-routing} to the multisections used in the
chosen expressions of $S_L$.  The coding $\chi$ determines the clopen
pieces used in that construction, while $c_G$ determines the fragment
moves on those pieces.  Together they determine every generator label
of $H_L$ from a bounded window.

Assume finally that $L=\Aalt(\mathcal G)$.  Since $G\leq H_L$, the original
fragmentation coding $c_G$ is read directly from the $H_L$-orbital graph.
Since $L\leq H_L$, every element of the fixed finite generating set $S_L$ is
a word of uniformly bounded length in the fragment generators of $H_L$.
Proposition~\ref{pro:clopen-readable} then recovers every atom of the finite
clopen coding $\chi$ from a bounded ball in the $H_L$-orbital graph.  Thus the
joint coding $c_G\vee \chi$ is locally derivable from the $H_L$-graph, proving
MLD.
\end{proof}

\begin{pro}[The full-group graph in an expansive presentation]
\label{pro:expansive-full-group-MLD}
Assume the dihedral action in
Theorem~\ref{thm:germ-preserving-containment} is expansive, and put
$L=\Aalt(G,X)$.  Let $S$ be any finite generating set of $L$.
For the chosen generating coding $\chi_X$, the
containing fragmentation $H$ can be chosen so that, on every regular
orbit and with its common vertex identification,
\[
 \Gamma(L\curvearrowright X,\cdot;S)
 \overset{\MLD}{\cong}c_G\vee \chi_X
 \overset{\MLD}{\cong}\Gamma(H\curvearrowright X,\cdot).
\]

\end{pro}
\begin{proof}
The group $L$ is finitely generated by
\cite[Theorem~2.7 and Proposition~5.2]{Nek2018}.  For $s\in S$, choose
a finite clopen partition $X=\bigsqcup_i U_{s,i}$ and bounded words
$w_{s,i}$ in $A\cup B$ with
$s|_{U_{s,i}}=w_{s,i}|_{U_{s,i}}$.  By Proposition~\ref{pro:continuous}, a finite $\chi_X$-window tells us
which piece $U_{s,i}$ contains the point.  We then follow the word
$w_{s,i}$ in the graph given by $c_G$.  This determines every $S$-edge
from a finite window of $c_G\vee\chi_X$.

Conversely, Proposition~\ref{pro:clopen-readable} reads every $\chi_X$-atom
from a bounded $S$-ball.  Fix $h\in A\cup B$.  By
Corollary~\ref{pro:same-germ-groupoid}, for each $x\in X$ the germ of
$h$ at $x$ is represented by an element of $L$.  Compactness and a
clopen refinement give
\[
 X=\bigsqcup_{j=1}^{m_h}V_{h,j},\qquad
 h|_{V_{h,j}}=\alpha_{h,j}|_{V_{h,j}},\qquad\alpha_{h,j}\in L.
\]
A bounded $S$-ball determines which $V_{h,j}$ contains $x$, by
Proposition~\ref{pro:clopen-readable}.  Following the corresponding
$S$-word $\alpha_{h,j}$ then locates $hx$, including the case $hx=x$.
Repeating this for the finitely many original fragments recovers $c_G$.  On a regular orbit, the two neighboring
dihedral moves are then selected from $A$ and $B$ through
$Ax=\{x,ax\}$ and $Bx=\{x,bx\}$.  This gives the reverse local rule
on the same orbit, with uniform bounds in both directions.
\end{proof}

\begin{rem}[The general expansive-action comparison]
\label{rem:general-full-group-MLD}
Let $\Lambda\curvearrowright X$ be a minimal expansive action of any
finitely generated group on an infinite Cantor set, and let
$L=\Aalt(\Lambda,X)$ be the alternating full group of its germ groupoid.
Fix a finite symmetric generating set $T$ of $\Lambda$, a finite generating
set $S$ of $L$ (available by \cite[Theorem~5.6]{NekDyn}), and the finite
clopen generating partition $\mathcal Q_X$ with coding $\chi_X$ from
Convention~\ref{conv:expansive-presentations}.  In the bounded-degree graph
sense of Definition~\ref{def:local-derivation}, write
$\Gamma(\Lambda\curvearrowright X,x;T)\vee\chi_X$ for the $T$-labeled
orbital graph whose vertex $y$ is marked by the central symbol
$\chi_X(y)(1_\Lambda)$, namely the $\mathcal Q_X$-atom containing $y$.
Then, for every $x\in X$,
\[
    \Gamma(L\curvearrowright X,x;S)
    \overset{\MLD}{\cong}
    \Gamma(\Lambda\curvearrowright X,x;T)\vee\chi_X.
\]
The vertex identification is the identity on $Lx=\Lambda x$ by
Corollary~\ref{pro:same-germ-groupoid}.  Indeed, each $s\in S$ agrees on
finitely many clopen pieces with fixed $T$-words.  The generating
$\chi_X$-name reads those pieces in a bounded window by
Proposition~\ref{pro:continuous}, so the decorated $\Lambda$-graph locally
recovers every $S$-edge.  Conversely, each $t\in T$ is locally represented
by elements of $L$; compactness gives a finite clopen partition on whose
pieces it agrees with fixed $S$-words.  Proposition~\ref{pro:clopen-readable}
reads those pieces and the $\mathcal Q_X$-color from a bounded $S$-ball.
The fixed words locate the recovered endpoints by bounded paths in both
directions.  This extends the proof of
Proposition~\ref{pro:expansive-full-group-MLD}; the comparison holds on
every orbit because the clopen coding is defined on all of $X$.
\end{rem}

\begin{pro}[Virtuality of all finite containing fragmentations]
\label{pro:virtuality}
Let $G\curvearrowright X$ be a minimal finite fragmentation, and write
$\Aalt(G)$ for the alternating full group of its germ groupoid.  Let $H$
be any finite fragmentation satisfying
\[
       \Aalt(G)\leq H\leq\mathbf F(G).
\]
Such an $H$ exists precisely in the expansive case.
Then
\[
       [H:\Aalt(G)]<\infty.
\]
Consequently every two such finite containing fragmentations are
commensurable through the common finite-index subgroup $\Aalt(G)$, and all
of them have the same growth type as $\Aalt(G)$.
\end{pro}

\smallskip
\noindent\emph{Existence of the containing fragmentation and expansivity.}
Such a finite containing fragmentation exists if and only if
$\Aalt(G,X)$ is finitely generated, or equivalently, if
$G\curvearrowright X$ is expansive.  Expansivity implies finite generation by
\cite[Proposition~5.5 and Theorem~5.6]{NekDyn}, and then
Theorem~\ref{thm:germ-preserving-containment} supplies $H$.  Conversely,
the finite index proved below makes $\Aalt(G,X)$ finitely generated if
$H$ exists.  To see that finite generation forces expansivity, choose a
nontrivial finite clopen partition recording the local $G$-word decisions
of a finite symmetric generating set of $\Aalt(G,X)$.  Equal
$G$-coding names force equal partition coding names under every alternating-full-group
word.  A local three-cycle can, however, send one of two distinct points
into a different atom while fixing the other, a contradiction.  Thus that
partition generates $G\curvearrowright X$.

\begin{proof}
Let $\mathbf F(G)_t$ be the subgroup of $\mathbf F(G)$ generated by its
torsion elements.  Since $H$ is generated by finitely many fragment
involutions,
\[
       H\leq\mathbf F(G)_t.
\]
Put $H'=[H,H]$.  Because $\Aalt(G)\leq H$, the characteristic subgroup
$H'$ is normalized by $\Aalt(G)$.  Moreover $H'\neq1$, since otherwise
$H$ would be abelian while it contains the non-abelian group $\Aalt(G)$.
The simplicity theorem \cite[Theorem~4.1]{NekDyn} therefore gives
\[
       \Aalt(G)\leq H'.
\]
Since $H\leq\mathbf F(G)_t$, it follows that
\[
       \Aalt(G)\leq H'\leq\mathbf F(G)_t'.
\]
By \cite[Theorem~2.7]{MatteBon16}, the group
$\mathbf F(G)_t'$ is simple.  The group $\Aalt(G)$ is normal in
$\mathbf F(G)$ by \cite[Theorem~4.1]{NekDyn}, and therefore it is a
nontrivial normal subgroup of $\mathbf F(G)_t'$.  Consequently
\[
       \mathbf F(G)_t'=\Aalt(G).
\]
The preceding inclusions then give $H'=\Aalt(G)$.  Thus
$H/\Aalt(G)$ is abelian.  It is generated by the images of finitely many
involutions, so it is a finite elementary abelian $2$-group.  This proves the
finite-index assertion and the growth comparison.
\end{proof}

\begin{rem}[The line associated with a containing fragmentation]\label{rem:linear-model}
The orbital graph of the alternating full group may have bounded jumps
and short cycles.  In the expansive case, a containing fragmentation
gives a decorated alternating line instead.  The containing group has
the same word-growth type as the alternating full group by
Proposition~\ref{pro:virtuality}; their orbital graphs are MLD by
Proposition~\ref{pro:expansive-full-group-MLD}.
This genuine line structure is needed in
Proposition~\ref{pro:entropy-growth}.
\end{rem}

\subsection*{Proof of Theorem~\ref{thm:compatibility}}
\phantomsection\label{proof:main-II}
\begin{proof}
Apply Theorem~\ref{thm:germ-preserving-containment} with
$L=\Aalt(G,X)$.  The resulting containing fragmentation $H$ has an
orbital graph locally derivable from the joint of
the natural fragmentation coding $c_G$ and a finite clopen coding $\chi$.  Under $(U)$, $\chi$ is locally derivable from $c_G$ by
Proposition~\ref{pro:U}, so the joint is also locally derivable from
$c_G$ by Lemma~\ref{lem:join-MLD}.  Hence $H$ has LR regular orbital
graphs.  Since its germ groupoid is unchanged, it has the same PNH
singularities.  The periodicity theorem
\cite[Theorem~1.1]{Nek2018} and the growth theorem
\cite[Theorem~6.6]{Nek2018} now give torsion and subexponential growth of $H$, and therefore $L$.
If the group had polynomial growth, Gromov's theorem  \cite{Gromov81} would make it virtually
nilpotent; a finitely generated virtually nilpotent torsion group is finite,
contradicting infinitude.

\end{proof}

\section{Entropy of the orbital graph as a \texorpdfstring{$\Z$}{Z}-subshift}
\label{sec:orbital-entropy}
\subsection{Entropy of codings and semicocycle extensions}\label{subsec:entropy}
Throughout this section $T=ab$, and logarithms are natural.  We normalize
entropy by translation steps; the entropy below is not an invariant of
the abstract fragmented group.

\begin{rem}[Translation versus dihedral normalization]\label{rem:entropy-normalization}
The subgroup $\langle T\rangle\cong\Z$ has index two in $\Dinf$, and
\[
 h_{\rm top}(\Dinf\curvearrowright X)=\tfrac12 h_{\rm top}(T\curvearrowright X)
\]
for the usual amenable-group normalization.  Indeed, the two-coset sets
$F_n=\{T^j,bT^j:0\leq j<n\}$ are F\o lner sets of cardinality $2n$.
For a finite clopen partition $\mathcal Q$, its refinement over $F_n$ is
$\bigvee_{j=0}^{n-1}T^{-j}(\mathcal Q\vee b^{-1}\mathcal Q)$.
Dividing its log-cardinality by $2n$ and taking the supremum over
$\mathcal Q$ gives the formula.  One $A$--$B$ cell corresponds to one
translation step, so cell entropy uses $h_{\rm top}(T)$.
\end{rem}

\begin{Def}[Orbital cell entropy]\label{def:orbital-cell-entropy}
For the orbital graph of a fragmentation $G=\langle A\cup B\rangle$ at a regular point, group one $A$-edge and the following
$B$-edge into a cell recording their crossing profiles; the
remaining fragment labels are loops.  If $C_G(n)$ counts length-$n$ cell
words, allowing both orientations, put
\[
 h_{\rm orb}(G\curvearrowright X)=\lim_{n\to\infty}\frac1n\log C_G(n).
\]
The limit exists because $C_G(n+m)\leq C_G(n)C_G(m)$.  The mutual local
containment of Lemma~\ref{lem:LR-one-all} makes $C_G(n)$, and hence the
entropy, independent of the regular orbit.  Changing the root or orientation
does not change its value.  LR gives $C_G(n)=O(n)$ and hence zero cell
entropy, since all length-$n$ patterns occur in an interval of length $O(n)$.
\end{Def}

\begin{pro}[Entropy of a subshift and of a clopen coding]\label{pro:clopen-coding-entropy}
For a subshift $Z\subseteq\mathcal C^{\Z}$, let $C_Z(n)$ be its
language complexity, the number of length-$n$ words occurring in $Z$.
Then
\[
 h_{\rm top}(\sigma\curvearrowright Z)=\lim_{n\to\infty}\frac1n\log C_Z(n).
\]
For a finite clopen partition $\mathcal Q$ of a Cantor $T$-system,
with central-symbol map $p_{\mathcal Q}$ and associated $T$-coding $\chi=c_{p_{\mathcal Q}}$, let $Z_\chi=\chi(X)$.  Then
\[
 h_{\rm top}(\sigma\curvearrowright Z_\chi)=h_{\rm top}(T,\mathcal Q)
 \leq h_{\rm top}(T\curvearrowright X),\qquad
 h_{\rm top}(T\curvearrowright X)=\sup_{\mathcal Q\text{ finite clopen}}
 h_{\rm top}(T,\mathcal Q).
\]
The inequality is equality for a generating partition.  In particular,
for a minimal dihedral system the positive-entropy condition may be written
\[
 h_{\rm top}(\Dinf\curvearrowright X)
 =\tfrac12 h_{\rm top}(T\curvearrowright X)>0.
\]
A partition generating for $\Dinf$ need not generate for $T$, but
$\mathcal Q\vee b^{-1}\mathcal Q$ does.  For the chosen
$b$-invariant generating partition $\mathcal Q_X$, however, the partition
itself generates for $T$, and $w_X(X)$ is already the required ordinary
one-track subshift.
\end{pro}
\begin{proof}
The nonempty atoms of $\bigvee_{j=0}^{n-1}T^{-j}\mathcal Q$ are in
bijection with the length-$n$ words in $Z_\chi$.  Thus their number has
exponential growth rate $h_{\rm top}(T,\mathcal Q)$.  In a subshift, every finite open cover is
refined by a cylinder partition on some $[-r,r]$.  Its first $n$ iterates
are refined by length-$(n+2r)$ cylinders, giving the same exponential rate.
This proves the subshift formula.  A generating name map is a conjugacy,
and every finite open cover of a Cantor set has a clopen partition refining
it, proving the remaining assertions.  The two-track statement follows
from $\Dinf=\{T^n:n\in\Z\}\sqcup\{bT^n:n\in\Z\}$.
\end{proof}

\begin{lem}[The free locus is universally conull]\label{lem:free-universally-conull}
For a minimal dihedral action on an infinite Cantor set and any finite
fragmentation $G$, the set
\[
 R=\operatorname{Free}(\Dinf\curvearrowright X)=\operatorname{Reg}(G\curvearrowright X)
\]
is an invariant dense $G_\delta$, is contained in $D_p$ for the natural fragmentation coding rule $p$, and has full measure
for every $T$-invariant Borel probability measure.
\end{lem}
\begin{proof}
A point fixed by a nontrivial power of $T$ would have a finite dihedral
orbit.  For a reflection $r$, the sets $T^n\Fix(r)$ are pairwise disjoint.
Indeed, a point in two of them would be fixed by two distinct conjugate
reflections, hence by their product, a nontrivial translation.  For any
$T$-invariant probability measure $\mu$, these disjoint sets all have
the same measure, so $\mu(\Fix(r))=0$.
The countably many reflection fixed sets are nowhere dense by
Proposition~\ref{pro:fragmentation-coding-general}; their complement is
therefore conull and a dense $G_\delta$.
At a free point, each fragment word agrees locally with a dihedral
word.  If the fragment word fixes the point, so does this dihedral
word.  Freeness forces the latter to be the identity, so the fragment
word has trivial germ there.  Conversely, move a
nonfree point along its common orbit to a point fixed by $a$ or $b$.
All corresponding fragments fix it, but they cannot all be locally
trivial: their finitely many identity neighborhoods would give an open
fixed set for the original reflection.  Thus the point is not $G$-regular.
Finally, all fragment profiles are continuous along a free orbit, so
$R\subseteq D_p$.
\end{proof}

For an invertible topological system $S\curvearrowright X$, write $\mathcal M_S(X)$ for the compact convex space of $S$-invariant Borel
probability measures, with the weak-* topology.
\begin{lem}[Injectivity on a universally conull set]\label{lem:universally-null-extension}
Let $f:(Y,S_Y)\to(X,S_X)$ be a continuous factor map of compact metrizable invertible
systems.  Suppose an invariant Borel set $E\subseteq Y$ is conull for
every $S_Y$-invariant measure and $f|_E$ is injective.  Then
\[
 f_*:\mathcal M_{S_Y}(Y)\longrightarrow\mathcal M_{S_X}(X)
\]
is an affine homeomorphism.  Each invariant measure system is isomorphic
via $f$ to its image, ergodic measures correspond, and
\[
 h_\nu(S_Y)=h_{f_*\nu}(S_X),\qquad
 h_{\rm top}(S_Y\curvearrowright Y)=h_{\rm top}(S_X\curvearrowright X).
\]
In particular this holds when every non-singleton fiber lies over a
universally null Borel subset $N$ of $X$, without any bound on fiber sizes.
For an aperiodic $S$, any countable union of $S$-orbits is universally null.
\end{lem}
\begin{proof}
By Lusin--Souslin \cite[Theorem~15.1]{Kechris95}, $f|_E$ has a Borel image
and a Borel inverse.  It is therefore an isomorphism on full-measure
invariant sets.  To lift an invariant measure $\mu$ on $X$, first approximate it by
finite atomic measures, lift these measures to $Y$, and take a weak-*
limit.  This gives a lift of $\mu$.  Averaging its translates and
passing to a further limit gives an invariant lift, since the
pushforward remains $\mu$ throughout the averaging.  Any two invariant lifts of $\mu$ give full measure to $E$.  On $E$,
both are determined by $\mu$ and the Borel inverse of $f$, so they
coincide.  Thus $f_*$ is a continuous affine bijection between compact
Hausdorff spaces and is therefore an affine homeomorphism.  It carries
extreme points to extreme points, which gives the correspondence of
ergodic measures.  Measure entropy is invariant
under isomorphism; the variational principle \cite[Chapter~8]{Walters82}
gives equality of topological entropies.
For the special case, first replace $N$ by its invariant saturation
$\bigcup_{n\in \mathbb{Z}}{S_X}^{-n}(N)$, which is still universally
null, and take $E=f^{-1}(X\setminus N)$.  Finally, an invariant
probability measure for an aperiodic system has no atoms, so every
countable union of orbits has measure zero.
\end{proof}

\begin{pro}[Clopen resolution and fragmentation entropy]\label{pro:resolved-entropy}
Let $\Dinf=\langle a,b\rangle\curvearrowright X$ be a minimal action
on an infinite Cantor set, and let
$G=\langle A\cup B\rangle\curvearrowright X$ be a finite fragmentation
of this action.  Let $p$ be the natural fragmentation coding rule, with projections
\[
 X\xleftarrow{\pi}\widehat X_p\xrightarrow{\rho_p}Z_p
\]
as in Proposition--Definition~\ref{pro:U}, where $\Lambda$ is $D_\infty$ (not $G$).  The graph extension $\pi$ is minimal
and almost one-to-one.  The fragmentation has a unique lift to $\widehat X_p$ such that $\pi$ is an
almost one-to-one extension of $G\curvearrowright \widehat X_p$ to
$G\curvearrowright X$.  The same formulas define a fragmentation
$G\curvearrowright Z_p$.  Furthermore, $G\subset \mathbf F(D_\infty \curvearrowright \widehat X_p)$.  Moreover,
\[
 h_{\rm orb}(G\curvearrowright X)=h_{\rm top}(T\curvearrowright Z_p)
 \leq h_{\rm top}(T\curvearrowright\widehat X_p)
 =h_{\rm top}(T\curvearrowright X).
 \tag{4.1}\label{eq:resolved-entropy}
\]
For the resolved coordinate partition
$\widehat{\mathcal P}=\{\widehat P_\alpha\}_{\alpha\in\mathcal A}$ of
$\widehat X_p$, one has
\[
 h_{\rm top}(T\curvearrowright\widehat X_p,\widehat{\mathcal P})
 =h_{\rm top}(T\curvearrowright Z_p).
\]
This clopen partition of $\widehat X_p$ is topologically generating
exactly under $(U)$.
The map $\pi_*$ is an affine homeomorphism of invariant-measure spaces,
with the measure isomorphisms and entropy preservation of
Lemma~\ref{lem:universally-null-extension}.
\end{pro}
\begin{proof}
Coordinatewise continuity on $D_p$ gives
\[
 \pi^{-1}(x)=\{(x,c_p(x))\}\quad(x\in D_p).
 \tag{4.2}\label{eq:resolved-singleton}
\]
The graph points over $R=\operatorname{Free}(\Dinf\curvearrowright X)$
are dense, since $R$ is dense in $D_p$ and the coding is continuous
there.  By minimality of $X$, any nonempty closed invariant subset of
$\widehat X_p$ projects onto $X$.  It must contain all the singleton
fibers over $R$, and hence all of $\widehat X_p$.  This proves minimality
of the extension.
For
$h\in A$ define its lift on $\widehat X_p$ to be $a(x,c)$ when the
$h$-coordinate of $c(e)$ is $1$, and $(x,c)$ otherwise.  These two domains
are clopen and $a$-invariant, because the corresponding profile is
$a$-invariant on the dense graph.  The analogous $B$-formula uses $b$.
All relations and fragmentation identities hold on the dense graph and
hence everywhere.  Projection is equivariant and makes the lifted action
faithful.  Forgetting $x$ gives the same clopen formulas on $Z_p$.

To compare the coding with the orbital graph, take a free point $x$
and order its dihedral orbit as
\[
 \cdots\xleftrightarrow{b}bT^{n-1}x\xleftrightarrow{a}T^nx
 \xleftrightarrow{b}bT^nx\xleftrightarrow{a}T^{n+1}x
 \xleftrightarrow{b}\cdots.
\]
Its $n$th $A$--$B$ cell is $\zeta_x(n)=p(T^nx)=(\alpha_n,\beta_n)$.
At the other vertices,
\[
 p(bT^nx)=(\alpha_{n+1},\beta_n).
 \tag{4.3}\label{eq:cell-profiles}
\]
Thus the one-track cell word determines the entire dihedral name by a
finite-window rule.  These identities pass to the hulls, giving a conjugacy
between $(Z_p,T)$ and the cell hull.  Both orientations occur by dihedral
minimality; reversal sends $\zeta_x$ to
$n\mapsto(\alpha_{1-n},\beta_{-n})$, which has the same entropy.
Proposition~\ref{pro:clopen-coding-entropy} proves the cell and partition
entropy assertions.  Finally, every invariant measure on $\widehat X_p$ gives full measure
to $\pi^{-1}(R)$ by Lemma~\ref{lem:free-universally-conull}.
The fibers over $R$ are singletons by \eqref{eq:resolved-singleton},
so we may apply Lemma~\ref{lem:universally-null-extension} to $\pi$.
The remaining inequality follows because $\rho_p$ is a factor map.
\end{proof}

\begin{pro}[Adjoining a clopen coding]\label{pro:joint-coding-entropy}
For any finite clopen coding $\chi$ of $X$, put
\[
 Z_{p,\chi}=\overline{\{(c_p(x),\chi(x)):x\in D_p\}}.
\]
Then
\[
 h_{\rm orb}(G\curvearrowright X)\leq h_{\rm top}(T\curvearrowright Z_{p,\chi})
 \leq h_{\rm top}(T\curvearrowright X).
 \tag{4.4}\label{eq:joint-entropy-bounds}
\]
The second inequality is equality when the dihedral $\chi$-coding is
generating, in particular for $\chi=\chi_X$ when this coding exists.
For an arbitrary coding $\chi=c_\ell$, the joint has the
$\Z$-indexed presentation
\[
 n\longmapsto\bigl(p(T^nx),\ell(T^nx),\ell(bT^nx)\bigr),
 \tag{4.5}\label{eq:joint-cell-coding}
\]
whose subshift entropy equals $h_{\rm top}(T\curvearrowright Z_{p,\chi})$.
For the $b$-invariant generating coding $\chi_X=c_{f_X}$ of a minimal
expansive dihedral action, the last two coordinates coincide; hence
this presentation is MLD with the one-track joint
\[
 j_{p,X}(x)(n):=\bigl(p(T^nx),f_X(T^nx)\bigr),
 \tag{4.5'}\label{eq:one-track-joint}
\]
without any change in entropy.
\end{pro}
\begin{proof}
The map $(x,c)\mapsto(c,\chi(x))$ is a continuous surjection from
$\widehat X_p$ onto $Z_{p,\chi}$, and $Z_{p,\chi}$ factors onto $Z_p$.
Proposition~\ref{pro:resolved-entropy} gives the inequalities.  If $\chi$ is
generating, the second coordinate also recovers $X$ continuously, giving
equality.  Formula~\eqref{eq:cell-profiles} recovers the missing
$p$-coordinates.  An arbitrary $\chi$ needs the two entries shown in
\eqref{eq:joint-cell-coding}, whereas for the invariant generator
they agree and give exactly \eqref{eq:one-track-joint}.
\end{proof}

\subsection{Positive entropy and exponential growth}
\label{subsec:entropy-growth}
We first show that $(W)$ preserves invariant measures and entropy.
Together with LR, it implies zero entropy for the joint coding.

\begin{pro}[Invariant measures and entropy under $(W)$]\label{pro:W-entropy}
If the natural fragmentation coding rule $p$ satisfies $(W)$, both arrows induce
affine homeomorphisms
\[
 \mathcal M_T(X)\xleftarrow{\pi_*}\mathcal M_T(\widehat X_p)
 \xrightarrow{(\rho_p)_*}\mathcal M_T(Z_p).
 \tag{4.6}\label{eq:W-measure-simplex}
\]
For each invariant measure $\lambda$ on $\widehat X_p$, the two factor maps
are measure-theoretic isomorphisms and
\[
 h_{\pi_*\lambda}(T)=h_\lambda(T)=h_{(\rho_p)_*\lambda}(T).
\]
In particular, ergodic measures correspond and the entropy functions on
these measure spaces agree under the displayed homeomorphisms.  Also,
\[
 h_{\rm orb}(G\curvearrowright X)=h_{\rm top}(T\curvearrowright Z_p)
 =h_{\rm top}(T\curvearrowright\widehat X_p)=h_{\rm top}(T\curvearrowright X).
 \tag{4.7}\label{eq:W-entropy}
\]
Every finite clopen joint $Z_{p,\chi}$ has this same entropy.  Consequently
$(W)+\mathrm{LR}$ forces zero entropy for $X$ and for every such joint.
\end{pro}
\begin{proof}
Set $E=\pi^{-1}(R)$.  It is invariant, Borel and conull for every invariant
measure on $\widehat X_p$.  Its points are $(x,c_p(x))$ by
\eqref{eq:resolved-singleton}.  Equal symbolic coordinates of two points
of $E$ imply equal first coordinates by $(W)$, since $R\subseteq D_p$.
Thus $\rho_p|_E$ is injective.  Lemma~\ref{lem:universally-null-extension}
applies to $\rho_p$, as it already does to $\pi$ in
Proposition~\ref{pro:resolved-entropy}.  This gives the claimed affine
homeomorphisms, measure-theoretic isomorphisms and entropy equalities.
The variational principle gives \eqref{eq:W-entropy}; the bounds in
\eqref{eq:joint-entropy-bounds} coincide.  LR gives zero cell entropy by
Definition~\ref{def:orbital-cell-entropy}.
\end{proof}

\begin{pro}[Zero entropy via a containing fragmentation]
\label{pro:alternating-orbital-entropy}
Suppose $\Aalt(G,X)$ is finitely generated, the natural fragmentation coding rule satisfies
$(W)$, and the generic $G$-graphs are LR.  There is a germ-preserving finite
fragmentation $H$ containing $G$ and $\Aalt(G,X)$ whose natural generic line
has zero cell entropy.  For any finite generating set $S$ of $\Aalt(G,X)$,
its orbital graph is graph-locally derivable from this line.  If $C_S(r)$
counts its rooted radius-$r$ ball types on regular orbits, then
\[
 \limsup_{r\to\infty}\frac1r\log C_S(r)=0.
\]
\end{pro}
\begin{proof}
Theorem~\ref{thm:germ-preserving-containment}, with $L=\Aalt(G,X)$, gives
$H$ whose generic line is MLD with $c_G\vee \chi$ for a finite clopen $\chi$.
Propositions~\ref{pro:W-entropy} and~\ref{pro:joint-coding-entropy} give zero
entropy for this joint and hence for the $H$-line.  Write each $s\in S$ as a word in the fragment generators of $H$.
Their lengths have a common bound, so an $H$-window of radius $O(r)$
determines every path in a rooted radius-$r$ $S$-ball.  Hence
$C_S(r)\leq C_0C_H(Cr+C)$ for fixed constants.  Since the $H$-line
has zero cell entropy, this bound gives the required limit.
\end{proof}

\begin{rem}[Measured and topological recovery]\label{rem:W-entropy-not-LR}
Under $(W)$, the invariant-measure simplices and the corresponding
measured systems agree, but the compact topologies need not agree.
Under $(U)$, $\rho_p$ is a homeomorphism and every clopen membership
test is determined by a finite window.  Thus $(W)$ with LR gives zero
joint entropy, whereas $(U)$ with LR gives LR of the joint.
Neither the measure correspondence nor zero entropy alone implies
subexponential word growth.  The example in
Subsection~\ref{sec:construction} distinguishes these conclusions.
\end{rem}

\begin{Def}[Admissible orbital lines]
\label{conv:admissible-orbital-lines}
Following \cite[Section~2.2]{NekSub25}, a bi-infinite sequence
$w=(E_i)_{i\in\Z}$ of nonempty subsets of a finite alphabet $S$ is
\emph{admissible} if $E_i\cap E_{i+1}=\varnothing$ for every $i$.
Its natural labeled line has an $s$-edge between $i$ and $i+1$ for each
$s\in E_i$, and an $s$-loop at $i$ when $s\notin E_{i-1}\cup E_i$.

An admissible sequence $w$ defines a continuous action \[
\langle S \mid s^2 = 1 \text{ for all } s \in S \rangle
 \curvearrowright Z_w \] on the compact subshift $Z_w$ generated by $w$.  If $w$ is aperiodic and
uniformly recurrent, then $Z_w$ is a Cantor set.  On the associated integer
line, the corresponding involutions act by
\[
 s(i)=
 \begin{cases}
  i+1,&s\in E_i,\\
  i-1,&s\in E_{i-1},\\
  i,&s\notin E_{i-1}\cup E_i.
 \end{cases}
\]
Write $\Lambda_w\leq\operatorname{Sym}(\Z)$ for the quotient group they generate, and $C_w(n)$ for its
length-$n$ word complexity.  A \emph{traversal} of a pattern
$P=(E_1^P,\ldots,E_n^P)$ is a choice $s_i\in E_i^P$ for each $i$;
its group element is $s_n\cdots s_1$, which carries the initial endpoint
of every occurrence to its other endpoint in exactly $n$ steps.
\end{Def}

\begin{pro}[Ordered traversals and growth of admissible orbital lines]
\label{pro:entropy-growth}
Let $w=(E_i)_{i\in\Z}$ be admissible over $S$, put
$r=\max_i|E_i|<|S|$, and write
$\gamma_{\Lambda_w,S}(n)=|B_{\Lambda_w,S}(n)|$.

We can define an injection
\[
\begin{gathered}
\Phi_n:\{P=(E_{k+1},\ldots,E_{k+n}):k\in\Z\}\longrightarrow B_{\Lambda_w,S}(n)^r,\\
\Phi_n(P)=\bigl(g_1(P),\ldots,g_r(P)\bigr).
\end{gathered}
\]
Therefore
\[
 C_w(n)\leq\gamma_{\Lambda_w,S}(n)^r\leq\gamma_{\Lambda_w,S}(n)^{|S|}
 \qquad(n\geq1).
 \tag{4.8}\label{eq:admissible-growth}
\]
The same estimate holds when words read in both orientations are counted.
Consequently, positive entropy of $Z_w$ forces exponential growth of
$\Lambda_w$.

In particular, if $H=\langle A\cup B\rangle$ is a finite fragmentation
of a (not necessarily minimal) dihedral action, positive entropy of any of
its natural generic orbital lines forces exponential word growth of $H$.
\end{pro}

\begin{proof}
Fix one total order on $S$.  For a nonempty $E\subseteq S$, let
$\alpha^{(k)}(E)$ be its $k$-th largest element when $k\leq|E|$, and
let $\alpha^{(k)}(E)=\min E$ otherwise.  For
$P=(E_1^P,\ldots,E_n^P)$, set
\[
 a_i^{(k)}(P)=\alpha^{(k)}(E_i^P),\qquad
 g_k(P)=a_n^{(k)}(P)\cdots a_1^{(k)}(P)
 \quad(1\leq k\leq r).
\]
By construction, each selected letter crosses the corresponding edge
of $P$.  Thus each of these words traverses $P$ in exactly $n$ steps.

Suppose $g_k(P)=g_k(Q)$ in $\Lambda_w$ for $1\leq k\leq r$.
At an occurrence of $Q$, $g_k(Q)$ moves the initial endpoint $n$
edges to the right.  Since $g_k(P)$ is the same group element, its
word of length $n$ must reach the same endpoint.  Each letter moves
by at most one edge, so there can be no step to the left and no loop:
every letter must cross the next edge to the right.  Reversing the roles of $P,Q$
gives
\[
 a_i^{(k)}(P)\in E_i^Q,\qquad
 a_i^{(k)}(Q)\in E_i^P
 \quad(1\leq i\leq n,\ 1\leq k\leq r).
\]
The first $r$ selected ranks enumerate all elements of each edge set:
\[
 E_i^P=\{a_i^{(k)}(P):1\leq k\leq r\},\qquad
 E_i^Q=\{a_i^{(k)}(Q):1\leq k\leq r\}.
\]
The inclusions above now give $E_i^P=E_i^Q$ for every $i$, and hence
$P=Q$.  Thus $\Phi_n$ is injective, proving
\eqref{eq:admissible-growth}.  For an occurrence read in the opposite
direction, the endpoint is $n$ edges to the left.  Reaching it in
$n$ steps again requires crossing every edge in order, so the same
argument counts words in both orientations.  The entropy formula gives
\[
 \liminf_{n\to\infty}\frac1n\log\gamma_{\Lambda_w,S}(n)
 \geq\frac1r h_{\rm top}(\sigma\curvearrowright Z_w),
\]
proving the positive-entropy assertion.  This argument uses neither
minimality nor repetitivity of $w$.

For the fragmentation $H$, each natural edge label is the set of
generators in $S$ crossing that edge.  These sets are nonempty and
adjacent sets are disjoint, so they form an admissible sequence on any
regular orbital line.  The group $H_w$ defined by any specific line $w$ is the image of
$H$ in its action on that orbit.  Write $C_w^{\mathrm{cell}}(n)$ for the
number of length-$n$ cell words on this line, allowing both orientations.
Each such cell word records $2n$ consecutive crossing sets.  The preceding
bound, allowing both orientations, therefore gives
\[
 C_w^{\mathrm{cell}}(n)\leq |B_{H_w,S}(2n)|^{|A|+|B|}
 \leq |B_{H,S}(2n)|^{|A|+|B|}.
\]
Positive cell entropy implies exponential word growth of $H$.
\end{proof}

\begin{rem}[Relation with the linear-Schreier-graph estimates]
Bartholdi--Nekrashevych--Zheng \cite{BNZ22} obtain upper bounds by
counting the possible portraits of short group elements.  Our estimate
is a lower bound, not a reversal of their portrait estimate.
In \cite[Proposition~4.3]{NekSub25}, prescribed labeled segments give
short elements distinguished by maximal geodesic traversals.  Here we
fix one total order on the generating set and choose finitely many
traversals whose group elements together determine every edge label.
\end{rem}

\subsubsection*{Proof of Theorem~\ref{thm:positive-entropy-criterion}}
\begin{proof}
Put $L=\Aalt(G,X)$ and
$\mathcal G=\operatorname{Germ}(G\curvearrowright X)$.
Expansivity gives finite generation of $L$ by
\cite[Theorem~2.7 and Proposition~5.2]{Nek2018}; simplicity follows from
\cite[Theorem~4.1]{NekDyn}.  By
Corollary~\ref{pro:same-germ-groupoid}, the $L$-orbits coincide with the
infinite $G$-orbits, so $L$ is infinite.

Using $\chi_X$, Theorem~\ref{thm:germ-preserving-containment} and
Proposition~\ref{pro:expansive-full-group-MLD} give a finite fragmentation
$H$ with
\[
 G\leq H,\qquad L\leq H\leq\mathbf F(G),\qquad
 \operatorname{Germ}(H\curvearrowright X)=\mathcal G,
\]
whose natural regular line is MLD with $c_G\vee\chi_X$.

The joint $\bigl(p(T^nx),f_X(T^nx)\bigr)_{n\in\Z}$ has entropy
$h_{\rm top}(ab\curvearrowright X)>0$ by
Proposition~\ref{pro:joint-coding-entropy}.  It is MLD with the $H$-line
using the same dihedral orbit coordinates, so that line also has
positive geometric cell entropy.  Proposition~\ref{pro:entropy-growth}
gives exponential growth of $H$.  Since $[H:L]<\infty$ by
Proposition~\ref{pro:virtuality}, $L$ has exponential growth as well.

Every finite fragmentation of a minimal
dihedral action is amenable by \cite[Proposition~5.1]{Nek2018}; amenability passes to the subgroup $L$.

Assume, in addition, that $G\curvearrowright X$ has a PNH singularity.
This property passes to $H\curvearrowright X$ by equality of the germ groupoids.
The periodicity theorem \cite[Theorem~4.1]{Nek2018} therefore makes $H$
torsion.  The group $L$ is torsion because it is contained in $H$.
\end{proof}

\begin{rem}[Realization of the hypotheses]
\label{rem:positive-entropy-realization}
Subsection~\ref{sec:positive-growth} constructs actions satisfying
the hypotheses of Theorem~\ref{thm:positive-entropy-criterion}.
In particular, Proposition~\ref{pro:positive-entropy-action} produces
a fragmented action satisfying \NekExpHyp\ that preserves the LR regular
orbital graphs of the underlying generalized Grigorchuk boundary action,
while its associated alternating full group has exponential word growth.
\end{rem}

\section{Construction via almost one-to-one extensions}
\label{sec:extensions}
We first prove the lifting lemmas common to the two constructions.
We then give the positive-entropy examples satisfying \NekExpHyp\ in
Subsection~\ref{sec:positive-growth}, and the zero-entropy examples
where LR does not pass to containing orbital graphs in
Subsection~\ref{sec:construction}.

Recall that the singleton-fiber locus of an almost one-to-one factor
map is an invariant dense $G_\delta$.  For minimal systems, its full
preimage is also dense.  This observation gives the following comparison
of regular orbital graphs.

\begin{lem}[Regular orbital graphs in almost one-to-one extensions]
\label{lem:regular-graph-extension}
Let $\pi:Y\to X$ be an almost one-to-one factor map between minimal
Cantor $\Lambda$-systems, where $\Lambda$ is finitely generated.  Let $x$
be regular for $\Lambda\curvearrowright X$.  Then every $y\in\pi^{-1}(x)$ is
regular and $\Lambda_y=\Lambda_x$.  For every finite generating set $S$, restriction
of $\pi$ therefore gives a rooted labeled graph isomorphism
\[
 \Gamma(\Lambda\curvearrowright Y,y;S)
 \cong\Gamma(\Lambda\curvearrowright X,x;S).
\]
No singleton-fiber assumption at $x$ is needed.  If all non-singleton
fibers lie over regular base points, then
\[
 \operatorname{Reg}(\Lambda\curvearrowright Y)
 =\pi^{-1}\bigl(\operatorname{Reg}(\Lambda\curvearrowright X)\bigr).
\]
In particular the regular graph language, its entropy and its repetitivity function
are unchanged under an almost one-to-one extension of minimal systems.
\end{lem}
\begin{proof}
Let $X_1$ be the invariant singleton-fiber locus.  Its full preimage $Y_1$
is dense in $Y$, since it contains the orbit of any point in a singleton
fiber and the action on $Y$ is minimal.  Equivariance gives $\Lambda_y\subseteq\Lambda_x$.  For the reverse
inclusion, take $g\in\Lambda_x$.  Since $x$ is regular, $g$ is the
identity on some open neighborhood $U$ of $x$.  Each point of
$Y_1\cap\pi^{-1}(U)$ is the unique lift of its image, so $g$ fixes it.
This set is dense in $\pi^{-1}(U)$, and continuity therefore makes
$g$ the identity on all of $\pi^{-1}(U)$.  Thus
$\Lambda_x\subseteq\Lambda_y$, and every stabilizer element has trivial
germ at $y$.  The orbit map $\Lambda y\to\Lambda x$ is a bijection;
equivariance shows that it preserves the labels in $S$.

If $x$ is nonregular and $\pi^{-1}(x)=\{y\}$, take $g\in \Lambda_x$ with a
nontrivial germ and moved points $x_n\to x$.  Choose lifts $y_n$ of $x_n$.
Every accumulation point of $(y_n)$ lies over $x$, so compactness and
$\pi^{-1}(x)=\{y\}$ give $y_n\to y$.  By equivariance, $g$ moves
every $y_n$ and fixes $y$.  Its germ at $y$ is therefore nontrivial,
so $y$ is nonregular.  This proves the
claimed equality of regular loci under the additional fiber hypothesis.
For the final assertion, a lift of any regular base point is regular and
has the same graph; Lemma~\ref{lem:LR-one-all} then identifies the regular
graph languages and repetitivity functions of the two minimal systems.
\end{proof}

\begin{lem}[Canonical lifts of a fragmentation]
\label{lem:canonical-fragmentation-lift}
Let $\pi:Y\to X$ be an almost one-to-one factor map between minimal
Cantor $\Dinf$-systems, and let
$G=\langle A\cup B\rangle\curvearrowright X$ be a finite fragmentation
with natural fragmentation coding rule $p$.  If every nonregular base point has a singleton
fiber, the named fragments have unique continuous lifts to a faithful
fragmentation of the same group $G$ on $Y$.  The lifted fragmentation action is minimal,
its regular orbital graphs are unchanged, and every original PNH point
has a singleton fiber whose point is again PNH.  If the dihedral action
on $Y$ is expansive, its lifted fragmentation is expansive as well.

More generally, existence and uniqueness of the lifted fragmentation
hold whenever the lifted original reflections fix their entire fibers
over the corresponding fixed sets:
\[
 a_Yy=y\quad\bigl(\pi(y)\in\Fix(a_X)\bigr),\qquad
 b_Yy=y\quad\bigl(\pi(y)\in\Fix(b_X)\bigr).
 \tag{5.1}\label{eq:fixed-reflection-fibers}
\]
On the inverse image of the base natural-coding domain, the lifted
crossing choices are the pullbacks of the coordinates of $p$.
\end{lem}
\begin{proof}
The singleton-fiber hypothesis implies
\eqref{eq:fixed-reflection-fibers}: a point fixed by an original
reflection is nonregular by Lemma~\ref{lem:free-universally-conull},
and that reflection must fix its unique lift.  We construct the lifts
under the more general condition \eqref{eq:fixed-reflection-fibers}.
For $h\in A$, define
\[
 \widetilde h(y)=
 \begin{cases}
 a_Yy,&h\pi(y)=a_X\pi(y)\ne\pi(y),\\
 y,&\text{otherwise}.
 \end{cases}
\]
The crossing choice is locally constant away from
$\pi^{-1}\Fix(a_X)$, so the lift is continuous there.  At a point
$y$ above $\Fix(a_X)$, both possible values tend to $y$, since
$a_Yy=y$ by \eqref{eq:fixed-reflection-fibers}.  This proves continuity
also at such points.  We lift $B$ by the same formula with $b_Y$.
On the dense singleton-fiber locus these are the unique possible lifts.
All relations of $G$ hold on this dense set, and therefore everywhere.
The lifted action is faithful: since $\pi$ is onto, an element acting
trivially upstairs must act trivially downstairs.  The identities
$\widetilde A y=\{y,a_Yy\}$ and
$\widetilde B y=\{y,b_Yy\}$ follow from the corresponding identities
downstairs and \eqref{eq:fixed-reflection-fibers}.  They give
minimality, since the two actions have the same point-orbits.
Uniqueness of the continuous lifts follows from their values on the
dense singleton-fiber locus.  Expansivity of the lifted fragmentation
follows from \cite[Proposition~5.2]{Nek2018}.

Return to the singleton-fiber hypothesis over nonregular points.
Lemma~\ref{lem:regular-graph-extension} identifies the regular loci and
orbital graphs.

Suppose $\pi^{-1}(x)=\{y\}$.  Equivariance identifies the stabilizers
of $x$ and $y$.  To identify their germ groups, consider a word which
is the identity on a neighborhood of $x$.  Its lift is the identity
on the full preimage, by density of singleton-fiber points.
Conversely, suppose the lift is the identity on a neighborhood $U$
of $y$.  Compactness and $\pi^{-1}(x)=\{y\}$ give a neighborhood $V$
of $x$ with $\pi^{-1}(V)\subset U$.  Equivariance then shows that the
original word is the identity on $V$.

It remains to check the accumulation of fixed-point interiors.  Let
$g$ represent a nontrivial germ at $x$ and suppose that
$\operatorname{Int}\operatorname{Fix}(g)$ accumulates at $x$.  Given a
neighborhood $U$ of $y$, choose a neighborhood $V$ of $x$ such that
$\pi^{-1}(V)\subset U$.  By the accumulation assumption, there is a
nonempty open set
\[
W\subset V\cap\operatorname{Int}\operatorname{Fix}(g).
\]
The lift of $g$ is the identity on the nonempty open set
$\pi^{-1}(W)\subset U$.  Since $U$ was arbitrary, its fixed-point
interior accumulates at $y$.  The germ groups were identified above,
so this argument applies to every nontrivial germ at a PNH point.
Its unique lift is therefore PNH as well.
\end{proof}

\begin{lem}[Compatible fragmentations over free fibers]
\label{lem:compatible-free-graphs}
Let $\pi:Y\to X$ be any factor map between minimal infinite Cantor
$\Dinf$-systems.  Suppose finite fragmentations
$G_X=\langle A_X\cup B_X\rangle$ and
$G_Y=\langle A_Y\cup B_Y\rangle$ have corresponding generator labels $A_X \leftrightarrow A_Y$ and $B_X \leftrightarrow B_Y$
and satisfy $\pi(g_Yy)=g_X\pi(y)$ for every such label.  Then, on the
dense invariant $G_\delta$ set
\[
 Y_0=\pi^{-1}\bigl(\operatorname{Free}(\Dinf\curvearrowright X)\bigr),
\]
restriction of $\pi$ gives canonical rooted labeled graph isomorphisms
\[
 \Gamma(G_Y\curvearrowright Y,y)
 \cong\Gamma(G_X\curvearrowright X,\pi(y)).
\]
In particular their generic regular graph languages and repetitivity
functions agree.  No almost-one-to-one hypothesis is needed.
\end{lem}
\begin{proof}
The free locus of $X$ is a nonempty invariant dense $G_\delta$ by
Lemma~\ref{lem:free-universally-conull}.  Its preimage is nonempty and
invariant, hence dense by minimality of $Y$, and is $G_\delta$.
If $x=\pi(y)$ is free, equivariance gives
$\operatorname{Stab}_{\Dinf}(y)\subseteq
\operatorname{Stab}_{\Dinf}(x)=\{1\}$.  Thus $\pi$ bijects the two
dihedral orbits, identifying both original orbital graphs with
$\operatorname{Cay}(\Dinf,\{a,b\})$.  Fragmentation does not change the orbit sets, and compatibility of
the generators makes the corresponding fragmented edges and loops
agree.  Both points are regular by
Lemma~\ref{lem:free-universally-conull}.  We can then apply
Lemma~\ref{lem:LR-one-all} in each system to obtain the same language
and repetitivity function on all generic regular orbits.
\end{proof}

\begin{rem}[Extending before or after fragmentation]
\label{rem:extension-behavior}
Write $K=\Dinf$.  Let $q:D\to\mathcal C$ be a finite-alphabet coding
rule on a common dense invariant subset
$D\subseteq\operatorname{Reg}(G\curvearrowright X)$, and form the two
marked graph semicocycle extensions
\[
 \widehat X_q^{\Lambda}
 =\overline{\{(x,c_q^{\Lambda}(x)):x\in D\}},\qquad
 c_q^{\Lambda}(x)(g)=q(gx),\qquad \Lambda\in\{G,K\}.
\]
There is a canonical continuous surjection over $X$
\[
 \Theta:\widehat X_q^G\longrightarrow\widehat X_q^K.
\]
Here the $K$-action on $\widehat X_q^G$ is recovered from the lifted
marked fragments: $a$ sends a point $y$ to the unique other point of
$\widetilde A y$ when that set has two elements, and fixes $y$ when it
is a singleton; similarly for $b$ and $\widetilde B$.
These sets have at most two points by density of the common graph
locus.  The recovered maps are continuous, because the fragment sets
are finite and all possible images coalesce at a singleton.  They are
involutions by density.  Reading the central $q$-symbol along this
recovered $K$-action defines $\Theta$.

The fragmentation lifts compatibly to $\widehat X_q^K$ precisely when
$\Theta$ is injective, equivalently when the two marked extensions are
isomorphic over $X$.  Indeed, compatible lifts on $\widehat X_q^K$
allow one to read its $G$-name and produce a continuous inverse to
$\Theta$; both compositions are the identity on the dense common graph
locus.  Conversely, a homeomorphism $\Theta$ transports the canonical
fragmentation on $\widehat X_q^G$.  In particular,
Lemma~\ref{lem:canonical-fragmentation-lift} gives a sufficient
condition for the two constructions to agree.
\end{rem}

In both Subsections~\ref{sec:positive-growth} and~\ref{sec:construction},
we start with the LR boundary fragmentation $G_\omega\curvearrowright X$,
where $X=\partial T_2\cong\Ztwo$, and its PNH orbit $\mathbb{Z}$.
We choose an expansive almost one-to-one Toeplitz extension
$\pi:\widehat X\to X$ of the dihedral adding machine action and lift
the fragmentation.  Then $\pi^*\Clop(X)\subseteq\Clop(\widehat X)$,
with strict inclusion when a fiber is nontrivial.  Thus the extension
adds clopen sets without changing the regular fragmented graphs.
In the following diagram, horizontal arrows indicate fragmentation
and vertical arrows are the compatible factor maps.
\[
\begin{tikzcd}[column sep=5em,row sep=3em]
 G_\omega\curvearrowright\widehat X \arrow[d,"\pi"'] &
 \begin{gathered}\Dinf\curvearrowright\widehat X\\[-.2em]
       \text{\small expansive Toeplitz system}\end{gathered}
 \arrow[l,"\text{lifted fragmentation}"']
 \arrow[d,"\pi\;\text{almost one-to-one}"]\\
 \begin{gathered}G_\omega\curvearrowright X\\[-.2em]
       \text{\small LR regular graphs; PNH}\end{gathered} &
 \begin{gathered}\Dinf\curvearrowright X\\[-.2em]
       X\cong\partial T_2\cong\Ztwo\end{gathered}
 \arrow[l,"\text{boundary fragmentation}"']
\end{tikzcd}
\]
\begin{samepage}
For the generating coding $\chi_{\widehat X}$,
Proposition~\ref{pro:expansive-full-group-MLD}
identifies the alternating-full-group graph and a containing fragmented
line with the joint of the lifted natural coding and this new clopen
coding:
\[
\begin{tikzpicture}[
  >=stealth,
  every node/.style={font=\small,inner sep=2pt,outer sep=1.5pt},
  baseline=(joint.base)
]
  % Common rectangular core; all four nodes refer to the same regular orbit.
  \node (joint) at (0,0) {$c_{G_\omega}\vee\chi_{\widehat X}$};
  \node (chi) at (4.7,0) {$\chi_{\widehat X}$};
  \node (cg) at (0,-2.8) {$c_{G_\omega}$};
  \node (cay) at (4.7,-2.8) {$\operatorname{Cay}(\Dinf)$};
  \node[anchor=north,inner sep=1pt] (chi-description)
    at ([yshift=-1pt]chi.south)
    {($\mathcal Q_{\widehat X}$-decorated $\Dinf$-line)};
  \node[anchor=north,inner sep=1pt] (cg-description)
    at ([yshift=-1pt]cg.south) {$(\mathrm{LR})$};

  % The five MLD comparisons are branches, not extra vertices of the core.
  \node[anchor=east] (ag) at (-1.8,0)
    {$\Gamma(\Aalt(G_\omega,\widehat X);S)$};
  \node[anchor=east] (g) at (-1.8,-2.8)
    {$\Gamma(G_\omega;S)$};
  \node[anchor=west] (ad) at (6.9,0)
    {$\Gamma(\Aalt(\Dinf,\widehat X);S)$};
  \node[anchor=west] (d) at (6.9,-2.8)
    {$\Gamma(\Dinf)$};
  \node (h) at (0,1.45) {$\Gamma(H)$};

  \draw[<->] (ag.east) -- node[above,inner sep=1.5pt] {$\MLD$} (joint.west);
  \draw[<->] (h.south) -- node[right,inner sep=1.5pt] {$\MLD$} (joint.north);
  \draw[<->] (chi.east) -- node[above,inner sep=1.5pt] {$\MLD$} (ad.west);
  \draw[<->] (g.east) -- node[above,inner sep=1.5pt] {$\MLD$} (cg.west);
  \draw[<->] (cay.east) -- node[above,inner sep=1.5pt] {$\MLD$} (d.west);

  \draw[->] (joint.east) -- (chi.west);
  \draw[->] (joint.south) -- (cg.north);
  \draw[->] (chi-description.south) -- (cay.north);
  \draw[->] (cg.east) -- (cay.west);
\end{tikzpicture}
\]
\end{samepage}
All graphs in this diagram are rooted at the same regular point,
which is suppressed from the notation; each occurrence of $S$ refers
to a separately chosen finite generating set of the indicated group.
The upper-right MLD comparison follows from
Proposition~\ref{pro:expansive-full-group-MLD} applied to the
unfragmented action, whose natural crossing profiles are constant.
The lower-left comparison is the natural fragmentation coding, and the
lower-right one is the identification of a regular dihedral orbital
graph with $\operatorname{Cay}(\Dinf,\{a,b\})$, followed, if necessary,
by a change of finite generators.

By Convention~\ref{conv:expansive-presentations}, the $b$-invariant
generator of $\widehat X$ gives a one-track name $w_{\widehat X}$ which
determines $\chi_{\widehat X}$.  We use full dihedral names in the
diagram so that $c_{G_\omega}$ and the spatial coding have the same
orbit coordinates.  As a sequence, their joint is the natural
fragmentation profile together with $w_{\widehat X}$.
In Subsection~\ref{sec:positive-growth}, we choose this clopen coding
to have positive entropy, giving a containing line of positive entropy
and a group of exponential growth.  In Subsection~\ref{sec:construction},
both coordinates are LR and $(W)$ holds, but their joint repetitivity
is uniformly larger than the prescribed function up to linear rescaling.

The two constructions differ at the singular fibers.
In Subsection~\ref{sec:positive-growth}, every nonregular point has a
singleton fiber, so Lemma~\ref{lem:canonical-fragmentation-lift} applies
and preserves its PNH germs.  In Subsection~\ref{sec:construction}, we
split the singular orbit into two-point fibers.  Both lifts
of the reflection-fixed point remain fixed by the lifted reflection, so
\eqref{eq:fixed-reflection-fibers} still gives the fragmentation;
condition (Q3) in Lemma~\ref{lem:q-construction} separately ensures PNH
at both lifts.

\begin{rem}[The alternating full group of the joint coding hull]
\label{rem:joint-hull-full-groups}
Let $p$ be the natural fragmentation coding rule on $\widehat X$, and put
\[
 Y:=\overline{\{(c_{G_\omega}(x),\chi_{\widehat X}(x)):x\in D_p\}}.
\]
Then $Y$ is an expansive $\Dinf$-subshift.  We use $w_{\widehat X}$
for the spatial coordinate and the full dihedral name for the possibly
partial fragmentation coding.  The generating
$\chi_{\widehat X}$-coordinate identifies $Y$ with the graph resolution
of $p$ and gives an almost one-to-one factor $\rho:Y\to\widehat X$.
Propositions~\ref{pro:resolved-entropy} and~\ref{pro:joint-coding-entropy}
therefore give a faithful lifted fragmentation $\widetilde G_\omega\curvearrowright Y$.
Lifting clopen local descriptions through $\rho$ gives a canonical embedding
$\iota:\mathbf F(G_\omega,\widehat X)\hookrightarrow
\mathbf F(\widetilde G_\omega,Y)$, and lifting multisections gives
\[
 \iota\bigl(\Aalt(G_\omega,\widehat X)\bigr)
 \leq\Aalt(\widetilde G_\omega,Y)=\Aalt(\Dinf,Y).
\]
The lift is well-defined because local equalities hold over the dense
singleton-fiber locus and extend by continuity.  The equality on the
right follows because the crossing profiles are clopen on $Y$: every
lifted fragment is locally $1$, $a$, or $b$, and each of $a,b$ locally
agrees with a lifted fragment.  Thus the two actions on $Y$ generate
the same germ groupoid.

The first inclusion can be strict, even though the regular orbital
graphs are MLD.  Indeed, for $x\in\operatorname{Reg}(G_\omega\curvearrowright
\widehat X)$ and its unique lift $y\in Y$, Proposition~\ref{pro:expansive-full-group-MLD}
and its unfragmented case give, with finite generating sets understood,
\[
 \Gamma(\Aalt(G_\omega,\widehat X)\curvearrowright\widehat X,x)
 \overset{\MLD}{\cong}c_{G_\omega}(x)\vee\chi_{\widehat X}(x)
 \overset{\MLD}{\cong}\Gamma(\Aalt(\Dinf,Y)\curvearrowright Y,y).
\]
In fact the first inclusion is equality precisely when $\rho$ is
one-to-one.  If $y_0\ne y_1$ lie over the same point $x$, choose
$g,h\in\Dinf$ with $x,gx,hx$ distinct, and a sufficiently small clopen
neighborhood $U$ of $y_0$ such that $U,gU,hU$ are disjoint and avoid $y_1$.
The three-cycle on these sets belongs to $\Aalt(\Dinf,Y)$, moves $y_0$
to $gy_0$, and fixes $y_1$; it therefore cannot descend through $\rho$.
Conversely, an injective $\rho$ is a conjugacy and identifies the clopen
multisections.  In the PNH examples below, the natural fragmentation coding rule is singular,
so $\rho$ has non-singleton fibers and the inclusion is strict.

We must therefore apply the full-group criterion on the same space.  Since $\widetilde G_\omega\leq\mathbf F(\Dinf,Y)$,
Corollary~\ref{pro:same-germ-groupoid} gives
\[
 \Aalt(\widetilde G_\omega,Y)=\Aalt(\Dinf,Y)
 \quad\Longleftrightarrow\quad
 \Dinf\leq\mathbf F(\widetilde G_\omega,Y).
\]
On this resolved hull the right-hand condition is automatic; it does
not identify $\Aalt(\Dinf,Y)$ with the lifted image of
$\Aalt(G_\omega,\widehat X)$, whose clopen constructions come from
$\widehat X$.  See also Remark~\ref{rem:unfragmented-alternating-full-group}
for the comparison with the unfragmented alternating full group on
the original space.
\end{rem}

\subsection{Positive-entropy constructions and exponential growth}
\label{sec:positive-growth}

\subsubsection*{Positive-entropy extensions preserving LR fragmented orbital graphs}
\label{subsec:positive-LR-extension}

We use the entropy normalization of Subsection~\ref{subsec:entropy}
and the Toeplitz factor theorem recalled in Section~\ref{sec:preliminaries}.

\begin{lem}[A pointed positive-entropy binary Toeplitz system]
\label{lem:self-contained-positive-Toeplitz}
For every $h\in(0,\log2)$ there are a Toeplitz point
$z\in\mathcal A^\Z$ with $\mathcal A=\{0,1\}$, its minimal orbit closure
\[
       Z=\overline{\{\sigma^nz:n\in\Z\}},
\]
and a continuous almost one-to-one factor map
\[
       \pi_Z:Z\longrightarrow\Ztwo
\]
such that
\begin{enumerate}[label=(\roman*)]
  \item $\pi_Z(\sigma y)=\pi_Z(y)+1$ and
  $\pi_Z^{-1}(0)=\{z\}$;
  \item $h_{\rm top}(\sigma\curvearrowright Z)>h$;
  \item every coordinate of $z$ has a period which is a power of two.
\end{enumerate}
\end{lem}

\begin{proof}
Choose $\eta$ with
\[
       \frac{h}{\log2}<\eta<1
\]
and numbers $\varepsilon_k\in(0,1)$ such that
\[
       \prod_{k\geq0}(1-\varepsilon_k)>\eta.
       \tag{5.2}\label{eq:hole-product}
\]
We recursively construct powers of two $p_k$ and periodic partial words
\[
       W_k\in\{0,1,*\}^{\Z/p_k\Z}.
\]
Write $J_k\subset\Z/p_k\Z$ for the set of holes $*$ of $W_k$ and
$\rho_k=|J_k|/p_k$.  Start with $p_0=1$ and $W_0=*$.  Suppose $W_k$ is
constructed.  Choose $m_k=2^{l_k}$ so large that
\[
       \frac{2^{|J_k|}+2k+3}{m_k}<\varepsilon_k,
       \tag{5.3}\label{eq:mk-large}
\]
and put $p_{k+1}=m_kp_k$.  Regard a period of length $p_{k+1}$ as
$m_k$ consecutive copies of $W_k$.

Choose $2^{|J_k|}$ copies and fill their holes in all possible ways,
using each binary filling $J_k\to\{0,1\}$ once.  In addition, for every
integer $n$ with $|n|\leq k+1$ whose residue is still a hole, complete the
copy containing that residue; on copies not already assigned a prescribed
binary filling, fill the remaining holes with zero.  Leave the holes in all
other copies unresolved.  Denote the resulting partial word by $W_{k+1}$.
Previously fixed coordinates are never altered, and at most
$2^{|J_k|}+2k+3$ copies have been completed (that is, contain no holes).  Therefore
\[
 \rho_{k+1}
 \geq \rho_k\left(1-\frac{2^{|J_k|}+2k+3}{m_k}\right)
 >\rho_k(1-\varepsilon_k).
\]
Together with \eqref{eq:hole-product}, this gives
\[
       \rho_k\geq\prod_{j<k}(1-\varepsilon_j)>\eta
       \qquad(k\geq1).
       \tag{5.4}\label{eq:positive-hole-density}
\]

At stage $k+1$, we fill every remaining hole in $[-k-1,k+1]$.
Thus every integer coordinate is fixed after finitely many stages.  Let
$z\in\{0,1\}^{\Z}$ be the resulting complete word.  If the coordinate $n$
is fixed at stage $k$, then its entire residue class modulo $p_k$ has been
fixed with the same value and is never changed later.  Thus $p_k$ is a
period of that coordinate.  Hence $z$ is Toeplitz and every one of its
coordinates has a power-of-two period.

Each finite block of $z$ is fixed at some finite stage.  It then repeats
with that stage's period, so its occurrences have bounded gaps.
Hence $z$ is uniformly recurrent and $Z$ is minimal.  To estimate
entropy, note that stage $k+1$ places all $2^{|J_k|}$ completions of
$W_k$ in the word, and these blocks are never changed afterwards.
Writing $C_Z(n)$ for the language complexity, we obtain
\[
       C_Z(p_k)\geq2^{|J_k|}.
\]
Since $\log C_Z(n)$ is subadditive and $p_k\to\infty$, the entropy limit
can be evaluated along $(p_k)$.  By \eqref{eq:positive-hole-density},
\[
 \begin{split}
 h_{\rm top}(\sigma\curvearrowright Z)
 &=\lim_{k\to\infty}\frac1{p_k}\log C_Z(p_k)\\
 &\geq\liminf_{k\to\infty}\frac{|J_k|}{p_k}\log2
 \geq\eta\log2>h.
 \end{split}
 \tag{5.5}\label{eq:entropy-track}
\]

The positive entropy just proved implies that $z$ is nonperiodic.
Choose an essential period structure $(s_j)$ with each $s_j$ dividing one
of the powers of two $p_k$.  The sequence $(s_j)$ is unbounded: a bounded
essential period structure would make the whole word $z$ periodic.
Consequently $s_j=2^{r_j}$ with $r_j\to\infty$, and its inverse-limit
odometer is $\Ztwo$.  The Toeplitz factor theorem recalled in
Section~\ref{sec:preliminaries} gives an almost one-to-one factor
$\rho:Z\to\Ztwo$ with $\rho^{-1}(\rho(z))=\{z\}$.  Put
\[
       \pi_Z(y)=\rho(y)-\rho(z),
       \qquad \pi_Z^{-1}(0)=\{z\}.
       \tag{5.6}\label{eq:normalized-Toeplitz-factor}
\]
Translation on $\Ztwo$ commutes with addition by $1$, so
$\pi_Z(\sigma y)=\pi_Z(y)+1$.  This proves all three assertions.
\end{proof}

\begin{rem}[Relation with Bratteli--Vershik constructions]
One can also describe the recursion by towers.  The word $W_k$
labels a tower of height $p_k$.  To pass to height $p_{k+1}=m_kp_k$,
cut the tower into $m_k$ subcolumns and stack them.  In the completed
subcolumns we use every binary filling of $J_k$; the other subcolumns
retain their holes.  Their proportion is controlled by \eqref{eq:positive-hole-density},
which yields the entropy bound \eqref{eq:entropy-track}.  See \cite{GJ00} for the Bratteli--Vershik description of Toeplitz flows.
\end{rem}

\begin{pro}[A positive-entropy action satisfying \NekExpHyp]
\label{pro:positive-entropy-action}
For every LR parameter $\omega$ and every $h\in(0,\log2)$ there exist a minimal expansive Toeplitz dihedral subshift
$D_\infty \curvearrowright \Xc$ and a finite fragmentation
\[
  G_\omega\curvearrowright \Xc
\]
which is a faithful almost one-to-one extension of its boundary fragmentation action $G_\omega\curvearrowright \partial T_2$,
with the same LR generic regular orbital graphs and a purely non-Hausdorff
point, such that the translation $t=ab$ has entropy greater than $h$.
Consequently $\Aalt(G_\omega,\Xc)$ is a finitely generated, infinite,
simple, amenable torsion group of exponential word growth. The fragmented
action satisfies \NekExpHyp.
\end{pro}

\begin{proof}
Apply Lemma~\ref{lem:self-contained-positive-Toeplitz} and retain its
notation $(Z,\sigma,z,\pi_Z)$.  Put
\[
  \widehat z(n)=(z(n),z(-n)),\qquad
  \Xc=\overline{\{\sigma^m\widehat z:m\in\Z\}}.
\]
The point $\widehat z$ is Toeplitz: at each position, the two coordinate
periods are powers of two, so one is a multiple of the other and gives
a common period.  Thus $\Xc$ is minimal, and it is expansive because
it is a finite-alphabet subshift.  Projection to the first coordinate
maps $\Xc$ onto $Z$, so
\[
       h_{\rm top}(\sigma\curvearrowright\Xc)
       \geq h_{\rm top}(\sigma\curvearrowright Z)>h.
       \tag{5.7}\label{eq:Xhat-entropy}
\]
Define
\[
  b(y)(n)=\operatorname{swap}(y(-n)),\qquad a=\sigma b.
\]
Then $b\sigma b=\sigma^{-1}$,
\[
       b^2=a^2=1,\qquad ab=\sigma,
\]
and $b\widehat z=\widehat z$.  Thus $a$ and $b$ define a minimal expansive
$\Dinf$-action on $\Xc$.

Let $p_1:\Xc\to Z$ be first-coordinate projection and define
\[
       \pi=\pi_Z\circ p_1:\Xc\longrightarrow\Ztwo.
\]
On the dense orbit of $\widehat z$ one has
\[
       \pi(\sigma^n\widehat z)=n,
\]
so, by continuity,
\[
       \pi(by)=-\pi(y),\qquad
       \pi(ay)=1-\pi(y).
       \tag{5.8}\label{eq:dihedral-odometer-factor}
\]
In particular, $\pi$ is an equivariant factor onto the boundary odometer.
Its fiber over $0$ is the singleton $\{\widehat z\}$.  Indeed, suppose
$\sigma^{n_j}\widehat z\to y$ and $\pi(y)=0$.  Then $n_j\to0$ in $\Ztwo$.  Fix a coordinate.  Its periods in $z$
and in the reversal of $z$ are powers of two, and $n_j$ is eventually
divisible by both.  Thus each fixed coordinate eventually agrees with
that of $\widehat z$, giving
\[
       \sigma^{n_j}\widehat z\longrightarrow\widehat z.
\]
Thus $y=\widehat z$.

Since every symbol occurs infinitely often in $\omega$, the three
fragmentation pieces accumulate at $0$, which is a PNH point of the
boundary fragmentation $G_\omega$.  Equivariance makes every fiber over
the nonregular boundary orbit $\Z$ a singleton.  Lemma~\ref{lem:canonical-fragmentation-lift}
therefore gives a faithful expansive lift with PNH at $\widehat z$.
Explicitly, put
\[
 P_i^\omega=\bigcup_{\omega_k=i}W_k\qquad(i=0,1,2).
\]
Then $b_{ij}$ acts as $b$ over $P_i^\omega\cup P_j^\omega$ and as the
identity elsewhere.  Since $\omega$ is LR,
Proposition~\ref{pro:valuationLR} gives linear repetitivity of the
regular boundary orbital graphs.
Lemma~\ref{lem:regular-graph-extension} identifies every regular lifted
graph with its boundary graph, including over non-singleton regular
fibers.  Hence their language and LR repetitivity function are unchanged.
Together with \eqref{eq:Xhat-entropy}, this proves the dynamical
assertions.  The algebraic conclusions for $\Aalt(G_\omega,\Xc)$ now follow
from Theorem~\ref{thm:positive-entropy-criterion}.

\end{proof}

\begin{rem}[The unfragmented alternating full group]
\label{rem:unfragmented-alternating-full-group}
For this positive-entropy subshift, $\Aalt(\Dinf,\Xc)$ itself has
exponential growth.  Let $\mathcal G_\sigma$ and $\mathcal G_{\Dinf}$ be the
groupoids of the translation and dihedral actions.  Since
$\mathcal G_\sigma$ is an open subgroupoid of $\mathcal G_{\Dinf}$, every
compact-open degree-three multisection of $\mathcal G_\sigma$ is also one of
$\mathcal G_{\Dinf}$, and hence
\[
   \Aalt(\mathcal G_\sigma)\leq\Aalt(\mathcal G_{\Dinf}).
\]
For a minimal subshift, $\Aalt(\mathcal G_\sigma)$ is the commutator
subgroup of the topological full group, and Matui proved that this group has
exponential growth \cite{Matui12}.  Consequently every finite fragmentation
containing $\Aalt(\Dinf,\Xc)$ is exponential as well.  Notice that this
particular conclusion does not require positive entropy; aperiodicity and
minimality of the subshift suffice.

This unfragmented inclusion argument does not by itself imply exponential
growth of $\Aalt(G_\omega,\Xc)$: the fragmented and unfragmented
actions have different germ groupoids at the purely non-Hausdorff point, and
there is generally no inclusion between their alternating full groups.
For the fragmented action, exponential growth follows instead from
Proposition~\ref{pro:virtuality}, which transfers growth from a
germ-preserving containing fragmentation to $\Aalt(G,\Xc)$.
\end{rem}

\Needspace{5\baselineskip}
\subsubsection*{Sparse-marker extensions and nonisomorphic groups}
\label{subsec:continuum-groups}

We now extend the positive-entropy example again, keeping the generic
fragmented graph language unchanged while making the repetitivity of
the alternating full group's orbital graphs uniformly larger than a
prescribed function up to linear rescaling.

In Proposition~\ref{pro:general-expansive-marker}, we construct a
three-valued sparse-marker coding rule over a general expansive
fragmentation.  The resulting extension has three-point fibers over
one chosen regular orbit and singleton fibers over the nonregular
points.  In Theorem~\ref{thm:continuum-simple-torsion}, we make two
choices of clopen neighborhood at each stage.  These give continuum
many symbolic extensions of one fixed base, distinct as subsets of
the same shift.  By Lemma~\ref{lem:countable-presentations}, only
countably many can have alternating full groups of any one abstract
isomorphism type.  Hence the family contains continuum many pairwise
nonisomorphic alternating full groups.

\begin{pro}[Three-state sparse-marker extension]
\label{pro:general-expansive-marker}
Let $\Dinf\curvearrowright X$ and $G=\langle A\cup B\rangle\curvearrowright
X$ satisfy \NekExpHyp.  For every prescribed $F$ there is a minimal expansive
almost one-to-one graph semicocycle extension
\[
 \pi:\Xc\longrightarrow X
\]
whose fibers have exactly three points over one free regular dihedral
orbit and are singletons elsewhere.  The lifted fragmentation again
satisfies \NekExpHyp, with the same regular orbital graphs; every
original PNH point has a singleton fiber whose lift is PNH.  The orbital graph of $\Aalt(G,\Xc)$ at every regular point and the
orbital graph of the containing fragmentation with the same germ have
repetitivity uniformly larger than $F$ up
to linear rescaling.
\end{pro}

\begin{proof}
Replace $F$ by $F^\uparrow$ and let $\widehat F$ be the majorant from
Lemma~\ref{lem:MLD-repetitivity}.  Fix the alphabet $\mathcal B_X$ and
the generating partition $\mathcal Q_X=\{Q_\beta: \beta\in\mathcal B_X\}$
from Convention~\ref{conv:expansive-presentations};
write $p_X:X\to\mathcal B_X$ for its central-symbol map, so that $\chi_X=c_{p_X}$.
Choose a PNH point $\xi$ and a $D_\infty$-free $G$-regular point $\eta$ outside its
orbit.  We construct nested clopen neighborhoods
\[
 X=C_0\supsetneq C_1\supsetneq C_2\supsetneq\cdots,
 \qquad\bigcap_j C_j=\{\eta\},
\]
and strictly increasing integers $R_j\geq j+1$ such that, for every point $x$, membership in
$C_j$ is determined by the centered radius-$R_j$ window of $w_X(x)$.
Put $n_j=2R_j+1$.  Start with any clopen proper neighborhood $C_1$ and
choose $R_0,R_1$.  Inductively, for $j\geq1$, after $C_j$ and $R_j$ are fixed, choose an
integer
\[
 N_{j+1}\geq4\widehat F(n_j)+8n_j+8
\]
and a proper clopen neighborhood $C_{j+1}\subsetneq C_j$ of $\eta$
with diameter at most $2^{-j-1}$ such that
\[
 gC_{j+1}\cap C_{j+1}=\varnothing
 \qquad(1\neq g\in B_{\Dinf}(N_{j+1})).
\]
We can choose such a neighborhood because only finitely many
translates are involved and none fixes $\eta$.  Then choose
$R_{j+1}>R_j$ so that membership in $C_{j+1}$ is determined by the
centered radius-$R_{j+1}$ window of $w_X$.

Define the marker and joint coding rules
\[
 m(x)=j\pmod3\quad\text{on }C_j\setminus C_{j+1},
 \qquad
 \varphi=(p_X,m):X\setminus\{\eta\}\longrightarrow
 \mathcal C:=\mathcal B_X\times\{0,1,2\}.
\]
Its coding is defined on the invariant dense set
$D=X\setminus\Dinf\eta$.  Form the graph semicocycle extension and
its symbolic hull:
\[
 \widehat X_\varphi
 =\overline{\{(x,c_\varphi(x)):x\in D\}}
 \subset X\times\mathcal C^{\Dinf},
 \qquad
 Z_\varphi=\overline{c_\varphi(D)}.
\]
Put $\Xc=\widehat X_\varphi$ and let $\pi:\Xc\to X$ and
$\rho:\Xc\to Z_\varphi$ be the two coordinate projections.

For $x\in D$, the coding is continuous at $x$, since every
point of its orbit avoids the only discontinuity of $m$.  Thus the fiber
over $x$ is a singleton.  Over $\eta$, the $\chi_X$-coordinates and all
marker coordinates indexed by $g\in D_\infty \setminus \{1\}$ have unique limits, since
$g\eta\neq\eta$.  Only the central marker can vary.  Its three values $0,1,2$ all occur
arbitrarily close to $\eta$, because each residue class of annular
indices occurs infinitely often and every annulus meets $D$.  Thus $\pi^{-1}(\eta)$ has exactly
three points, and equivariance gives the same fiber size over
$\Dinf\eta$.

By construction, the singleton graph points above $D$ are dense in
$\Xc$.  A nonempty closed invariant subset must project onto $X$,
since the base action is minimal.  It must therefore contain the unique
lift of every point of $D$, and hence all of $\Xc$.  Thus the extension
is minimal and $\pi$ is almost one-to-one.  The coding $c_\varphi$ has $\chi_X|_D$ as its first coordinate, so its coding name at $x$ recovers $x$ continuously.
Equivalently, $\varphi$ satisfies $(U)$ in
Proposition--Definition~\ref{pro:U}; $\rho$ is a conjugacy with the
finite-alphabet subshift $Z_\varphi$.  The lifted dihedral action is therefore expansive.  By
Lemma~\ref{lem:reflection-invariant-generator}, choose an invariant
generating presentation $w_{\Xc}$.

The non-singleton fibers lie over a free regular orbit.
Lemma~\ref{lem:canonical-fragmentation-lift} gives the unique faithful
lifted fragmentation and preserves every original PNH point in its
singleton fiber.  Lemma~\ref{lem:regular-graph-extension} identifies all
regular lifted orbital graphs with their base orbital graphs, preserving their LR
language and repetitivity function.

For $x\in D$, use the cell coding name
\[
 \zeta_\varphi(x)(k)=\bigl(\varphi(T^kx),\varphi(bT^kx)\bigr),
 \qquad T=ab.
\]
It extends to every $y\in\Xc$ by the formula
$\zeta(y)=J(\rho(y))$, with $J$ as in Convention~\ref{conv:cell-subshift}.
Let $R_\zeta$ denote repetitivity in this cell presentation of the regular
orbital lines, allowing both orientations.
For each $j$, choose
\[
 x_j\in(C_{j+2}\setminus C_{j+3})\cap D
       \cap\operatorname{Free}(\Dinf\curvearrowright X),
 \qquad P_j=\zeta_\varphi(x_j)|_{[-R_j,R_j]}.
\]
Such a point exists because the annulus is nonempty and open and the free
locus with the orbit $\Dinf\eta$ removed is dense.  The pattern $P_j$ has
length $n_j=2R_j+1$.  Its $p_X$-coordinates detect membership in $C_j$,
and its central marker is $j+2\pmod3$.  On the dense singleton-fiber
locus these two tests force the projected center to lie in $C_{j+2}$:
the preceding annuli $C_j\setminus C_{j+1}$ and
$C_{j+1}\setminus C_{j+2}$ carry the other two marker values.
This also holds at boundary fibers: approximate a point in the
$P_j$-cylinder by singleton graph points.  Since the cylinder is clopen,
the approximating points eventually lie in it, and their projections
lie in $C_{j+2}$.  Closedness of $C_{j+2}$ gives the same conclusion
at the limit.  Thus $\zeta_\varphi(x)|_{[-R_j,R_j]}=P_j$ implies $x\in C_{j+2}$.
Distinct occurrence centers of $P_j$ on a regular orbital line are therefore more
than $N_{j+2}$ dihedral edges apart, by the separation imposed on
$C_{j+2}$, and hence more than $N_{j+2}/2$ cells apart in a fixed
orientation.  Minimality and continuity of the clopen symbolic
coordinates ensure that $P_j$ occurs on every regular orbital line,
allowing both orientations.  Therefore
\[
 R_\zeta(n_j)\geq N_{j+2}/2-n_j>\widehat F(n_{j+1}),
\]
where the last inequality follows from
$N_{j+2}\geq4\widehat F(n_{j+1})+8n_{j+1}+8$.
If $n_j\leq n<n_{j+1}$, monotonicity now gives
\[
 R_\zeta(n)\geq R_\zeta(n_j)>\widehat F(n_{j+1})\geq\widehat F(n).
\]
Since $\rho$ is a conjugacy, the $\Z$-coding $\zeta=J\circ\rho$
is generating for $T$ and is therefore MLD with the invariant one-track
coding $w_{\Xc}$.  Forgetting the natural fragmentation coordinate
recovers $w_{\Xc}$ locally from the joint $c_G\vee\chi_{\Xc}$.
Thus Lemma~\ref{lem:MLD-repetitivity} transfers the all-scale lower
bound to the joint, up to linear rescaling.
Finally, Proposition~\ref{pro:expansive-full-group-MLD} gives
\[
 \Gamma(\Aalt(G,\Xc),\cdot;S)
 \overset{\MLD}{\cong}c_G\vee\chi_{\Xc}
 \overset{\MLD}{\cong}\Gamma(H\curvearrowright\Xc,\cdot)
\]
for every finite generating set $S$ and a containing
fragmentation $H$.  Lemma~\ref{lem:MLD-repetitivity} gives uniform
domination up to every fixed linear rescaling, with thresholds allowed
to depend on $S$ or $H$ and on that rescaling.
\end{proof}

\begin{cor}[Positive-entropy sparse-marker extensions]
\label{cor:positive-sparse-marker}
Assume the hypotheses of Proposition~\ref{pro:general-expansive-marker} and,
in addition, that $h_{\rm top}(ab\curvearrowright X)>0$.  For every
increasing $F$, its extension $\pi:\Xc\to X$ can be chosen so that
$\Aalt(G,\Xc)$ is finitely generated, infinite, simple, torsion, amenable and
of exponential word growth, while every finite-generator generic orbital
graph of $\Aalt(G,\Xc)$ has repetitivity uniformly larger than $F$ up
to linear rescaling.
\end{cor}

\begin{proof}
Take the extension of Proposition~\ref{pro:general-expansive-marker}.
Its exceptional fibers lie over a countable union of translation orbits,
so Lemma~\ref{lem:universally-null-extension} gives
\[
 h_{\rm top}(T\curvearrowright\Xc)=h_{\rm top}(T\curvearrowright X)>0.
\]
The lifted action is expansive and retains a PNH point by that proposition.
Theorem~\ref{thm:positive-entropy-criterion} gives all the stated
algebraic properties, and Proposition~\ref{pro:general-expansive-marker}
gives the repetitivity claim.
\end{proof}

\begin{lem}[Countably many symbolic presentations of one alternating full group]
\label{lem:countable-presentations}
Fix a finite alphabet $\mathcal C$ and an abstract finitely generated group
$K$.  Among the infinite minimal dihedral subshifts
$Y\subseteq\mathcal C^{\Dinf}$ admitting a finite fragmentation
$G_Y=\langle A_Y\cup B_Y\rangle$ with $\Aalt(G_Y,Y)\cong K$, there are
at most countably many \emph{literally distinct} subshifts
$Y\subseteq\mathcal C^{\Dinf}$.  This counts the subsets themselves,
not merely their conjugacy classes.
\end{lem}
\begin{proof}
Fix one such presentation $K_*=\Aalt(G_*,Y_*)$.  Local three-cycles show
that an alternating full group of a minimal germ groupoid is locally
minimal: inside any nonempty clopen set $U$, the subgroup supported in $U$
sends each point to every other point of its groupoid orbit in $U$, and
that orbit is dense in $U$.  Rubin's theorem in
\cite[Theorem~3.8]{NekDyn} therefore implements every abstract isomorphism
$\Aalt(G_Y,Y)\to K_*$ by a homeomorphism $f:Y\to Y_*$.  In particular,
$f$ conjugates their germ groupoids, consistently with
\cite[Theorem~3.9]{NekDyn}.

The group $\mathcal F_*:=\mathbf F(K_*\curvearrowright Y_*)$ is countable.
To see this, choose a countable clopen basis of $Y_*$.  By compactness,
every clopen set is a finite union of basis elements, so there are only
countably many finite clopen partitions.  An element of $\mathcal F_*$
is specified by one such partition and an element of the countable group
$K_*$ on each atom.  There are only countably many choices.

By Corollary~\ref{pro:same-germ-groupoid}, $G_Y$ and $\Aalt(G_Y,Y)$ have
the same germ groupoid.  Thus the finite sets
\[
 \widetilde A=fA_Yf^{-1},\qquad \widetilde B=fB_Yf^{-1}
\]
are subsets of $\mathcal F_*$.  They determine the transported original
involutions $\widetilde a=faf^{-1}$ and $\widetilde b=fbf^{-1}$ uniquely:
\[
 \widetilde A y=\{y,\widetilde a y\},\qquad
 \widetilde B y=\{y,\widetilde b y\}.
\]
If the first set is a singleton then $\widetilde a y=y$; otherwise its
unique point different from $y$ is $\widetilde a y$.  The second equality
recovers $\widetilde b$ in the same way.  These involutions need not belong
to $\mathcal F_*$; uniqueness from the finite fragment sets is enough.

Transport the zero-coordinate coding rule $p_Y$ to the clopen coding rule
$\widetilde p=p_Y\circ f^{-1}:Y_*\to\mathcal C$.  There are only
countably many choices of the finite data
$(\widetilde A,\widetilde B,\widetilde p)$.  These data recover $Y$ as the
image
\[
 \left\{\bigl(\widetilde p(\widetilde d\,y)\bigr)_{d\in\Dinf}:
                 y\in Y_*\right\}\subseteq\mathcal C^{\Dinf},
\]
where $\widetilde d$ is the word in $\widetilde a,\widetilde b$
corresponding to $d$.  By equivariance of $f$, these are precisely the coding names of points
of $Y$, so the image is $Y$ itself.  Each choice of the finite data
therefore determines one subset of the full shift.  Since there are
only countably many choices, there are only countably many such
subshifts.
\end{proof}

\begin{rem}[Reconstruction and the fragmentation marking]
\label{rem:rigidity-not-marked-dihedral}
An abstract isomorphism $\Phi:K_Y\to K_*$ is implemented by a
homeomorphism $f$ with $\Phi(g)=fgf^{-1}$.  This conjugates the
alternating-full-group actions, but need not give
$f a_Yf^{-1}=a_*$ or $f b_Yf^{-1}=b_*$ for the chosen dihedral
generators.  Rubin reconstruction does not recover this marking.
In Lemma~\ref{lem:countable-presentations}, we use the fragmentation
identities $A_Yy=\{y,a_Yy\}$ and $B_Yy=\{y,b_Yy\}$ to determine the
\emph{transported} involutions from countably many finite choices.
These identities are not available for an arbitrary action with the
same alternating full group.  Including the transported clopen coding
rule then determines the symbolic presentation as a subset, without
claiming conjugacy of the original marked $\Dinf$-systems.
\end{rem}

\begin{thm}[Continuum many simple amenable torsion groups of exponential growth over one fixed base]
\label{thm:continuum-simple-torsion}
Fix an LR parameter $\omega$ and the positive-entropy fragmentation
$G=G_\omega\curvearrowright X$ supplied by
Proposition~\ref{pro:positive-entropy-action} for this parameter.  For every $F:\N\to\N$
there are $2^{\aleph_0}$ almost one-to-one equivariant dihedral extensions
$\pi_\lambda:X_\lambda\to X$ such that the lifted $G$-actions satisfy
\NekExpHyp. The groups $L_\lambda=\Aalt(G,X_\lambda)$ are pairwise
abstractly nonisomorphic finitely generated infinite simple amenable
torsion groups of exponential word growth.

All the extensions have exactly three-point fibers over one fixed free
regular dihedral orbit and singleton fibers elsewhere.  Their regular
$G$-orbital graphs are the graphs of their factor images in $\partial T_2$;
hence their regular languages, complexity functions and LR repetitivity
functions are the same as those of $G\curvearrowright X$.
Writing $T=ab$ and $T_\lambda=a_\lambda b_\lambda$, the maps
\[
 (\pi_\lambda)_*:\mathcal M_{T_\lambda}(X_\lambda)
 \longrightarrow\mathcal M_T(X)
\]
are affine homeomorphisms.  For every $T_\lambda$-invariant measure $\nu$,
$\pi_\lambda$ is a measure-theoretic isomorphism from
$(X_\lambda,\nu,T_\lambda)$ to $(X,(\pi_\lambda)_*\nu,T)$.
In particular, all these systems have the same positive translation
entropy as $X$.  Every point of the original PNH orbit has a singleton
fiber whose point is PNH upstairs.  The $L_\lambda$-orbital graph at every
regular point and the orbital graph of a containing fragmentation with the
same germ have repetitivity uniformly larger than $F$ up to linear rescaling.
Thus the family has the stated common orbital-graph properties.
\end{thm}
\begin{proof}
Replace $F$ by $F^\uparrow$ and use the majorant $\widehat F$ from
Lemma~\ref{lem:MLD-repetitivity}.  Let $\xi$ be the PNH point in the
singleton fiber over $0$ from
Proposition~\ref{pro:positive-entropy-action}.  Fix a free regular point
$\eta$ outside its dihedral orbit.  We make two choices of clopen
neighborhood at each positive stage of
Proposition~\ref{pro:general-expansive-marker}.  Choose
\[
 X\supsetneq C_1^0\supsetneq C_1^1\supsetneq C_2^0
 \supsetneq C_2^1\supsetneq\cdots,
 \qquad \eta\in C_j^1,\qquad \diam C_j^0\leq2^{-j}.
\]
Put $C_0^0=C_0^1=X$.  Choose strictly increasing radii
$R_j\geq j+1$ such that, for every point $x$, membership in
$C_j^0$ and $C_j^1$ is determined by the centered radius-$R_j$
window of $w_X(x)$, and put $n_j=2R_j+1$.  Once the two choices at
stage $j$ and their common $R_j$ are fixed, choose
\[
 N_{j+1}\geq4\widehat F(n_j)+8n_j+8
\]
and $C_{j+1}^0\subsetneq C_j^1$ so small that
\[
 gC_{j+1}^0\cap C_{j+1}^0=\varnothing
 \qquad(1\neq g\in B_{\Dinf}(N_{j+1})).
\]
Then choose $C_{j+1}^1\subsetneq C_{j+1}^0$ and their common determining
radius $R_{j+1}$.

For $\varepsilon=(\varepsilon_j)_{j\geq1}\in\{0,1\}^{\N_{\geq1}}$,
put $C_0^\varepsilon=X$, $C_j^\varepsilon=C_j^{\varepsilon_j}$ for
$j\geq1$, and
\[
 m_\varepsilon(x)=j\pmod3
 \quad\bigl(x\in C_j^\varepsilon\setminus C_{j+1}^\varepsilon\bigr).
\]
With $D=X\setminus\Dinf\eta$ and the fixed alphabet
$\mathcal C=\mathcal B_X\times\{0,1,2\}$, set
\[
 X_\varepsilon=
 \overline{\{c_{(p_X,m_\varepsilon)}(x):x\in D\}}
 \subseteq\mathcal C^{\Dinf},
\]
and denote by $\pi_\varepsilon : X_\varepsilon \to X$ the almost one-to-one factor. The first coordinate is the generating $\chi_X$-name, so it recovers
$x$ continuously and identifies $X_\varepsilon$ with the corresponding
graph semicocycle extension.  Every branch satisfies the choices in
Proposition~\ref{pro:general-expansive-marker}, with the same sequence
$(n_j)$ and separation bounds.  Thus each $X_\varepsilon$ is minimal
and expansive, has exactly three-point fibers over $\Dinf\eta$ and
singleton fibers elsewhere, and admits the faithful lifted fragmentation.
Every point of the protected PNH orbit $\Dinf\xi$ has a singleton
fiber whose lift is PNH.  The bounds in that proposition give
the required eventual uniform repetitivity for every finite generating
set, and for containing fragmentations with the same germ.

The split orbit $D_\infty \eta$ is regular, so
Lemma~\ref{lem:regular-graph-extension} identifies all regular lifted
graphs with their base graphs.  The non-singleton fibers of $\pi_\varepsilon$
lie over $\Dinf\eta$, a countable union of aperiodic translation orbits
in the base.  Lemma~\ref{lem:universally-null-extension}
therefore makes the factor $X_\varepsilon\to X$ induce an affine
homeomorphism of translation-invariant measure spaces and an isomorphism
of each corresponding measure-preserving system.  In particular,
\[
 h_{\rm top}(T\curvearrowright X_\varepsilon)
 =h_{\rm top}(T\curvearrowright X)>0.
\]
Corollary~\ref{cor:positive-sparse-marker} gives the stated algebraic
properties, while the regular-graph preservation shows that each
lifted fragmentation satisfies \NekExpHyp.

We now check that different branches give distinct subsets of
$\mathcal C^{\Dinf}$.  Take $\varepsilon\neq\varepsilon'$ and choose
$j$ with $\varepsilon_j\neq\varepsilon'_j$.  Choose a point
\[
 x\in(C_j^0\setminus C_j^1)\cap D.
\]
The nesting makes $m_\varepsilon(x)$ and $m_{\varepsilon'}(x)$ equal to
$j$ and $j-1$ modulo three, in one order or the other, so they are different.
For each branch $\delta$, let
\[
 P_\delta:X_\delta\longrightarrow\chi_X(X)
\]
be first-coordinate projection.  Since $\chi_X$ is generating, it
identifies $X$ with $\chi_X(X)$, and $P_\delta=\chi_X\circ\pi_\delta$.
Continuity of each branch's coding rule along the orbit of $x\in D$ gives
\[
 \begin{aligned}
 P_\varepsilon^{-1}(\chi_X(x))
   &=\{c_{(p_X,m_\varepsilon)}(x)\},\\
 P_{\varepsilon'}^{-1}(\chi_X(x))
   &=\{c_{(p_X,m_{\varepsilon'})}(x)\}.
 \end{aligned}
\]
The two names differ at the central marker.  If
$X_\varepsilon=X_{\varepsilon'}$ as subsets of $\mathcal C^{\Dinf}$,
restriction of the first-coordinate projection would be the same map
on both.  Its fiber over $\chi_X(x)$ could not be these two different
singletons.  Thus the subshifts are distinct.

There are thus $2^{\aleph_0}$ literally distinct hulls $X_\varepsilon$.
By Lemma~\ref{lem:countable-presentations}, at most countably many can
have one abstract alternating-full-group isomorphism type.  A partition
of a set of cardinality $2^{\aleph_0}$ into countable classes has
$2^{\aleph_0}$ classes: for every infinite cardinal $\kappa$,
$\kappa\cdot\aleph_0=\kappa$.  Selecting one branch from each class
gives the required family $\{X_\lambda\}$ within the sparse-marker
constructions.

Finally, let $\pi_0:X\to\partial T_2$ be the fixed factor from
Proposition~\ref{pro:positive-entropy-action}.  Its fiber over $0$ is
$\{\xi\}$, and $\pi_\lambda^{-1}(\xi)$ is also a singleton.
Thus $\pi_0\pi_\lambda$ has a singleton fiber over $0$.
We also need minimality for the translation action.  Let a nonempty
closed $T_\lambda$-invariant subset of $X_\lambda$ be given.  Its
image is all of $X$, since $T\curvearrowright X$ is minimal, so it
contains every singleton fiber of $\pi_\lambda$.  These fibers are
dense, and the subset is therefore all of $X_\lambda$.  This proves
minimality of $T_\lambda\curvearrowright X_\lambda$.
It is also expansive, since its finite-alphabet $\Dinf$-subshift
presentation is expansive, and expansivity passes to the finite-index
translation subgroup $\langle T\rangle$.  Equivariance and minimality now make
$\pi_0\pi_\lambda$ an almost one-to-one odometer factor of this minimal
subshift.  By the Toeplitz factor characterization recalled in
Section~\ref{sec:preliminaries}, the $X_\lambda$ are Toeplitz systems.
These composite maps intertwine the lifted named generators with the
generalized Grigorchuk boundary action
$G_\omega\curvearrowright\partial T_2$, and the regular graph
identifications compose as well.
This proves Theorem~\ref{thm:positive-entropy} in its boundary-extension form,
including the uniform repetitivity assertion.
\end{proof}

\begin{rem}[An alternative proof of the continuum conclusion]
\label{rem:alternative-continuum}
For uncountability alone in Theorem~\ref{thm:continuum-simple-torsion},
Lemma~\ref{lem:countable-presentations} and the arbitrarily large
repetitivity of Corollary~\ref{cor:positive-sparse-marker} suffice;
we do not need the binary branching used in the proof.

To obtain continuum many groups, one may instead use a
descriptive-set-theoretic argument: the set of actions satisfying
Theorem~\ref{thm:continuum-simple-torsion} is Borel in a suitable
standard coding.

\end{rem}

\begin{rem}[A Toeplitz alternative for large repetitivity]
\label{rem:toeplitz-bad-repetitivity}
One may also obtain alternating full groups with the algebraic properties
stated in Theorem~\ref{thm:positive-entropy} and arbitrarily bad
regular-orbital repetitivity by choosing the positive-entropy dyadic
Toeplitz word in Proposition~\ref{pro:positive-entropy-action} to have
repetitivity larger than a prescribed function.  Its Toeplitz coordinate
is locally recoverable from the generating coding, and hence from the
regular orbital graphs of $\Aalt(G_\omega,\Xc)$ by
Proposition~\ref{pro:expansive-full-group-MLD}; the lower bounds transfer
up to linear rescaling.  This approach is specific to the Toeplitz
boundary construction, whereas Proposition~\ref{pro:general-expansive-marker}
provides a general extension mechanism for any base fragmentation
satisfying \NekExpHyp, without requiring a special Toeplitz model.
\end{rem}

Nekrashevych constructs continuum many pairwise different growth types
of simple finitely generated groups in \cite[Theorem~4.4]{NekSub25}.
In Theorem~\ref{thm:continuum-simple-torsion}, all groups have exponential
growth.  What varies is their abstract isomorphism type, while the
regular fragmented graph language stays fixed.

\subsection{Zero-entropy examples with non-LR joint codings}\label{sec:construction}
We now prove the full version of Theorem~\ref{thm:main-transfer-strong}.
We use the valuation notation and LR criterion from
Section~\ref{sec:preliminaries}.  The dyadic construction is similar to
the positive-entropy example, but now both coordinate words are LR and
their joint has zero entropy.

\begin{thm}[The refined two-LR-coordinate construction]
\label{thm:refined-full}
Let $F:\N\to\N$ be arbitrary and let $\omega$ be an LR parameter.  There
exist a binary level sequence $q=(q_k)$, its simple Toeplitz subshift
$X$, a minimal expansive
$\Dinf$-action, an almost one-to-one factor
\[
       \pi:X\longrightarrow\partial T_2\cong\Ztwo,
\]
and a faithful lifted fragmentation
$G_\omega\curvearrowright X$ such that, for every regular
$z\in\Ztwo\setminus\Z$ with unique lift $\bar z$:
\begin{enumerate}[label=(\roman*)]
  \item The non-singleton fibers are exactly those over
  $\Dinf\cdot0^\infty\cong\Z$, and each has two points.  The maximal
  equicontinuous factor of the dihedral action is the boundary odometer.

  \item The invariant one-track coding $w_X(\bar z)\overset{\MLD}{\cong}u_q^z$ is LR.  For every finite
  generating set of $\Aalt(\Dinf,X)$ its orbital graph is MLD with
  $u_q^z$, and there is a finite fragmentation
  $H_0\supseteq\Aalt(\Dinf,X)$ with
  \[
     \Gamma(H_0\curvearrowright X,\bar z)
       \overset{\MLD}{\cong}u_q^z.
  \]

  \item The action $G_\omega\curvearrowright X$ is faithful, minimal,
  expansive and LR, and
  \[
     \Gamma(G_\omega\curvearrowright X,\bar z)
       \overset{\MLD}{\cong}u_\omega^z.
  \]
  Every point above the exceptional base orbit is purely non-Hausdorff.  The natural
  fragmentation coding rule satisfies $(W)$, which here means precisely that
  \[
      x\longmapsto
      \Gamma(G_\omega\curvearrowright X,x;A\cup B)
  \]
  is injective on the Nekrashevych-regular locus: taking the rooted labeled
  orbital graph determines the regular point.

  \item Condition $(U)$ fails.  The joint $u_\omega^z\vee u_q^z$ has
  at most $12n$ distinct length-$n$ patterns for every $n\geq1$, but its
  repetitivity is uniformly larger than $F$ up to linear rescaling.

  \item There is a finite fragmentation $H$ of the same dihedral action,
  containing both $G_\omega$ and $\Aalt(G_\omega,X)$ and having the same
  germ as $G_\omega \curvearrowright X$, such that
  \[
     \Gamma(H\curvearrowright X,\bar z)
       \overset{\MLD}{\cong}u_\omega^z\vee u_q^z.
  \]
  Moreover, for every finite generating set $S$ of
  $\Aalt(G_\omega,X)$,
  \[
     \Gamma(\Aalt(G_\omega,X)\curvearrowright X,\bar z;S)
       \overset{\MLD}{\cong}u_\omega^z\vee u_q^z,
  \]
  and this graph has repetitivity uniformly larger than $F$ up to linear
  rescaling.  Its number $C_S(R)$ of rooted radius-$R$ ball types is
  $O(R)$; the containing fragmentation $H$ also has linear
  rooted-ball complexity.  The implied constants may depend on the
  generating sets.  Every finite fragmentation $H'$ of the same dihedral action
  with $\Aalt(G_\omega,X)\leq H'$ has the same large repetivity property.  If, in addition,
  $\operatorname{Germ}(H'\curvearrowright X)$ is the original fragmented
  germ groupoid, then $H'$ contains $\Aalt(G_\omega,X)$ with finite
  index and therefore has the same growth type.
\end{enumerate}
\end{thm}

Fix an LR parameter $\omega$, replace $F$ by $F^\uparrow$, and let
$\widehat F$ be as in Lemma~\ref{lem:MLD-repetitivity}.  Choose $M$
such that every length-$M$ block of $\omega$ contains $0,1,2$.
Choose levels
\[
 0=\ell_0<\ell_1<\ell_2<\cdots,\qquad
 \ell_{j+1}-\ell_j\geq2M,
 \qquad n_j=2^{\ell_j+1}+1,
\]
with $\ell_1\geq2M$ and, successively,
\[
 2^{\ell_{j+2}}>\widehat F(n_{j+1})+2n_{j+1}
 \qquad(j\geq0).
 \tag{5.9}\label{eq:phase-bound}
\]
The right-hand side is already fixed when $\ell_{j+2}$ is chosen.
Write $I_j=[\ell_j,\ell_{j+1})\cap\N$.

\begin{lem}[A cyclic binary level sequence]
\label{lem:q-construction}
There is $q=(q_k)\in\{0,1\}^{\N}$ with:
\begin{enumerate}[label=(Q\arabic*)]
 \item for $k\in I_j$, $q_k=1$ exactly when $\omega_k=j\pmod3$;
 \item both symbols of $\{0,1\}$ occur in every length-$2M$ block of $q$;
 \item every pair $(i,\varepsilon)\in\{0,1,2\}\times\{0,1\}$ occurs
 infinitely often as $(\omega_k,q_k)$.
\end{enumerate}
\end{lem}
\begin{proof}
Define $q$ by (Q1).  Consider an interval of $2M$ consecutive levels.
It crosses at most one boundary between the $I_j$, since these
boundaries are at least $2M$ apart.  Hence it contains $M$ consecutive
levels within one $I_j$.  All three symbols of $\omega$ occur there:
the symbol $j\pmod3$ gives $q=1$, and either of the other two gives
$q=0$.  This proves (Q2).
For each $i$, the pair $(i,1)$ occurs in every interval of type
$j\equiv i\pmod3$, whereas $(i,0)$ occurs in both other types of
interval.  Each type occurs infinitely often, proving (Q3).
\end{proof}

\begin{rem}[The dyadic picture]\label{rem:dyadic-picture}
For every level sequence $q$, the positive part of the valuation word begins
\[
\begin{array}{c|cccccccccccccccc}
 n&1&2&3&4&5&6&7&8&9&10&11&12&13&14&15&16\\ \hline
 u_q^0(n)&q_0&q_1&q_0&q_2&q_0&q_1&q_0&q_3&q_0&q_1&q_0&q_2&q_0&q_1&q_0&q_4.
\end{array}
\]
The labels at the low levels repeat periodically.  In an interval of
length $2^r$, at most one position has level at least $r$.  The rule on
each $I_j$ keeps $q$ recurrent, but a suitable letter of the joint is
absent throughout two consecutive level intervals.  We use this absence
to obtain the repetitivity lower bounds below.
\end{rem}

Choose $\theta\in\Ztwo\setminus\Z$ and put
\[
 x_\theta(n)=q_{\nu_2(\theta+n)},\qquad
 X=\overline{\{\sigma^n x_\theta:n\in\Z\}}.
\]

The period and phase terminology is as in
Section~\ref{sec:preliminaries}.

\begin{pro}[The simple Toeplitz extension]
\label{pro:Toeplitz}
The subshift $X$ is minimal and LR.  There is an almost one-to-one factor
\[
 \pi:X\longrightarrow\Ztwo
\]
with singleton fibers off $\Z$, two-point fibers over $\Z$, and maximal
equicontinuous factor equal to the dyadic odometer.  Moreover,
$p_X(\sigma^n\bar z)=u_q^z(n)$.
\end{pro}

\begin{proof}
If $k=\nu_2(\theta+n)$, coordinate $n$ has period $2^{k+1}$; thus $x_\theta$
is Toeplitz.  Proposition~\ref{pro:valuationLR} and (Q2) give LR and
minimality.  To rule out a nonzero period $m$, put $r=\nu_2(m)$.  Choose $k>r$
with $q_k\ne q_r$, and then $n$ with $\nu_2(\theta+n)=k$.
Adding $m$ changes this valuation to $r$, so
$x_\theta(n+m)\ne x_\theta(n)$.
The essential dyadic period structure is therefore unbounded, and the
Toeplitz factor theorem gives the dyadic odometer, normalized by
$\pi(x_\theta)=\theta$.
For $s\notin\Z$, we have $s+n\ne0$ at every coordinate $n$.
The valuation is locally constant there, so that coordinate must be
$q_{\nu_2(s+n)}$.  Thus the fiber over $s$ is a singleton.
At phase $0$, the same argument determines every noncentral coordinate.
The central coordinate has the two cluster values $0,1$, both attained
by (Q2).  The other fibers over $\Z$ also have two points, by
equivariance.  For the shift action, the coding name of a subshift point under the zero-coordinate rule is the point itself.
\end{proof}

Define the dihedral action by
\[
 b(x)(n)=x(-n),\qquad a=\sigma b,
\]
so that $ab=\sigma$.

The points $\xi_0,\xi_1$ above phase $0$ satisfy
\[
 \xi_\varepsilon(0)=\varepsilon,\qquad
 \xi_\varepsilon(n)=q_{\nu_2(n)}\quad(n\ne0),
\]
so both are fixed by $b$.  On the dyadic base $a$ has no fixed point
and $\Fix(b)=\{0\}$; hence
\eqref{eq:fixed-reflection-fibers} holds.
Lemma~\ref{lem:canonical-fragmentation-lift} gives the unique faithful
lift of the boundary fragmentation $G_\omega$.  Its pieces are
\[
 P_i^\circ=\pi^{-1}\Bigl(\bigcup_{\omega_k=i}W_k\Bigr)
           \cap\pi^{-1}(\Ztwo\setminus\Z),\qquad
 P_i=\overline{P_i^\circ}.
\]
The fragment $b_{ij}$ acts as $b$ on $P_i\cup P_j$ and as the identity
on $P_k$, where $\{i,j,k\}=\{0,1,2\}$; these prescriptions agree at
the reflection-fixed lifts.

\begin{pro}[Lifted action, regular graphs and $(W)$]
\label{pro:lifted-all}
The lifted action is faithful, minimal and expansive.  Its regular orbital
graph at $\bar z$ is canonically the original boundary graph and hence MLD
with $u_\omega^z$.  The two points above $0$ and their orbits are purely
non-Hausdorff.  The natural fragmentation coding rule satisfies $(W)$, and on the
regular locus this means that
\[
 x\longmapsto\Gamma(G_\omega\curvearrowright X,x;A\cup B)
\]
is injective.
\end{pro}

\begin{proof}
The lifting and graph-comparison lemmas
\ref{lem:canonical-fragmentation-lift} and
\ref{lem:compatible-free-graphs} give faithfulness, minimality,
expansivity and the stated LR boundary graphs.  Condition (Q3) makes every $P_i^\circ$ accumulate on
both $\xi_0$ and $\xi_1$.  Each nontrivial Klein-four fragment therefore
has a nontrivial germ at each point and has a complementary open fixed
piece accumulating there.  These points and their orbits are PNH.
All other base fibers are singleton and free, so there are no further
singular points upstairs.

If $z\ne z'$ are regular phases and $r=\nu_2(z-z')$, choose $k>r$
with $\omega_k\ne\omega_r$ and an integer $n$ with
$\nu_2(z+n)=k$.  Then $\nu_2(z'+n)=r$, so their natural fragmentation
coding names differ.  The regular fibers are singleton and all singular points
are non-Hausdorff, so $D_p$ is precisely the regular locus here.
This proves $(W)$ and the stated graph injectivity.
\end{proof}

Put $d_k=(\omega_k,q_k)$ and $u_d^z=u_\omega^z\vee u_q^z$.
For each $j$, the letter
\[
 \beta_j=((j+2)\bmod3,1)\in \{0,1,2\}\times \{0,1\}
\]
is absent at every level $\ell_j\leq k<\ell_{j+2}$.  Let $k_j$ be
its first occurrence at or after $\ell_j$.  By (Q1) and the bound for
$\omega$,
\[
 \ell_{j+2}\leq k_j<\ell_{j+2}+M.
\]

\begin{pro}[Linear complexity and uniformly large repetitivity of the joint]
\label{pro:joint-nonLR}
Let $C_{\rm joint}(n)$ count the length-$n$ words of
$u_d^z=u_\omega^z\vee u_q^z$.  Uniformly in
$z\in\Ztwo\setminus\Z$,
\[
 C_{\rm joint}(n)\leq12n\qquad(n\geq1).
\]
For $n_j=2^{\ell_j+1}+1$ and every $z\in\Ztwo\setminus\Z$,
\[
 R_{u_d^z}(n_j)\geq2^{k_j}-n_j
 \geq2^{\ell_{j+2}}-n_j>\widehat F(n_{j+1}).
\]
Consequently $R_{u_d^z}(n)>\widehat F(n)$ for every sufficiently large $n$,
and the joint has repetitivity uniformly larger than $F$ up to linear
rescaling.  Both coordinate words are LR, whereas their joint is not LR.
\end{pro}
\begin{proof}
Fix $n\geq1$ and choose $r$ with $n<2^r\leq2n$.  A length-$n$
interval contains at most one position of valuation at least $r$.
Its starting residue modulo $2^r$ determines the location of that
position and all other labels from $d_0,\ldots,d_{r-1}$.
There are at most six possibilities for the remaining label, since
$d_k\in\{0,1,2\}\times\{0,1\}$.  Thus
\[
 C_{\rm joint}(n)\leq6\cdot2^r\leq12n.
\]
The estimate uses only that, once the starting residue is fixed,
at most one label remains undetermined.  See also \cite{Sell20} for
the combinatorics of simple Toeplitz subshifts.

To bound repetitivity from below, take a block of radius $2^{\ell_j}$
centered at a position of valuation $k_j$ in $u_d^z$.  We show that it
cannot recur at a center of valuation $s<k_j$.
If $s\geq\ell_j$, the central symbol is not $\beta_j$, by the choice
of $k_j$.  If $s<\ell_j$, agreement of the central symbols would give
$d_s=\beta_j$.  Write $v=z+m$ for this proposed center and choose an
integer $t$ with
\[
 |t|\leq2^{\ell_j},\qquad
 v+t\equiv2^{\ell_j}\pmod{2^{\ell_j+1}}.
\]
Then $\nu_2(v+t)=\ell_j$ and $\nu_2(t)=s$.  At the original high
center, the same offset has valuation $s$, so its symbol is
$d_s=\beta_j$; at the proposed center it is
$d_{\ell_j}\neq\beta_j$.  Thus these blocks cannot agree.

We have shown that every occurrence center has valuation at least
$k_j$.  All such centers lie in one residue class modulo $2^{k_j}$,
so distinct occurrences are separated by at least $2^{k_j}$.  Subtracting
the pattern length $n_j$ gives
$R_{u_d^z}(n_j)\geq2^{k_j}-n_j$.
Equation~\eqref{eq:phase-bound}, together with Lemma~\ref{lem:MLD-repetitivity}, gives the stated domination under every
fixed linear rescaling.
The individual words are LR by Proposition~\ref{pro:valuationLR} and
(Q2), while $R_{u_d^z}(n)>\widehat F(n)\geq n^2$ rules out LR of the
joint.
\end{proof}

\begin{cor}[Failure of $(U)$]
\label{cor:not-U}
The natural fragmentation coding rule satisfies $(W)$ but not $(U)$.
\end{cor}

\begin{proof}
If $(U)$ held, $\chi_{X}$ would be locally derivable from the LR
fragmentation coding.  Their joint would then be LR by
Lemma~\ref{lem:join-MLD}, contradicting
Proposition~\ref{pro:joint-nonLR}.
\end{proof}

\phantomsection\label{proof:main-Ib}
\begin{proof}[Proof of Theorem~\ref{thm:refined-full} and Theorem~\ref{thm:main-transfer-strong}]
Lemma~\ref{lem:q-construction} and Propositions~\ref{pro:Toeplitz}
and~\ref{pro:lifted-all} give (i) and (iii).
Proposition~\ref{pro:joint-nonLR} and Corollary~\ref{cor:not-U} give
(iv), including the eventual all-scale lower bound.

For (ii), apply Proposition~\ref{pro:expansive-full-group-MLD} to the
unfragmented action, whose natural crossing profiles are constant.
By the chosen invariant coding, the single track $w_X(\bar z)$ in
\eqref{eq:one-track-joint} is MLD with $u_q^z$.  The proposition
therefore identifies this LR word with every finite-generator graph of
$\Aalt(\Dinf,X)$ and with a containing fragmentation $H_0$.

For (v), apply the same proposition to $G_\omega\curvearrowright X$.
At $\bar z$, its natural coding name is MLD with $u_\omega^z$ and its clopen coding name with
$u_q^z$, so Lemma~\ref{lem:join-MLD} gives
\[
 \Gamma(\Aalt(G_\omega,X),\bar z;S)
 \overset{\MLD}{\cong}u_\omega^z\vee u_q^z
 \overset{\MLD}{\cong}\Gamma(H\curvearrowright X,\bar z)
\]
for every finite generating set $S$ and a containing germ-preserving fragmentation $H$.
A rooted radius-$R$ ball in either graph is determined by a joint window
of length at most $cR+c$, for a fixed integer $c$ depending on the
presentation.  Proposition~\ref{pro:joint-nonLR} therefore gives
\[
 C_S(R)\leq c\,C_{\rm joint}(cR+c)=O(R),
\]
where the constant also absorbs the two vertices per cell and the choice
of orientation.  The same argument gives linear rooted-ball complexity
for $H$.
For any other finite containing fragmentation $H'$, each generator
in $S$ is a word of bounded length in the fragment generators of $H'$.
Its orbital graph therefore also determines the joint locally, and
the line metrics are bi-Lipschitz comparable.  By
Proposition~\ref{pro:joint-nonLR} and Lemma~\ref{lem:MLD-repetitivity},
the same repetitivity lower bound, up to linear rescaling, holds for
all these graphs.  If $H'$ preserves the germ groupoid, then
$H'\leq\mathbf F(G_\omega)$ and Proposition~\ref{pro:virtuality} gives
the finite-index and growth assertions.
\end{proof}

\begin{rem}[Bad recurrence alone is not a growth lower bound]
Choose a generalized Grigorchuk parameter which is parity-periodic except
for a third symbol at levels with arbitrarily large gaps.  Its valuation
orbital graph can have arbitrarily bad repetitivity, while the corresponding
generalized Grigorchuk group still has intermediate growth
\cite{Gri85,Sam20}.
\end{rem}

\section*{Acknowledgments of AI assistance}
Generative AI tools, primarily ChatGPT (OpenAI) and Claude (Anthropic),
were used extensively in the development and preparation of this
manuscript, beyond linguistic editing. The principal mathematical ideas
and constructions of this work---including, but not limited to, the use
of positive-entropy almost one-to-one extensions and semicocycle coding
to control orbital graphs, the definitions and recovery properties of
conditions $(W)$ and $(U)$ (distinguishing global identification of points
from determination by uniform finite windows), and the MLD comparison
between the regular orbital graphs of alternating full groups and the joint
fragmentation and ambient dynamical codings---originated with the author.
Starting from
the author's original draft and detailed mathematical instructions,
AI tools assisted in producing an initial complete version of the
manuscript, which was subsequently revised through extensive iterative
discussions.

AI assistance also contributed to the mathematical development of the
work. In particular, it suggested reformulating the construction in
Section~5 using Toeplitz systems rather than the original
semicocycle-based presentation, helped identify a relevant result for
Proposition~3.20, and contributed to Remarks~5.7 and~5.15.
Proposition~4.12 also emerged from discussions with AI.

The author has critically examined and revised the resulting arguments
and takes full responsibility for all mathematical statements, proofs,
and any remaining errors.

\Needspace{20\baselineskip}

\end{document}